\documentclass{Article2023}

\usepackage{xr}
\newcommand{\Aca}{\mathcal{A}}

\newcommand{\Oca}{\mathcal{O}}
\newcommand{\Mca}{\mathcal{M}}

\newcommand{\Dbb}{\mathbb{D}}
\newcommand{\Ebb}{\mathbb{E}}
\newcommand{\Sbb}{\mathbb{S}}

\newcommand{\Bfr}{\mathfrak{b}}
\newcommand{\Cfr}{\mathfrak{c}}
\newcommand{\Dfr}{\mathfrak{d}}

\newcommand{\Tfr}{\mathfrak{t}}
\newcommand{\Sfr}{\mathfrak{s}}
\newcommand{\Pfr}{\mathfrak{p}}
\newcommand{\Qfr}{\mathfrak{q}}
\newcommand{\Rfr}{\mathfrak{r}}
\newcommand{\Ufr}{\mathfrak{u}}

\newcommand{\Att}{\mathtt{a}}
\newcommand{\Btt}{\mathtt{b}}
\newcommand{\Ctt}{\mathtt{c}}
\newcommand{\Dtt}{\mathtt{d}}
\newcommand{\Ett}{\mathtt{e}}

\newcommand{\As}{\mathsf{As}}
\newcommand{\BNC}{\mathsf{BNC}}

\newcommand{\Cli}{\mathsf{C}}

\newcommand{\RatFct}{\mathsf{RatFct}}

\newcommand{\Dendr}{\mathsf{Dendr}}
\newcommand{\NCP}{\mathsf{NCT}}
\newcommand{\LOp}{\mathsf{L}}

\newcommand{\Cro}{\mathsf{Cro}}

\newcommand{\NC}{\mathrm{NC}}

\newcommand{\Motzkin}{\mathrm{Mot}}

\newcommand{\FF}{\mathcal{F}\mathcal{F}}
\newcommand{\Motz}{\mathrm{Motz}}

\newcommand{\Unit}{\mathds{1}}
\newcommand{\Hilbert}{\mathcal{H}}
\newcommand{\Op}{\star}
\newcommand{\OpAssoc}{\odot}
\newcommand{\Leaf}{\perp}
\newcommand{\Nar}{\mathrm{nar}}
\DeclareMathOperator{\Rew}{\to}
\newcommand{\RewTrans}{\overset{*}{\Rew}}
\newcommand{\RewTransSym}{\overset{*}{\leftrightarrow}}
\newcommand{\CRew}{\Rightarrow}
\newcommand{\CRewTrans}{\overset{*}{\CRew}}
\newcommand{\CRewTransSym}{\overset{*}{\Leftrightarrow}}
\newcommand{\Free}{\mathrm{Free}}
\newcommand{\Gen}{\mathfrak{G}}
\newcommand{\Rel}{\mathfrak{R}}

\newcommand{\Eval}{\mathrm{ev}}
\newcommand{\Corolla}{\mathrm{c}}
\newcommand{\RelEq}{\equiv}
\newcommand{\Bubbles}{\mathcal{B}}
\newcommand{\Triangles}{\mathcal{T}}

\newcommand{\Cliques}{\mathcal{C}}
\newcommand{\BubbleTree}{\mathrm{bt}}
\newcommand{\Returned}{\mathrm{ref}}

\newcommand{\Alg}{\Aca}
\newcommand{\Id}{\mathrm{Id}}
\newcommand{\SeriesBubbles}{\mathrm{B}}
\newcommand{\SeriesElements}{\mathrm{F}}

\newcommand{\Border}{\mathrm{bor}}

\newcommand{\Frac}{\mathrm{F}}
\newcommand{\GDendr}{\prec}
\newcommand{\DDendr}{\succ}

\newcommand{\UnitClique}{
\begin{tikzpicture}[scale=.6,Centering]
    \node[CliquePoint](1)at(0,0){};
    \node[CliquePoint](2)at(.75,0){};
    \draw[CliqueEmptyEdge](1)edge[]node[CliqueLabel]{}(2);
\end{tikzpicture}}

\newcommand{\Triangle}[3]{
\begin{tikzpicture}[scale=.42,Centering]
    \node[CliquePoint](1)at(0,1){};
    \node[CliquePoint](2)at(0.87,-0.5){};
    \node[CliquePoint](3)at(-0.87,-0.5){};
    \draw[CliqueEdge](1)edge[]node[CliqueLabel]{$#3$}(2);
    \draw[CliqueEdge](1)edge[]node[CliqueLabel]{$#2$}(3);
    \draw[CliqueEdge](2)edge[]node[CliqueLabel]{$#1$}(3);
\end{tikzpicture}}

\newcommand{\TriangleEXX}[3]{
\begin{tikzpicture}[scale=.42,Centering]
    \node[CliquePoint](1)at(0,1){};
    \node[CliquePoint](2)at(0.87,-0.5){};
    \node[CliquePoint](3)at(-0.87,-0.5){};
    \draw[CliqueEdge](1)edge[]node[CliqueLabel]{$#3$}(2);
    \draw[CliqueEdge](1)edge[]node[CliqueLabel]{$#2$}(3);
    \draw[CliqueEmptyEdge](2)edge[]node[CliqueLabel]{}(3);
\end{tikzpicture}}

\newcommand{\TriangleXEX}[3]{
\begin{tikzpicture}[scale=.42,Centering]
    \node[CliquePoint](1)at(0,1){};
    \node[CliquePoint](2)at(0.87,-0.5){};
    \node[CliquePoint](3)at(-0.87,-0.5){};
    \draw[CliqueEdge](1)edge[]node[CliqueLabel]{$#3$}(2);
    \draw[CliqueEmptyEdge](1)edge[]node[CliqueLabel]{}(3);
    \draw[CliqueEdge](2)edge[]node[CliqueLabel]{$#1$}(3);
\end{tikzpicture}}

\newcommand{\TriangleXXE}[3]{
\begin{tikzpicture}[scale=.42,Centering]
    \node[CliquePoint](1)at(0,1){};
    \node[CliquePoint](2)at(0.87,-0.5){};
    \node[CliquePoint](3)at(-0.87,-0.5){};
    \draw[CliqueEmptyEdge](1)edge[]node[CliqueLabel]{}(2);
    \draw[CliqueEdge](1)edge[]node[CliqueLabel]{$#2$}(3);
    \draw[CliqueEdge](2)edge[]node[CliqueLabel]{$#1$}(3);
\end{tikzpicture}}

\newcommand{\TriangleXEE}[3]{
\begin{tikzpicture}[scale=.42,Centering]
    \node[CliquePoint](1)at(0,1){};
    \node[CliquePoint](2)at(0.87,-0.5){};
    \node[CliquePoint](3)at(-0.87,-0.5){};
    \draw[CliqueEmptyEdge](1)edge[]node[CliqueLabel]{}(2);
    \draw[CliqueEmptyEdge](1)edge[]node[CliqueLabel]{}(3);
    \draw[CliqueEdge](2)edge[]node[CliqueLabel]{$#1$}(3);
\end{tikzpicture}}

\newcommand{\TriangleEEX}[3]{
\begin{tikzpicture}[scale=.42,Centering]
    \node[CliquePoint](1)at(0,1){};
    \node[CliquePoint](2)at(0.87,-0.5){};
    \node[CliquePoint](3)at(-0.87,-0.5){};
    \draw[CliqueEdge](1)edge[]node[CliqueLabel]{$#3$}(2);
    \draw[CliqueEmptyEdge](1)edge[]node[CliqueLabel]{}(3);
    \draw[CliqueEmptyEdge](2)edge[]node[CliqueLabel]{}(3);
\end{tikzpicture}}

\newcommand{\TriangleEXE}[3]{
\begin{tikzpicture}[scale=.42,Centering]
    \node[CliquePoint](1)at(0,1){};
    \node[CliquePoint](2)at(0.87,-0.5){};
    \node[CliquePoint](3)at(-0.87,-0.5){};
    \draw[CliqueEmptyEdge](1)edge[]node[CliqueLabel]{}(2);
    \draw[CliqueEdge](1)edge[]node[CliqueLabel]{$#2$}(3);
    \draw[CliqueEmptyEdge](2)edge[]node[CliqueLabel]{}(3);
\end{tikzpicture}}

\newcommand{\TriangleEEE}[3]{
\begin{tikzpicture}[scale=.42,Centering]
    \node[CliquePoint](1)at(0,1){};
    \node[CliquePoint](2)at(0.87,-0.5){};
    \node[CliquePoint](3)at(-0.87,-0.5){};
    \draw[CliqueEmptyEdge](1)edge[]node[CliqueLabel]{}(2);
    \draw[CliqueEmptyEdge](1)edge[]node[CliqueLabel]{}(3);
    \draw[CliqueEmptyEdge](2)edge[]node[CliqueLabel]{}(3);
\end{tikzpicture}}

\newcommand{\SquareN}[4]{
\begin{tikzpicture}[scale=.6,Centering]
    \node[CliquePoint](1)at(-0.71,0.71){};
    \node[CliquePoint](2)at(0.71,0.71){};
    \node[CliquePoint](3)at(0.71,-0.71){};
    \node[CliquePoint](4)at(-0.71,-0.71){};
    \draw[CliqueEdge](1)edge[]node[CliqueLabel]{$#2$}(2);
    \draw[CliqueEdge](1)edge[]node[CliqueLabel]{$#1$}(4);
    \draw[CliqueEdge](2)edge[]node[CliqueLabel]{$#3$}(3);
    \draw[CliqueEdge](3)edge[]node[CliqueLabel]{$#4$}(4);
\end{tikzpicture}}

\newcommand{\SquareLeft}[5]{
\begin{tikzpicture}[scale=.6,Centering]
    \node[CliquePoint](1)at(-0.71,0.71){};
    \node[CliquePoint](2)at(0.71,0.71){};
    \node[CliquePoint](3)at(0.71,-0.71){};
    \node[CliquePoint](4)at(-0.71,-0.71){};
    \draw[CliqueEdge](1)edge[]node[CliqueLabel]{$#2$}(2);
    \draw[CliqueEdge](1)edge[]node[CliqueLabel]{$#1$}(4);
    \draw[CliqueEdge](2)edge[]node[CliqueLabel]{$#3$}(3);
    \draw[CliqueEdge](3)edge[]node[CliqueLabel]{$#4$}(4);
    \draw[CliqueEdge](1)edge[]node[CliqueLabel]{$#5$}(3);
\end{tikzpicture}}

\newcommand{\SquareRight}[5]{
\begin{tikzpicture}[scale=.6,Centering]
    \node[CliquePoint](1)at(-0.71,0.71){};
    \node[CliquePoint](2)at(0.71,0.71){};
    \node[CliquePoint](3)at(0.71,-0.71){};
    \node[CliquePoint](4)at(-0.71,-0.71){};
    \draw[CliqueEdge](1)edge[]node[CliqueLabel]{$#2$}(2);
    \draw[CliqueEdge](1)edge[]node[CliqueLabel]{$#1$}(4);
    \draw[CliqueEdge](2)edge[]node[CliqueLabel]{$#3$}(3);
    \draw[CliqueEdge](3)edge[]node[CliqueLabel]{$#4$}(4);
    \draw[CliqueEdge](2)edge[]node[CliqueLabel]{$#5$}(4);
\end{tikzpicture}}

\newcommand{\SquareMotz}{
\begin{tikzpicture}[scale=.53,Centering]
    \node[CliquePoint](1)at(-0.71,0.71){};
    \node[CliquePoint](2)at(0.71,0.71){};
    \node[CliquePoint](3)at(0.71,-0.71){};
    \node[CliquePoint](4)at(-0.71,-0.71){};
    \draw[CliqueEmptyEdge](2)edge[]node[CliqueLabel]{}(3);
    \draw[CliqueEmptyEdge](1)edge[]node[CliqueLabel]{}(4);
    \draw[CliqueEmptyEdge](3)edge[]node[CliqueLabel]{}(4);
    \draw[CliqueEdge](1)edge[]node[CliqueLabel]{$0$}(2);
\end{tikzpicture}}

\newcommand{\TriangleOp}[3]{\;
\begin{tikzpicture}[scale=.35,Centering]
    \node[shape=coordinate](1)at(0,1){};
    \node[shape=coordinate](2)at(0.87,-0.5){};
    \node[shape=coordinate](3)at(-0.87,-0.5){};
    \draw[draw=ColB!90](1)edge[]node[CliqueLabel,font=\tiny]{$#3$}(2);
    \draw[draw=ColB!90](1)edge[]node[CliqueLabel,font=\tiny]{$#2$}(3);
    \draw[draw=ColB!90](2)edge[]node[CliqueLabel,font=\tiny]{$#1$}(3);
\end{tikzpicture}\;}

\tikzstyle{CliqueEdge}=[draw=ColB!90,thick]
\tikzstyle{CliqueEdgeColorA}=[CliqueEdge,draw=ColA!80,fill=ColA!8]
\tikzstyle{CliqueEmptyEdge}=[draw=ColA!90,thick,densely dashed]
\tikzstyle{CliqueLabel}=[midway,inner sep=1pt,fill=ColA!0,font=\scriptsize]
\tikzstyle{CliquePoint}=[circle,inner sep=1pt,fill=ColB!25,draw=ColB!70]
\tikzstyle{CliqueEdgeGray}=[ColA!30,draw,cap=round]
\tikzstyle{CliqueEdgeBlue}=[ColA!80,thick,draw,cap=round]
\tikzstyle{CliqueEdgeRed}=[ColB!80,thick,draw,cap=round,dotted]
\tikzstyle{EdgeLabel}=[fill=ColA!0,inner sep=1pt,font=\tiny,midway]

\title[Operads of decorated cliques II]
    {Operads of decorated cliques II: Noncrossing cliques \vspace{-4ex}}
\date{\today}
\author[S. Giraudo]{Samuele Giraudo}
\address{\scriptsize
    Université du Québec à Montréal, LACIM, \\
    Pavillon Président-Kennedy, 201 Avenue du Président-Kennedy, Montréal, H2X~3Y7, Canada.
    \\ \vspace{-1ex}
    {\tt \href{mailto:giraudo.samuele@uqam.ca}{giraudo.samuele@uqam.ca}}}
\thanks{%
    This research has been partially supported by the projects CARPLO (ANR-20-CE40-0007)
    and ALCOHOL (ANR-19-CE40-0006) of the Agence nationale de la recherche.}

\begin{document}

\begin{abstract}
    A complete study of an operad $\NC\Mca$ of noncrossing configurations of chords
    introduced in previous work of the author is performed. This operad is defined on the
    linear span of all noncrossing $\Mca$-cliques. These are noncrossing configurations of
    chords with arcs labeled by a unitary magma $\Mca$. The magmatic product of $\Mca$
    intervenes for the computation of the operadic composition of $\Mca$-cliques. We show
    that this operad is binary, quadratic, and Koszul by considering techniques coming from
    rewrite systems on trees. We also compute a presentation for its Koszul dual. Finally,
    we explain how $\NC\Mca$ allows one to obtain alternative constructions of already known
    operads like operads of formal fractions and the operad of bicolored noncrossing
    configurations.
\end{abstract}

\MakeFirstPage

\section{Introduction}
Noncrossing configurations of chords on regular polygons are combinatorial objects appearing
in various contexts (see for instance~\cite{FN99,CP92,DRS10,PR14}). A simple generalization
of these consists in allowing several colors for the arcs of the configurations. This has
been considered in previous work~\cite{Cliques1} of the author, where the arcs are labeled
by elements of a unitary magma $\Mca$. More precisely, given a unitary magma $\Mca$,
$\Cliques_\Mca$ is the set of regular polygons $\Pfr$ where all possible arcs of $\Pfr$
are labeled by elements of~$\Mca$. When an arc is labeled by the unit $\Unit_\Mca$ of
$\Mca$, it is considered as missing, so that the notion of noncrossing configurations makes
sense. More precisely, noncrossing configurations have no diagonals labeled by elements
different from $\Unit_\Mca$ crossing another such diagonal. These objects have been named
noncrossing $\Mca$-cliques and are the main combinatorial objects studied in the present
article.

The linear span $\Cli\Mca$ of all $\Mca$-cliques (by relaxing the noncrossing condition)
over any field of characteristic zero forms an operad with rich algebraic and combinatorial
properties, introduced and studied in~\cite{Cliques1}. The operad structure added on these
objects allows one to compose two $\Mca$-cliques $\Pfr$ and $\Qfr$ by gluing a special arc
of $\Qfr$ (called the base) onto a selected edge of $\Pfr$. The magmatic product of $\Mca$
encodes how to relabel some edges of the resulting $\Mca$-clique. We showed in this previous
work that the subspace $\NC\Mca$ of $\Cli\Mca$ generated by all noncrossing $\Mca$-cliques
forms a suboperad of $\Cli\Mca$. The purpose of the present paper is to perform a complete
study of this operad. A first particularity of $\NC\Mca$ motivating this objective is that
it has a special status among all the suboperads of $\Cli\Mca$. Indeed, if $\Mca$ is
nontrivial, then $\Cli\Mca$ is not a binary operad, and $\NC\Mca$ is precisely the biggest
binary suboperad of~$\Cli\Mca$.

Let us now give an overview of the main properties of $\NC\Mca$ and the main results
contained in this article. First, by considering the dual trees of noncrossing
$\Mca$-cliques, we can view each noncrossing $\Mca$-clique of $\NC\Mca$ as a Schröder tree
(an ordered tree where all internal nodes have two or more children) with edges labeled by
$\Mca$ satisfying some conditions. Under this point of view, $\NC\Mca$ is an operad of such
Schröder trees endowed with a composition operation which is essentially a grafting
operation with a relabeling of edges or a contraction of edges. As a consequence of this
alternative combinatorial realization of $\NC\Mca$, we obtain a formula for its dimensions
involving Narayana numbers~\cite{Nar55}. To complete the study of $\NC\Mca$, a natural and
important question is to exhibit one of its presentations by generators and relations. In
order to compute the space of relations of $\NC\Mca$, we use techniques of rewrite
systems of trees~\cite{BN98}. Thus, we define a convergent rewrite rule $\Rew$ and show that
the space induced by $\Rew$ is the space of relations of $\NC\Mca$, leading to a
presentation by generators and relations of $\NC\Mca$. As an important consequence, this
proves that $\NC\Mca$ is always quadratic, regardless of $\Mca$. The existence of such a
convergent orientation of the space of relations of $\NC\Mca$ implies also by~\cite{Hof10}
that this operad is Koszul.

We also study some structures related with $\NC\Mca$. This includes the suboperads of
$\NC\Mca$ generated by some finite families of bubbles (the latter being noncrossing
$\Mca$-cliques with no diagonals). Under some conditions on the considered sets of bubbles,
we can describe the Hilbert series of these suboperads of $\NC\Mca$ by a system of algebraic
equations. We give two examples of suboperads of $\NC\Mca$ generated by some subsets of
triangles, including one which is a suboperad of $\NC\Dbb_0$ ($\Dbb_0$ is the multiplicative
monoid on $\{0, 1\}$) isomorphic to the operad $\Motz$ of Motzkin paths defined
in~\cite{Gir15}. Moreover, since $\NC\Mca$ is a binary and quadratic operad, its Koszul dual
$\NC\Mca^!$ is well-defined~\cite{GK94}. We compute its presentation, present an algebraic
equation for its Hilbert series, give a formula for its dimensions, and establish a
combinatorial realization of $\NC\Mca^!$ as a graded space involving dual $\Mca$-cliques,
which are $\Mca^2$-cliques with some constraints for the labels of their arcs.

Furthermore, by selecting appropriate unitary magmas $\Mca$, it is possible to provide
alternative constructions of already known operads as suboperads of $\NC\Mca$. We hence
construct the operad $\NCP$ of based noncrossing trees~\cite{Cha07,Ler11}, the suboperad
$\FF_4$ of the operad of formal fractions $\FF$~\cite{CHN16}, and the operad of bicolored
noncrossing configurations $\BNC$~\cite{CG14}. As a consequence of this last construction,
all the suboperads of $\BNC$ can be obtained from the construction $\NC$. This includes for
example the operad of noncrossing plants~\cite{Cha07}, the dipterous
operad~\cite{LR03,Zin12}, and the $2$-associative operad~\cite{LR06,Zin12}.

This text is organized as follows. Section~\ref{sec:definitions_trees} sets our notations
about trees, syntax trees, rewrite rules on trees, free operads, and Koszul duality of
operads. In Section~\ref{sec:operad_noncrossing}, we perform the aforementioned study of the
operad $\NC\Mca$ and in Section~\ref{sec:dual_NC_M}, of $\NC\Mca^!$. Finally, in
Section~\ref{sec:concrete_constructions_nc}, we use the construction $\NC$ to provide
alternative definitions of some known operads.

This paper is an extended version of~\cite{Gir17}, containing the proofs of the presented
results. It is also a sequel of~\cite{Cliques1}.

\subsubsection*{General notations and conventions}
All the algebraic structures of this article have a field of characteristic zero $\K$ as
ground field. For a set $S$, $\K \Angle{S}$ denotes the linear span of the elements of $S$.
For integers $a$ and $c$, $[a, c]$ denotes the set $\{b \in \N : a \leq b \leq c\}$ and
$[n]$ is short for the set $[1, n]$. The cardinality of a finite set $S$ is denoted by~$\#
S$. For a set $A$, $A^*$ denotes the set of finite sequences, called words, of elements
of~$A$. For an integer $n \geq 0$, $A^n$ (resp.\ $A^{\geq n}$) is the set of words on $A$ of
length $n$ (resp.\ at least~$n$). The word of length $0$ is the empty word denoted by
$\epsilon$. If $u$ is a word, its letters are indexed from left to right from $1$ to its
length $|u|$. For $i \in [|u|]$, $u_i$ is the letter of $u$ at position $i$. If $a$ is a
letter and $n$ is a nonnegative integer, $a^n$ denotes the word consisting in $n$
occurrences of $a$. For a letter $a$, $|u|_a$ denotes the number of occurrences of $a$
in~$u$.

\section{Elementary definitions and tools} \label{sec:definitions_trees}
The main purposes of this section are to provide tools to compute presentations and to prove
Koszulity of operads. For this, it is important to have precise definitions about free
operads, trees, and rewrite rules on trees at hand.

\subsection{Trees and rewrite rules}
Unless otherwise specified, we use in the sequel the standard terminology ({\em i.e.},
\Def{node}, \Def{edge}, \Def{root}, \Def{child}, {\em etc.})\ about ordered
trees~\cite{Knu97}. For the sake of completeness, we recall the most important definitions
and set our notations.

\subsubsection{Trees}
Let $\Tfr$ be an ordered tree. The \Def{arity} of a node of $\Tfr$ is its number of
children. An \Def{internal node} (resp.\ a \Def{leaf}) of $\Tfr$ is a node with a nonzero
(resp.\ null) arity. Internal nodes can be \Def{labeled}, that is, each internal node of a
tree is associated with an element of a certain set. Given an internal node $x$ of $\Tfr$,
the children of $x$ are by definition totally ordered from left to right
and are thus indexed from $1$ to the arity $\ell$ of $x$. For $i \in [\ell]$, the
\Def{$i$th subtree} of $\Tfr$ is the tree rooted at the $i$th child of $\Tfr$. Similarly,
the leaves of $\Tfr$ are totally ordered from left to right and thus are indexed from $1$ to
the number of its leaves. In our graphical representations, each ordered tree is depicted so
that its root is the uppermost node. Since we consider in the sequel only ordered rooted
trees, we shall call these simply \Def{trees}.

\subsubsection{Syntax trees}
Let $G := \bigsqcup_{n \geq 1} G(n)$ be a graded set. The \Def{arity} of an element $x$ of
$G$ is $n$ provided that $x \in G(n)$. A \Def{syntax tree} on $G$ is a tree such that its
internal nodes of arity $n$ are labeled by elements of arity $n$ of~$G$. The \Def{degree}
(resp.\ \Def{arity}) of a syntax tree is its number of internal nodes (resp.\ leaves). For
instance, if $G := G(2) \sqcup G(3)$ with $G(2) := \{\Att, \Ctt\}$ and $G(3) := \{\Btt\}$,
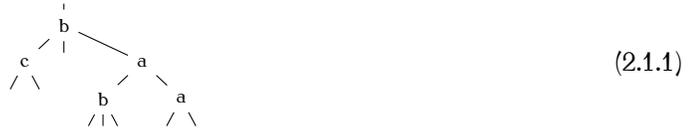
\begin{equation}
    \begin{tikzpicture}[xscale=.26,yscale=.15,Centering]
        \node(0)at(0.00,-6.50){};
        \node(10)at(8.00,-9.75){};
        \node(12)at(10.00,-9.75){};
        \node(2)at(2.00,-6.50){};
        \node(4)at(3.00,-3.25){};
        \node(5)at(4.00,-9.75){};
        \node(7)at(5.00,-9.75){};
        \node(8)at(6.00,-9.75){};
        \node[NodeST](1)at(1.00,-3.25){$\Ctt$};
        \node[NodeST](11)at(9.00,-6.50){$\Att$};
        \node[NodeST](3)at(3.00,0.00){$\Btt$};
        \node[NodeST](6)at(5.00,-6.50){$\Btt$};
        \node[NodeST](9)at(7.00,-3.25){$\Att$};
        \node(r)at(3.00,2.75){};
        \draw(0)--(1); \draw(1)--(3); \draw(10)--(11); \draw(11)--(9);
        \draw(12)--(11); \draw(2)--(1); \draw(4)--(3); \draw(5)--(6);
        \draw(6)--(9); \draw(7)--(6); \draw(8)--(6); \draw(9)--(3);
        \draw(r)--(3);
    \end{tikzpicture}
\end{equation}
is a syntax tree on $G$ of degree $5$ and arity $8$. Its root is labeled by $\Btt$ and has
arity~$3$.

A syntax tree $\Sfr$ is a \Def{subtree} of a syntax tree $\Tfr$ if it is possible to fit
$\Sfr$ at a certain place of $\Tfr$, by possibly superimposing leaves of $\Sfr$ and internal
nodes of $\Tfr$. In this case, we say that $\Tfr$ \Def{admits an occurrence} of (the
\Def{pattern}) $\Sfr$. Conversely, we say that $\Tfr$ \Def{avoids} $\Sfr$ if there is no
occurrence of $\Sfr$ in~$\Tfr$.

\subsubsection{Rewrite rules}
Let $S$ be a set of trees. A \Def{rewrite rule} on $S$ is a binary relation $\Rew$ on $S$
which has the property that, if $\Sfr \Rew \Sfr'$ for two trees $\Sfr$ and $\Sfr'$, then
$\Sfr$ and $\Sfr'$ have the same number of leaves. We say that a tree $\Tfr$ is
\Def{rewritable in one step} into $\Tfr'$ by $\Rew$ if there exist two trees $\Sfr$ and
$\Sfr'$ satisfying $\Sfr \Rew \Sfr'$ and $\Tfr$ has a subtree $\Sfr$ such that, by replacing
$\Sfr$ by $\Sfr'$ in $\Tfr$, we obtain $\Tfr'$. We denote this property by $\Tfr \CRew
\Tfr'$, so that $\CRew$ is a binary relation on $S$. If $\Tfr = \Tfr'$ or if there
exists a sequence of trees $(\Tfr_1, \dots, \Tfr_{k - 1})$ with $k \geq 1$ such that
\begin{math}
    \Tfr \CRew \Tfr_1 \CRew \cdots \CRew \Tfr_{k - 1} \CRew \Tfr',
\end{math}
we say that $\Tfr$ is \Def{rewritable} by $\CRew$ into $\Tfr'$, and we denote this property
by $\Tfr \CRewTrans \Tfr'$. In other words, $\CRewTrans$ is the reflexive and transitive
closure of $\CRew$. We write $\RewTrans$ for the reflexive and transitive closure of $\Rew$,
and we write $\RewTransSym$ (resp.\ $\CRewTransSym$) for the reflexive, transitive, and
symmetric closure of $\Rew$ (resp.\ $\CRew$). The \Def{vector space induced} by $\Rew$ is
the subspace of the linear span $\K \Angle{S}$ of all trees of $S$ generated by the family
of all $\Tfr - \Tfr'$ such that~$\Tfr \RewTransSym \Tfr'$.

For instance, let $S$ be the set of trees where internal nodes are labeled by $\{\Att, \Btt,
\Ctt\}$ and consider the rewrite rule $\Rew$ on $S$ satisfying
\vspace{-2ex}
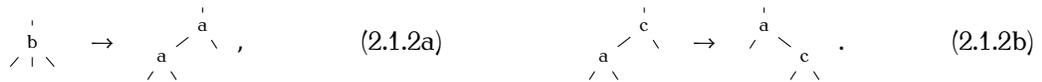
\begin{multicols}{2}
\begin{subequations}
\begin{equation}
    \begin{tikzpicture}[xscale=.4,yscale=.23,Centering]
        \node(0)at(0.00,-2.00){};
        \node(2)at(1.00,-2.00){};
        \node(3)at(2.00,-2.00){};
        \node[NodeST](1)at(1.00,0.00){$\Btt$};
        \draw(0)--(1);
        \draw(2)--(1);
        \draw(3)--(1);
        \node(r)at(1.00,1.75){};
        \draw(r)--(1);
    \end{tikzpicture}
    \enspace \Rew \enspace
    \begin{tikzpicture}[xscale=.27,yscale=.25,Centering]
        \node(0)at(0.00,-3.33){};
        \node(2)at(2.00,-3.33){};
        \node(4)at(4.00,-1.67){};
        \node[NodeST](1)at(1.00,-1.67){$\Att$};
        \node[NodeST](3)at(3.00,0.00){$\Att$};
        \draw(0)--(1);
        \draw(1)--(3);
        \draw(2)--(1);
        \draw(4)--(3);
        \node(r)at(3.00,1.5){};
        \draw(r)--(3);
    \end{tikzpicture}\,,
\end{equation}

\begin{equation}
    \begin{tikzpicture}[xscale=.27,yscale=.25,Centering]
        \node(0)at(0.00,-3.33){};
        \node(2)at(2.00,-3.33){};
        \node(4)at(4.00,-1.67){};
        \node[NodeST](1)at(1.00,-1.67){$\Att$};
        \node[NodeST](3)at(3.00,0.00){$\Ctt$};
        \draw(0)--(1);
        \draw(1)--(3);
        \draw(2)--(1);
        \draw(4)--(3);
        \node(r)at(3.00,1.5){};
        \draw(r)--(3);
    \end{tikzpicture}
    \enspace \Rew \enspace
    \begin{tikzpicture}[xscale=.27,yscale=.25,Centering]
        \node(0)at(0.00,-1.67){};
        \node(2)at(2.00,-3.33){};
        \node(4)at(4.00,-3.33){};
        \node[NodeST](1)at(1.00,0.00){$\Att$};
        \node[NodeST](3)at(3.00,-1.67){$\Ctt$};
        \draw(0)--(1);
        \draw(2)--(3);
        \draw(3)--(1);
        \draw(4)--(3);
        \node(r)at(1.00,1.5){};
        \draw(r)--(1);
    \end{tikzpicture}\,.
\end{equation}
\end{subequations}
\end{multicols}
\noindent We then have the following steps of rewritings by $\Rew$:
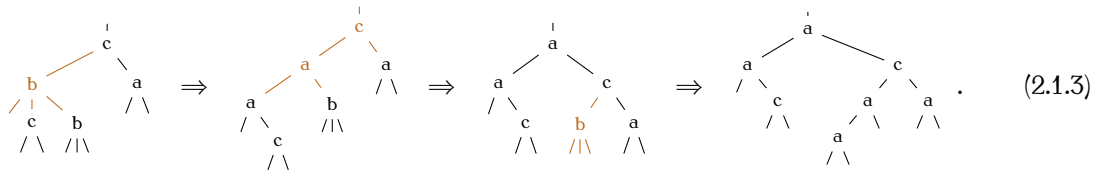
\begin{equation}
    \begin{tikzpicture}[xscale=.2,yscale=.16,Centering]
        \node(0)at(0.00,-6.50){};
        \node(10)at(8.00,-6.50){};
        \node(12)at(10.00,-6.50){};
        \node(2)at(1.00,-9.75){};
        \node(4)at(3.00,-9.75){};
        \node(5)at(4.00,-9.75){};
        \node(7)at(5.00,-9.75){};
        \node(8)at(6.00,-9.75){};
        \node[NodeST,text=ColB!90](1)at(2.00,-3.25){$\Btt$};
        \node[NodeST](11)at(9.00,-3.25){$\Att$};
        \node[NodeST](3)at(2.00,-6.50){$\Ctt$};
        \node[NodeST](6)at(5.00,-6.50){$\Btt$};
        \node[NodeST](9)at(7.00,0.00){$\Ctt$};
        \draw[draw=ColB!90](0)--(1);
        \draw[draw=ColB!90](1)--(9);
        \draw(10)--(11);
        \draw(11)--(9);
        \draw(12)--(11);
        \draw(2)--(3);
        \draw[draw=ColB!90](3)--(1);
        \draw(4)--(3);
        \draw(5)--(6);
        \draw[draw=ColB!90](6)--(1);
        \draw(7)--(6);
        \draw(8)--(6);
        \node(r)at(7.00,2.5){};
        \draw(r)--(9);
    \end{tikzpicture}
    \enspace \CRew \enspace
    \begin{tikzpicture}[xscale=.18,yscale=.18,Centering]
        \node(0)at(0.00,-8.40){};
        \node(11)at(10.00,-5.60){};
        \node(13)at(12.00,-5.60){};
        \node(2)at(2.00,-11.20){};
        \node(4)at(4.00,-11.20){};
        \node(6)at(6.00,-8.40){};
        \node(8)at(7.00,-8.40){};
        \node(9)at(8.00,-8.40){};
        \node[NodeST](1)at(1.00,-5.60){$\Att$};
        \node[NodeST,text=ColB!90](10)at(9.00,0.0){$\Ctt$};
        \node[NodeST](12)at(11.00,-2.80){$\Att$};
        \node[NodeST](3)at(3.00,-8.40){$\Ctt$};
        \node[NodeST,text=ColB!90](5)at(5.00,-2.80){$\Att$};
        \node[NodeST](7)at(7.00,-5.60){$\Btt$};
        \draw(0)--(1);
        \draw[draw=ColB!90](1)--(5);
        \draw(11)--(12);
        \draw[draw=ColB!90](12)--(10);
        \draw(13)--(12);
        \draw(2)--(3);
        \draw(3)--(1);
        \draw(4)--(3);
        \draw[draw=ColB!90](5)--(10);
        \draw(6)--(7);
        \draw[draw=ColB!90](7)--(5);
        \draw(8)--(7);
        \draw(9)--(7);
        \node(r)at(9.00,2.25){};
        \draw[draw=ColB!90](r)--(10);
    \end{tikzpicture}
    \enspace \CRew \enspace
    \begin{tikzpicture}[xscale=.18,yscale=.15,Centering]
        \node(0)at(0.00,-7.00){};
        \node(11)at(10.00,-10.50){};
        \node(13)at(12.00,-10.50){};
        \node(2)at(2.00,-10.50){};
        \node(4)at(4.00,-10.50){};
        \node(6)at(6.00,-10.50){};
        \node(8)at(7.00,-10.50){};
        \node(9)at(8.00,-10.50){};
        \node[NodeST](1)at(1.00,-3.50){$\Att$};
        \node[NodeST](10)at(9.00,-3.50){$\Ctt$};
        \node[NodeST](12)at(11.00,-7.00){$\Att$};
        \node[NodeST](3)at(3.00,-7.00){$\Ctt$};
        \node[NodeST](5)at(5.00,0.00){$\Att$};
        \node[NodeST,text=ColB!90](7)at(7.00,-7.00){$\Btt$};
        \draw(0)--(1);
        \draw(1)--(5);
        \draw(10)--(5);
        \draw(11)--(12);
        \draw(12)--(10);
        \draw(13)--(12);
        \draw(2)--(3);
        \draw(3)--(1);
        \draw(4)--(3);
        \draw[draw=ColB!90](6)--(7);
        \draw[draw=ColB!90](7)--(10);
        \draw[draw=ColB!90](8)--(7);
        \draw[draw=ColB!90](9)--(7);
        \node(r)at(5.00,2.75){};
        \draw(r)--(5);
    \end{tikzpicture}
    \enspace \CRew \enspace
    \begin{tikzpicture}[xscale=.2,yscale=.16,Centering]
        \node(0)at(0.00,-6.00){};
        \node(10)at(10.00,-9.00){};
        \node(12)at(12.00,-9.00){};
        \node(14)at(14.00,-9.00){};
        \node(2)at(2.00,-9.00){};
        \node(4)at(4.00,-9.00){};
        \node(6)at(6.00,-12.00){};
        \node(8)at(8.00,-12.00){};
        \node[NodeST](1)at(1.00,-3.00){$\Att$};
        \node[NodeST](11)at(11.00,-3.00){$\Ctt$};
        \node[NodeST](13)at(13.00,-6.00){$\Att$};
        \node[NodeST](3)at(3.00,-6.00){$\Ctt$};
        \node[NodeST](5)at(5.00,0.00){$\Att$};
        \node[NodeST](7)at(7.00,-9.00){$\Att$};
        \node[NodeST](9)at(9.00,-6.00){$\Att$};
        \draw(0)--(1);
        \draw(1)--(5);
        \draw(10)--(9);
        \draw(11)--(5);
        \draw(12)--(13);
        \draw(13)--(11);
        \draw(14)--(13);
        \draw(2)--(3);
        \draw(3)--(1);
        \draw(4)--(3);
        \draw(6)--(7);
        \draw(7)--(9);
        \draw(8)--(7);
        \draw(9)--(11);
        \node(r)at(5.00,2.25){};
        \draw(r)--(5);
    \end{tikzpicture}\,.
\end{equation}

We shall use the standard terminology (\Def{terminating}, \Def{normal form},
\Def{confluent}, \Def{convergent}, {\em etc.})\ about rewrite rules~\cite{BN98}. Let us
recall now the most important definitions. Let $\Rew$ be a rewrite rule on a set $S$ of
trees. We say that $\Rew$ is \Def{terminating} if there is no infinite chain $\Tfr \CRew
\Tfr_1 \CRew \Tfr_2 \CRew \cdots$. In this case, any tree $\Tfr$ of $S$ that cannot be
rewritten by $\Rew$ is a \Def{normal form} for $\Rew$. We say that $\Rew$ is \Def{confluent}
if for any trees $\Tfr$, $\Rfr_1$, and $\Rfr_2$ such that $\Tfr \CRewTrans \Rfr_1$ and $\Tfr
\CRewTrans \Rfr_2$, there exists a tree $\Tfr'$ such that $\Rfr_1 \CRewTrans \Tfr'$ and
$\Rfr_2 \CRewTrans \Tfr'$. If $\Rew$ is both terminating and confluent, then $\Rew$ is
\Def{convergent}.

\subsection{Free operads and Koszul duality}
All notations and conventions about operads come from Section~\ref{subsec:ns_operads}
of~\cite{Cliques1}.

\subsubsection{Free operads}
Let $\Gen := \bigoplus_{n \geq 1} \Gen(n)$ be a graded vector space. In particular, $\Gen$
is a graded set so that we can consider syntax trees on $\Gen$. The \Def{free operad} over
$\Gen$ is the operad $\Free(\Gen)$ wherein, for $n \geq 1$, $\Free(\Gen)(n)$ is the
linear span of the syntax trees on $\Gen$ of arity $n$. The labeling of the internal nodes
of the trees of $\Free(\Gen)$ is linear in the sense that if $\Tfr$ is a syntax tree on
$\Gen$ having an internal node labeled by $x + \lambda y \in \Gen$, $\lambda \in \K$, then,
in $\Free(\Gen)$, we have $\Tfr = \Tfr_x + \lambda \Tfr_y$, where $\Tfr_x$ (resp.\ $\Tfr_y$)
is the tree obtained by labeling by $x$ (resp.\ $y$) the considered node labeled by $x +
\lambda y$ in $\Tfr$. The partial composition $\Sfr \circ_i \Tfr$ of $\Free(\Gen)$ of two
syntax trees $\Sfr$ and $\Tfr$ on $\Gen$ consists in grafting the root of $\Tfr$ on the
$i$th leaf of $\Sfr$. The unit $\Leaf$ of $\Free(\Gen)$ is the tree consisting in one leaf.
For instance, by setting $\Gen := \K \Angle{G}$ where $G$ is the graded set defined in the
previous example, in~$\Free(\Gen)$ we have
\begin{equation}
    \begin{tikzpicture}[xscale=.28,yscale=.17,Centering]
        \node(0)at(0.00,-5.33){};
        \node(2)at(2.00,-5.33){};
        \node(4)at(4.00,-5.33){};
        \node(6)at(5.00,-5.33){};
        \node(7)at(6.00,-5.33){};
        \node[NodeST,text=ColA](1)at(1.00,-2.67){$\Att$};
        \node[NodeST,text=ColA](3)at(3.00,0.00){$\Att$};
        \node[NodeST,text=ColA](5)at(5.00,-2.67){$\Btt$};
        \node(r)at(3.00,2.25){};
        \draw[draw=ColA](0)--(1); \draw[draw=ColA](1)--(3);
        \draw[draw=ColA](2)--(1); \draw[draw=ColA](4)--(5);
        \draw[draw=ColA](5)--(3); \draw[draw=ColA](6)--(5);
        \draw[draw=ColA](7)--(5); \draw[draw=ColA](r)--(3);
    \end{tikzpicture}
    \enspace \circ_3 \enspace
    \begin{tikzpicture}[xscale=.28,yscale=.28,Centering]
        \node(0)at(0.00,-1.67){};
        \node(2)at(2.00,-3.33){};
        \node(4)at(4.00,-3.33){};
        \node[NodeST,text=ColB](1)at(1.00,0.00){$\Ctt$};
        \node[NodeST,text=ColB](3)at(3.00,-1.67){$\Att + \Ctt$};
        \node(r)at(1.00,1.5){};
        \draw[draw=ColB](0)--(1); \draw[draw=ColB](2)--(3);
        \draw[draw=ColB](3)--(1); \draw[draw=ColB](4)--(3);
        \draw[draw=ColB](r)--(1);
    \end{tikzpicture}
    \enspace = \enspace
    \begin{tikzpicture}[xscale=.25,yscale=.18,Centering]
        \node(0)at(0.00,-4.80){};
        \node(10)at(9.00,-4.80){};
        \node(11)at(10.00,-4.80){};
        \node(2)at(2.00,-4.80){};
        \node(4)at(4.00,-7.20){};
        \node(6)at(6.00,-9.60){};
        \node(8)at(8.00,-9.60){};
        \node[NodeST,text=ColA](1)at(1.00,-2.40){$\Att$};
        \node[NodeST,text=ColA](3)at(3.00,0.00){$\Att$};
        \node[NodeST,text=ColB](5)at(5.00,-4.80){$\Ctt$};
        \node[NodeST,text=ColB](7)at(7.00,-7.20){$\Att$};
        \node[NodeST,text=ColA](9)at(9.00,-2.40){$\Btt$};
        \node(r)at(3.00,2.40){};
        \draw[draw=ColA](0)--(1);
        \draw[draw=ColA](1)--(3);
        \draw[draw=ColA](10)--(9);
        \draw[draw=ColA](11)--(9);
        \draw[draw=ColA](2)--(1);
        \draw[draw=ColB](4)--(5);
        \draw[draw=ColB](5)--(9);
        \draw[draw=ColB](6)--(7);
        \draw[draw=ColB](7)--(5);
        \draw[draw=ColB](8)--(7);
        \draw[draw=ColA](9)--(3);
        \draw[draw=ColA](r)--(3);
    \end{tikzpicture}
    \enspace + \enspace
    \begin{tikzpicture}[xscale=.25,yscale=.18,Centering]
        \node(0)at(0.00,-4.80){};
        \node(10)at(9.00,-4.80){};
        \node(11)at(10.00,-4.80){};
        \node(2)at(2.00,-4.80){};
        \node(4)at(4.00,-7.20){};
        \node(6)at(6.00,-9.60){};
        \node(8)at(8.00,-9.60){};
        \node[NodeST,text=ColA](1)at(1.00,-2.40){$\Att$};
        \node[NodeST,text=ColA](3)at(3.00,0.00){$\Att$};
        \node[NodeST,text=ColB](5)at(5.00,-4.80){$\Ctt$};
        \node[NodeST,text=ColB](7)at(7.00,-7.20){$\Ctt$};
        \node[NodeST,text=ColA](9)at(9.00,-2.40){$\Btt$};
        \node(r)at(3.00,2.40){};
        \draw[draw=ColA](0)--(1);
        \draw[draw=ColA](1)--(3);
        \draw[draw=ColA](10)--(9);
        \draw[draw=ColA](11)--(9);
        \draw[draw=ColA](2)--(1);
        \draw[draw=ColB](4)--(5);
        \draw[draw=ColB](5)--(9);
        \draw[draw=ColB](6)--(7);
        \draw[draw=ColB](7)--(5);
        \draw[draw=ColB](8)--(7);
        \draw[draw=ColA](9)--(3);
        \draw[draw=ColA](r)--(3);
    \end{tikzpicture}\,.
\end{equation}

We denote by $\Corolla : \Gen \to \Free(\Gen)$ the inclusion map, sending any $x$ of $\Gen$
to the \Def{corolla} labeled by $x$, that is, the syntax tree consisting in a single
internal node labeled by $x$ attached to a required number of leaves. In the sequel, if
required by the context, we shall implicitly view an element $x$ of $\Gen$ as the corolla
$\Corolla(x)$ of $\Free(\Gen)$. For instance, given two elements $x$ and $y$ of $\Gen$, we
shall denote the syntax tree $\Corolla(x) \circ_i \Corolla(y)$ simply by $x \circ_i y$ for
all valid integers~$i$.

Free operads satisfy a universality property. Indeed, $\Free(\Gen)$ is the unique operad (up
to isomorphism) such that for an operad $\Oca$ and a linear map $f : \Gen \to \Oca$
respecting the arities, there exists a unique operad morphism $\phi : \Free(\Gen) \to \Oca$
such that $f = \phi \circ \Corolla$.

\subsubsection{Evaluations and treelike expressions}
Let us first fix a notation. If $\Oca$ is an operad, the \Def{complete composition map} of
$\Oca$ is the linear map
\begin{equation}
    \circ : \Oca(n) \otimes \Oca(m_1) \otimes \dots \otimes \Oca(m_n)
    \to \Oca(m_1 + \dots + m_n),
\end{equation}
defined, for $x \in \Oca(n)$ and $y_1, \dots, y_n \in \Oca$, by
\begin{equation}
    x \circ [y_1, \dots, y_n]
    := (\dots ((x \circ_n y_n) \circ_{n - 1} y_{n - 1}) \dots)
    \circ_1 y_1.
\end{equation}

For an operad $\Oca$, by viewing $\Oca$ as a graded vector space, $\Free(\Oca)$ is by
definition the free operad on $\Oca$. The \Def{evaluation map} of $\Oca$ is the linear map
\begin{equation}
    \Eval : \Free(\Oca) \to \Oca,
\end{equation}
defined recursively, for any syntax tree $\Tfr$ on $\Oca$, by
\begin{equation}
    \Eval(\Tfr) :=
    \begin{cases}
        \Unit \in \Oca,
            & \mbox{if } \Tfr = \Leaf, \\
        x \circ \left[\Eval\left(\Tfr_1\right), \dots,
        \Eval\left(\Tfr_k\right)\right],
            & \mbox{otherwise},
    \end{cases}
\end{equation}
where $x$ is the label of the root of $\Tfr$ and $\Tfr_1$, \dots, $\Tfr_k$ are, from left to
right, the subtrees of $\Tfr$. This map is the unique surjective operad morphism from
$\Free(\Oca)$ to $\Oca$ satisfying $\Eval(\Corolla(x)) = x$ for all $x \in \Oca$. If $S$ is
a subspace of $\Oca$, a \Def{treelike expression} on $S$ of $x \in \Oca$ is a tree $\Tfr$ of
$\Free(\Oca)$ such that $\Eval(\Tfr) = x$ and all internal nodes of $\Tfr$ are labeled
by~$S$.

\subsubsection{Presentations by generators and relations}
Let $G := \bigsqcup_{n \geq 1} G(n)$ be a graded set. Setting $\Gen := \K \Angle{G}$, we
denote the operad ideal of $\Free(\Gen)$ generated by the subspace $\Rel$ of $\Free(\Gen)$
by $\Angle{\Rel}$. Given an operad $\Oca$, the pair $(G, \Rel)$ is a \Def{presentation} of
$\Oca$ if $\Oca$ is isomorphic to $\Free(\Gen)/_{\langle \Rel \rangle}$. In this case, we
call $\Gen$ the \Def{space of generators} and $\Rel$ the \Def{space of relations} of $\Oca$.
We say that $\Oca$ is \Def{quadratic} if there is a presentation $(G, \Rel)$ of $\Oca$ such
that $\Rel$ is a homogeneous subspace of $\Free(\Gen)$ consisting in syntax trees of degree
$2$. Furthermore, we say that $\Oca$ is \Def{binary} if there is a presentation $(G, \Rel)$
of $\Oca$ such that $\Gen$ is concentrated in arity~$2$. Furthermore, if $\Oca$ admits a
presentation $(G, \Rel)$ and $\Rew$ is a rewrite rule on $\Free(\Gen)$ such that the space
induced by $\Rew$ is $\Rel$, we say that $\Rew$ is an \Def{orientation} of~$\Rel$.

\subsubsection{Koszul duality and Koszulity}%
\label{subsubsec:koszul_duality_koszulity_criterion}
In~\cite{GK94}, Ginzburg and Kapranov extended the notion of Koszul duality of quadratic
associative algebras to quadratic operads. Starting with an operad $\Oca$ admitting a binary
and quadratic presentation $(G, \Rel)$ where $G$ is finite, the \Def{Koszul dual} of $\Oca$
is the operad $\Oca^!$, isomorphic to the operad admitting the presentation $\left(G,
\Rel^\perp\right)$ where $\Rel^\perp$ is the annihilator of $\Rel$ in $\Free(\Gen)$, $\Gen$
being the space $\K \Angle{G}$, with respect to the bilinear map
\begin{equation}
    \langle -, - \rangle : \Free(\Gen)(3) \otimes \Free(\Gen)(3) \to \K
\end{equation}
defined, for all $x, x', y, y' \in \Gen(2)$, by
\begin{equation} \label{equ:scalar_product_koszul}
    \left\langle x \circ_i y, x' \circ_{i'} y' \right\rangle :=
    \begin{cases}
        1, & \mbox{if }
            x = x', y = y', \mbox{ and } i = i' = 1, \\
        -1, & \mbox{if }
            x = x', y = y', \mbox{ and } i = i' = 2, \\
        0, & \mbox{otherwise}.
    \end{cases}
\end{equation}
Then, with knowledge of a presentation of $\Oca$, one can compute a presentation
of~$\Oca^!$.

Recall that a quadratic operad $\Oca$ is \Def{Koszul} if its Koszul complex is
acyclic~\cite{GK94,LV12}. Furthermore, if $\Oca$ is Koszul and admits an Hilbert series,
then the Hilbert series of $\Oca$ and of its Koszul dual $\Oca^!$ are related~\cite{GK94} by
\begin{equation} \label{equ:Hilbert_series_Koszul_operads}
    \Hilbert_\Oca\left(-\Hilbert_{\Oca^!}(-t)\right) = t.
\end{equation}
Relation~\eqref{equ:Hilbert_series_Koszul_operads} can be either used to prove that an
operad is not Koszul (this is the case when the coefficients of the hypothetical Hilbert
series of the Koszul dual admit coefficients that are not nonnegative integers) or to
compute the Hilbert series of the Koszul dual of a Koszul operad.

Here, to prove the Koszulity of an operad $\Oca$, we shall make use of a tool
introduced by Dotsenko and Khoroshkin~\cite{DK10} in the context of Gröbner bases for
operads, which in our context can be reformulated in the following way by using rewrite
rules.
\begin{Lemma} \label{lem:koszulity_criterion_pbw}
    Let $\Oca$ be an operad admitting a quadratic presentation $(G, \Rel)$. If there exists
    an orientation $\Rew$ of $\Rel$ such that $\Rew$ is a convergent rewrite rule, then
    $\Oca$ is Koszul.
\end{Lemma}

If $\Rew$ satisfies the conditions contained in the statement of
Lemma~\ref{lem:koszulity_criterion_pbw}, then the set of normal forms of $\Rew$ forms a
basis of $\Oca$, called \Def{Poincaré--Birkhoff--Witt basis}. These bases arise from the work
of Hoffbeck~\cite{Hof10} (see also~\cite{LV12}).

\subsubsection{Algebras over operads}
An operad $\Oca$ encodes a category of algebras whose objects are called
\Def{$\Oca$-algebras}. An $\Oca$-algebra $\Alg_\Oca$ is a vector space endowed with a linear
left action
\begin{equation}
    \cdot : \Oca(n) \otimes \Alg_\Oca^{\otimes n} \to \Alg_\Oca,
    \qquad n \geq 1,
\end{equation}
satisfying the relations imposed by the structure of $\Oca$, which are
\begin{multline} \label{equ:algebra_over_operad}
    (x \circ_i y) \cdot
    \left(a_1 \otimes \dots \otimes a_{n + m - 1}\right)
    \\
    =
    x \cdot \left(a_1 \otimes \dots
        \otimes a_{i - 1} \otimes
        y \cdot \left(a_i \otimes \dots \otimes a_{i + m - 1}\right)
        \otimes a_{i + m} \otimes
        \dots \otimes a_{n + m - 1}\right),
\end{multline}
for all $x \in \Oca(n)$, $y \in \Oca(m)$, $i \in [n]$, and
\begin{math}
    a_1 \otimes \dots \otimes a_{n + m - 1} \in \Alg_\Oca^{\otimes {n + m - 1}}.
\end{math}

Notice that, by~\eqref{equ:algebra_over_operad}, if $G$ is a generating set of $\Oca$, it is
enough to define the action of each $x \in G$ on $\Alg_\Oca^{\otimes |x|}$ to wholly
define~$\cdot$. In other words, any element $x$ of $\Oca$ of arity $n$ plays the role of a
linear operation
\begin{equation}
    x : \Alg_\Oca^{\otimes n}  \to \Alg_\Oca,
\end{equation}
taking $n$ elements of $\Alg_\Oca$ as inputs and computing an element of $\Alg_\Oca$. By a
slight but convenient abuse of notation, for $x \in \Oca(n)$, we shall write $x(a_1,
\dots, a_n)$, or $a_1 \, x \, a_2$ if $x$ has arity $2$, for the element $x \cdot (a_1
\otimes \dots \otimes a_n)$ of $\Alg_\Oca$, for any $a_1 \otimes \dots \otimes a_n \in
\Alg_\Oca^{\otimes n}$. Observe that by~\eqref{equ:algebra_over_operad} an associative
element of $\Oca$ gives rise to an associative operation on~$\Alg_\Oca$.

\section{Operads of noncrossing decorated cliques} \label{sec:operad_noncrossing}
We use all the notations and definitions of Sections~\ref{subsec:configurations}
and~\ref{sec:construction_Cli} of~\cite{Cliques1} about decorated $\Mca$-cliques and the
$\Mca$-clique operad $\Cli\Mca$. We perform here a complete study of the suboperad
$\Cro_0\Mca$ of noncrossing $\Mca$-cliques defined in
Section~\ref{subsubsec:quotient_Cli_M_crossings} of the aforementioned paper. For
simplicity, this operad is denoted in the sequel as $\NC\Mca$ and called the
\Def{noncrossing $\Mca$-clique operad}. The process which produces from a unitary magma
$\Mca$ the operad $\NC\Mca$ is called the \Def{noncrossing clique construction}.

\subsection{General properties}
As shown in~\cite{Cliques1}, $\NC\Mca$ is an operad defined on the linear span of all
noncrossing $\Mca$-cliques and can be seen as a suboperad of $\Cli\Mca$ restrained on
$\Mca$-cliques with $0$ as crossing number. By definition of $\NC\Mca$, the partial
composition $\Pfr \circ_i \Qfr$ of two noncrossing $\Mca$-cliques $\Pfr$ and $\Qfr$ in
$\NC\Mca$ is equal to the partial composition $\Pfr \circ_i \Qfr$ in $\Cli\Mca$. Recall that
the partial composition $\Pfr \circ_i \Qfr$ is the noncrossing $\Mca$-clique obtained by
gluing the base of $\Qfr$ onto the $i$th edge of $\Pfr$ and by relabeling the common arcs
between $\Pfr$ and $\Qfr$, respectively the arcs $(i, i + 1)$ and $(1, m + 1)$, by $\Pfr_i
\Op \Qfr_0$, where $\Op$ is the magmatic product of $\Mca$. For instance, in $\Cli\Z$, we
have
\begin{equation}
    \begin{tikzpicture}[scale=.85,Centering]
        \node[CliquePoint](1)at(-0.50,-0.87){};
        \node[CliquePoint](2)at(-1.00,-0.00){};
        \node[CliquePoint](3)at(-0.50,0.87){};
        \node[CliquePoint](4)at(0.50,0.87){};
        \node[CliquePoint](5)at(1.00,0.00){};
        \node[CliquePoint](6)at(0.50,-0.87){};
        \draw[CliqueEdge](1)edge[]node[CliqueLabel]{$1$}(2);
        \draw[CliqueEdge](1)edge[bend left=30]node[CliqueLabel]{$-2$}(5);
        \draw[CliqueEmptyEdge](1)edge[]node[CliqueLabel]{}(6);
        \draw[CliqueEdge](2)edge[]node[CliqueLabel]{$-2$}(3);
        \draw[CliqueEmptyEdge](3)edge[]node[CliqueLabel]{}(4);
        \draw[CliqueEdge](3)edge[bend right=30]node[CliqueLabel]{$1$}(5);
        \draw[CliqueEmptyEdge](4)edge[]node[CliqueLabel]{}(5);
        \draw[CliqueEmptyEdge](5)edge[]node[CliqueLabel]{}(6);
    \end{tikzpicture}
    \enspace \circ_2 \enspace
    \begin{tikzpicture}[scale=.65,Centering]
        \node[CliquePoint](1)at(-0.71,-0.71){};
        \node[CliquePoint](2)at(-0.71,0.71){};
        \node[CliquePoint](3)at(0.71,0.71){};
        \node[CliquePoint](4)at(0.71,-0.71){};
        \draw[CliqueEmptyEdge](1)edge[]node[CliqueLabel]{}(2);
        \draw[CliqueEdge](1)edge[]node[CliqueLabel]{$1$}(3);
        \draw[CliqueEdge](1)edge[]node[CliqueLabel]{$3$}(4);
        \draw[CliqueEmptyEdge](2)edge[]node[CliqueLabel]{}(3);
        \draw[CliqueEdge](3)edge[]node[CliqueLabel]{$2$}(4);
    \end{tikzpicture}
    \enspace = \enspace
    \begin{tikzpicture}[scale=1.1,Centering]
        \node[CliquePoint](1)at(-0.38,-0.92){};
        \node[CliquePoint](2)at(-0.92,-0.38){};
        \node[CliquePoint](3)at(-0.92,0.38){};
        \node[CliquePoint](4)at(-0.38,0.92){};
        \node[CliquePoint](5)at(0.38,0.92){};
        \node[CliquePoint](6)at(0.92,0.38){};
        \node[CliquePoint](7)at(0.92,-0.38){};
        \node[CliquePoint](8)at(0.38,-0.92){};
        \draw[CliqueEdge](1)edge[]node[CliqueLabel]{$1$}(2);
        \draw[CliqueEdge](1)edge[bend left=30]node[CliqueLabel]{$-2$}(7);
        \draw[CliqueEmptyEdge](1)edge[]node[CliqueLabel]{}(8);
        \draw[CliqueEmptyEdge](2)edge[]node[CliqueLabel]{}(3);
        \draw[CliqueEdge](2)edge[bend right=30]node[CliqueLabel]{$1$}(4);
        \draw[CliqueEdge](2)edge[bend right=30]node[CliqueLabel]{$1$}(5);
        \draw[CliqueEmptyEdge](3)edge[]node[CliqueLabel]{}(4);
        \draw[CliqueEdge](4)edge[]node[CliqueLabel]{$2$}(5);
        \draw[CliqueEmptyEdge](5)edge[]node[CliqueLabel]{}(6);
        \draw[CliqueEdge](5)edge[bend right=30]node[CliqueLabel]{$1$}(7);
        \draw[CliqueEmptyEdge](6)edge[]node[CliqueLabel]{}(7);
        \draw[CliqueEmptyEdge](7)edge[]node[CliqueLabel]{}(8);
    \end{tikzpicture}\,.
\end{equation}
We call \Def{fundamental basis} of $\NC\Mca$ the fundamental basis of $\Cli\Mca$ restricted
to noncrossing $\Mca$-cliques. Observe that the fundamental basis of $\NC\Mca$ is a
set-operad basis.

To study $\NC\Mca$, we begin by establishing the fact that $\NC\Mca$ inherits some
properties of~$\Cli\Mca$. Then we shall describe a realization of $\NC\Mca$ in terms of
decorated Schröder trees, compute a minimal generating set of $\NC\Mca$, and compute its
dimensions.

\subsubsection{First properties}

\begin{Proposition} \label{prop:inherited_properties_NC_M}
    Let $\Mca$ be a unitary magma. Then,
    \begin{enumerate}[label={(\it\roman*)}]
        \item \label{item:inherited_properties_NC_M_1}
        the associative elements of $\NC\Mca$ are the ones of $\Cli\Mca$;
        \item \label{item:inherited_properties_NC_M_2}
        the group of symmetries of $\NC\Mca$ contains the map $\Returned$ (defined
        by~\eqref{equ:returned_map_Cli_M} in~\cite{Cliques1}) and all the maps $\Cli\theta$
        where $\theta$ are unitary magma automorphisms of $\Mca$;
        \item \label{item:inherited_properties_NC_M_3}
        the fundamental basis of $\NC\Mca$ is a basic set-operad basis if and only if $\Mca$
        is right cancelable;
        \item \label{item:inherited_properties_NC_M_4}
        the map $\rho$ (defined by~\eqref{equ:rotation_map_Cli_M} in~\cite{Cliques1}) is a
        rotation map of $\NC\Mca$ endowing it with a cyclic operad structure.
    \end{enumerate}
\end{Proposition}
\begin{proof}
    First, since $\NC\Mca$ is a suboperad of $\Cli\Mca$, each associative element of
    $\NC\Mca$ is an associative element of $\Cli\Mca$. Moreover, since all $\Mca$-bubbles
    are in $\NC\Mca$ and, as shown in~\cite{Cliques1}, all associative elements of
    $\Cli\Mca$ are linear combinations of $\Mca$-bubbles, each associative element of
    $\Cli\Mca$ belongs to $\NC\Mca$. This shows~\ref{item:inherited_properties_NC_M_1}.
    Moreover, since for any noncrossing $\Mca$-clique $\Pfr$, $\Returned(\Pfr)$ (resp.\
    $\rho(\Pfr)$) is still noncrossing and $\Returned$ belongs to the group of symmetries of
    $\Cli\Mca$ (resp.\ $\rho$ is a rotation map of $\Cli\Mca$),
    \ref{item:inherited_properties_NC_M_2} (resp.\ \ref{item:inherited_properties_NC_M_4})
    holds. Finally, again since $\NC\Mca$ is a suboperad of $\Cli\Mca$, since the
    fundamental basis of $\Cli\Mca$ is a basic set-operad basis, and since $\NC\Mca(2) =
    \Cli\Mca(2)$, \ref{item:inherited_properties_NC_M_3} holds.
\end{proof}

\subsubsection{Treelike expressions on bubbles} \label{subsubsec:treelike_bubbles}
Let $\Pfr$ be a noncrossing $\Mca$-clique of arity $n \geq 2$, and $(x, y)$ be a diagonal or
the base of $\Pfr$. Let $\{z_1, \dots, z_k\}$ be the set of vertices of $\Pfr$ such that $x
= z_1 < \dots < z_k = y$ and for any $i \in [k - 1]$, $z_{i + 1}$ is the greatest vertex of
$\Pfr$ such that $(z_i, z_{i + 1})$ is a solid diagonal or a (not necessarily solid) edge of
$\Pfr$. The \Def{area} of $\Pfr$ adjacent to $(x, y)$ is the $\Mca$-bubble $\Qfr$ of arity
$k$ whose base is labeled by $\Pfr(x, y)$ and $\Qfr_i = \Pfr(z_i, z_{i + 1})$ for all $i \in
[k]$. From a geometric point of view, $\Qfr$ is the unique maximal component of $\Pfr$
adjacent to the arc $(x, y)$, without solid diagonals, and bounded by solid diagonals or
edges of~$\Pfr$. For instance, for the noncrossing $\Z$-clique
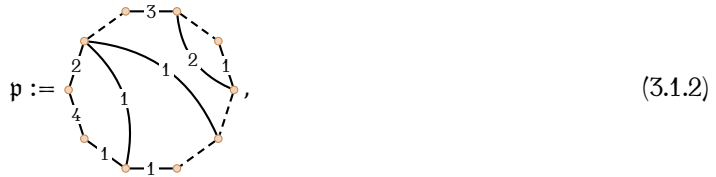
\begin{equation}
    \Pfr :=
    \begin{tikzpicture}[scale=1.1,Centering]
        \node[CliquePoint](1)at(-0.31,-0.95){};
        \node[CliquePoint](2)at(-0.81,-0.59){};
        \node[CliquePoint](3)at(-1.00,-0.00){};
        \node[CliquePoint](4)at(-0.81,0.59){};
        \node[CliquePoint](5)at(-0.31,0.95){};
        \node[CliquePoint](6)at(0.31,0.95){};
        \node[CliquePoint](7)at(0.81,0.59){};
        \node[CliquePoint](8)at(1.00,0.00){};
        \node[CliquePoint](9)at(0.81,-0.59){};
        \node[CliquePoint](10)at(0.31,-0.95){};
        \draw[CliqueEdge](1)edge[]node[CliqueLabel]{$1$}(2);
        \draw[CliqueEdge](1)edge[]node[CliqueLabel]{$1$}(10);
        \draw[CliqueEdge](2)edge[]node[CliqueLabel]{$4$}(3);
        \draw[CliqueEdge](1)edge[bend right=30]node[CliqueLabel]{$1$}(4);
        \draw[CliqueEdge](3)edge[]node[CliqueLabel]{$2$}(4);
        \draw[CliqueEmptyEdge](4)edge[]node[]{}(5);
        \draw[CliqueEdge](5)edge[]node[CliqueLabel]{$3$}(6);
        \draw[CliqueEdge](4)edge[bend left=30]node[CliqueLabel]{$1$}(9);
        \draw[CliqueEmptyEdge](6)edge[]node[]{}(7);
        \draw[CliqueEdge](6)edge[bend right=30]node[CliqueLabel]{$2$}(8);
        \draw[CliqueEdge](7)edge[]node[CliqueLabel]{$1$}(8);
        \draw[CliqueEmptyEdge](8)edge[]node[]{}(9);
        \draw[CliqueEmptyEdge](9)edge[]node[]{}(10);
    \end{tikzpicture}\,,
\end{equation}
the path associated with the diagonal $(4, 9)$ of $\Pfr$ is
$(4, 5, 6, 8, 9)$. For this reason, the area of $\Pfr$ adjacent to
$(4, 9)$ is the $\Z$-bubble
\begin{equation}
    \begin{tikzpicture}[scale=.6,Centering]
        \node[CliquePoint](1)at(-0.59,-0.81){};
        \node[CliquePoint](2)at(-0.95,0.31){};
        \node[CliquePoint](3)at(-0.00,1.00){};
        \node[CliquePoint](4)at(0.95,0.31){};
        \node[CliquePoint](5)at(0.59,-0.81){};
        \draw[CliqueEmptyEdge](1)edge[]node[]{}(2);
        \draw[CliqueEdge](1)edge[]node[CliqueLabel]{$1$}(5);
        \draw[CliqueEdge](2)edge[]node[CliqueLabel]{$3$}(3);
        \draw[CliqueEdge](3)edge[]node[CliqueLabel]{$2$}(4);
        \draw[CliqueEmptyEdge](4)edge[]node[]{}(5);
    \end{tikzpicture}\,.
\end{equation}

\begin{Proposition} \label{prop:unique_decomposition_NC_M}
    Let $\Mca$ be a unitary magma and $\Pfr$ be a noncrossing $\Mca$-clique of arity greater
    than $1$. Then there is a unique $\Mca$-bubble $\Qfr$ with a maximal arity $k \geq 2$
    such that $\Pfr = \Qfr \circ [\Rfr_1, \dots, \Rfr_k]$, where each $\Rfr_i$, $i \in [k]$,
    is a noncrossing $\Mca$-clique with a base labeled by~$\Unit_\Mca$.
\end{Proposition}
\begin{proof}
    Let $\Qfr'$ be the area of $\Pfr$ adjacent to its base and $k'$ be the arity of $\Qfr'$.
    By definition of the partial composition of $\NC\Mca$, for all $\Mca$-cliques $\Ufr$,
    $\Ufr'$, $\Ufr_1$, and $\Ufr_2$, if $\Ufr = \Ufr' \circ_i \Ufr_1 = \Ufr' \circ_i \Ufr_2$
    and $\Ufr_1$ and $\Ufr_2$ have bases labeled by $\Unit_\Mca$, then $\Ufr_1 = \Ufr_2$.
    This implies in particular that there are unique noncrossing $\Mca$-cliques $\Rfr'_i$,
    $i \in [k']$, with bases labeled by $\Unit_\Mca$ such that $\Pfr = \Qfr' \circ
    \left[\Rfr'_1, \dots, \Rfr'_{k'}\right]$. Finally, the fact that $\Qfr'$ is the area of
    $\Pfr$ adjacent to its base implies the maximality for the arity of $\Qfr'$. The
    statement of the proposition follows.
\end{proof}

Consider the linear map
\begin{equation}
    \BubbleTree :
    \NC\Mca \to \Free\left(\K \Angle{\Bubbles_\Mca}\right)
\end{equation}
defined recursively by $\BubbleTree(\UnitClique) := \Leaf$ and, for a noncrossing
$\Mca$-clique $\Pfr$ of arity greater than~$1$, by
\begin{equation}
    \BubbleTree(\Pfr) :=
    \Corolla(\Qfr) \circ
    \left[\BubbleTree(\Rfr_1), \dots, \BubbleTree(\Rfr_k)\right],
\end{equation}
where $\Pfr = \Qfr \circ [\Rfr_1, \dots, \Rfr_k]$ is the unique decomposition of $\Pfr$
stated in Proposition~\ref{prop:unique_decomposition_NC_M}. We call $\BubbleTree(\Pfr)$ the
\Def{bubble tree} of~$\Pfr$. For instance, in~$\NC\Z$,
\begin{equation} \label{equ:bubble_tree_example}
    \begin{tikzpicture}[scale=1.45,Centering]
        \node[CliquePoint](1)at(-0.31,-0.95){};
        \node[CliquePoint](2)at(-0.81,-0.59){};
        \node[CliquePoint](3)at(-1.00,-0.00){};
        \node[CliquePoint](4)at(-0.81,0.59){};
        \node[CliquePoint](5)at(-0.31,0.95){};
        \node[CliquePoint](6)at(0.31,0.95){};
        \node[CliquePoint](7)at(0.81,0.59){};
        \node[CliquePoint](8)at(1.00,0.00){};
        \node[CliquePoint](9)at(0.81,-0.59){};
        \node[CliquePoint](10)at(0.31,-0.95){};
        \draw[CliqueEdge](1)edge[]node[CliqueLabel]{$1$}(2);
        \draw[CliqueEdge](1)edge[bend left=20]node[CliqueLabel]{$2$}(5);
        \draw[CliqueEdge](1)edge[]node[CliqueLabel]{$1$}(10);
        \draw[CliqueEdge](2)edge[]node[CliqueLabel]{$4$}(3);
        \draw[CliqueEdge](2)edge[bend right=20]node[CliqueLabel]{$1$}(4);
        \draw[CliqueEdge](3)edge[]node[CliqueLabel]{$2$}(4);
        \draw[CliqueEmptyEdge](4)edge[]node[]{}(5);
        \draw[CliqueEdge](5)edge[]node[CliqueLabel]{$3$}(6);
        \draw[CliqueEdge](5)edge[bend right=30]node[CliqueLabel]{$3$}(9);
        \draw[CliqueEdge](5)edge[bend right=30]node[CliqueLabel]{$1$}(10);
        \draw[CliqueEmptyEdge](6)edge[]node[]{}(7);
        \draw[CliqueEdge](6)edge[bend right=30]node[CliqueLabel]{$2$}(9);
        \draw[CliqueEdge](7)edge[]node[CliqueLabel]{$1$}(8);
        \draw[CliqueEmptyEdge](8)edge[]node[]{}(9);
        \draw[CliqueEmptyEdge](9)edge[]node[]{}(10);
    \end{tikzpicture}
    \quad \xmapsto{\; \BubbleTree \;} \quad
    \begin{tikzpicture}[xscale=.4,yscale=.34,Centering]
        \node[](0)at(0.00,-6.00){};
        \node[](11)at(9.00,-12.00){};
        \node[](12)at(10.00,-12.00){};
        \node[](14)at(12.00,-6.00){};
        \node[](2)at(1.00,-9.00){};
        \node[](4)at(3.00,-9.00){};
        \node[](5)at(4.00,-6.00){};
        \node[](7)at(6.00,-9.00){};
        \node[](9)at(8.00,-12.00){};
        \node[](1)at(2.00,-3.00){
            \begin{tikzpicture}[scale=0.4]
                \node[CliquePoint](1_1)at(-0.71,-0.71){};
                \node[CliquePoint](1_2)at(-0.71,0.71){};
                \node[CliquePoint](1_3)at(0.71,0.71){};
                \node[CliquePoint](1_4)at(0.71,-0.71){};
                \draw[CliqueEdge](1_1)edge[]node[CliqueLabel]{$1$}(1_2);
                \draw[CliqueEmptyEdge](1_1)edge[]node[]{}(1_4);
                \draw[CliqueEdge](1_2)edge[]node[CliqueLabel]{$1$}(1_3);
                \draw[CliqueEmptyEdge](1_3)edge[]node[]{}(1_4);
            \end{tikzpicture}};
        \node[](10)at(9.00,-9.00){
            \begin{tikzpicture}[scale=0.4]
                \node[CliquePoint](10_1)at(-0.71,-0.71){};
                \node[CliquePoint](10_2)at(-0.71,0.71){};
                \node[CliquePoint](10_3)at(0.71,0.71){};
                \node[CliquePoint](10_4)at(0.71,-0.71){};
                \draw[CliqueEmptyEdge](10_1)edge[]node[]{}(10_2);
                \draw[CliqueEmptyEdge](10_1)edge[]node[]{}(10_4);
                \draw[CliqueEdge](10_2)edge[]node[CliqueLabel]{$1$}(10_3);
                \draw[CliqueEmptyEdge](10_3)edge[]node[]{}(10_4);
            \end{tikzpicture}};
        \node[](13)at(11.00,-3.00){
            \begin{tikzpicture}[scale=0.3]
                \node[CliquePoint](13_1)at(-0.87,-0.50){};
                \node[CliquePoint](13_2)at(-0.00,1.00){};
                \node[CliquePoint](13_3)at(0.87,-0.50){};
                \draw[CliqueEdge](13_1)edge[]node[CliqueLabel]{$3$}(13_2);
                \draw[CliqueEmptyEdge](13_1)edge[]node[]{}(13_3);
                \draw[CliqueEmptyEdge](13_2)edge[]node[]{}(13_3);
            \end{tikzpicture}};
        \node[](3)at(2.00,-6.00){
            \begin{tikzpicture}[scale=0.3]
                \node[CliquePoint](3_1)at(-0.87,-0.50){};
                \node[CliquePoint](3_2)at(-0.00,1.00){};
                \node[CliquePoint](3_3)at(0.87,-0.50){};
                \draw[CliqueEdge](3_1)edge[]node[CliqueLabel]{$4$}(3_2);
                \draw[CliqueEmptyEdge](3_1)edge[]node[]{}(3_3);
                \draw[CliqueEdge](3_2)edge[]node[CliqueLabel]{$2$}(3_3);
            \end{tikzpicture}};
        \node[](6)at(5.00,0.00){
            \begin{tikzpicture}[scale=0.3]
                \node[CliquePoint](6_1)at(-0.87,-0.50){};
                \node[CliquePoint](6_2)at(-0.00,1.00){};
                \node[CliquePoint](6_3)at(0.87,-0.50){};
                \draw[CliqueEdge](6_1)edge[]node[CliqueLabel]{$2$}(6_2);
                \draw[CliqueEdge](6_1)edge[]node[CliqueLabel]{$1$}(6_3);
                \draw[CliqueEdge](6_2)edge[]node[CliqueLabel]{$1$}(6_3);
            \end{tikzpicture}};
        \node[](8)at(7.00,-6.00){
            \begin{tikzpicture}[scale=0.3]
                \node[CliquePoint](8_1)at(-0.87,-0.50){};
                \node[CliquePoint](8_2)at(-0.00,1.00){};
                \node[CliquePoint](8_3)at(0.87,-0.50){};
                \draw[CliqueEdge](8_1)edge[]node[CliqueLabel]{$3$}(8_2);
                \draw[CliqueEmptyEdge](8_1)edge[]node[]{}(8_3);
                \draw[CliqueEdge](8_2)edge[]node[CliqueLabel]{$2$}(8_3);
            \end{tikzpicture}};
        \draw[Edge](0)--(1);
        \draw[Edge](1)--(6);
        \draw[Edge](10)--(8);
        \draw[Edge](11)--(10);
        \draw[Edge](12)--(10);
        \draw[Edge](13)--(6);
        \draw[Edge](14)--(13);
        \draw[Edge](2)--(3);
        \draw[Edge](3)--(1);
        \draw[Edge](4)--(3);
        \draw[Edge](5)--(1);
        \draw[Edge](7)--(8);
        \draw[Edge](8)--(13);
        \draw[Edge](9)--(10);
        \node(r)at(5.00,2.5){};
        \draw[Edge](r)--(6);
    \end{tikzpicture}\,.
\end{equation}

\begin{Lemma} \label{lem:map_NC_M_bubble_tree_treelike_expression}
    Let $\Mca$ be a unitary magma. For a noncrossing $\Mca$-clique $\Pfr$,
    $\BubbleTree(\Pfr)$ is a treelike expression on $\Bubbles_\Mca$ of~$\Pfr$.
\end{Lemma}
\begin{proof}
    We proceed by induction on the arity $n$ of $\Pfr$. If $n = 1$, since $\Pfr =
    \UnitClique$ and $\BubbleTree(\UnitClique) = \Leaf$, the statement of the lemma
    immediately follows. Otherwise, we have
    \begin{math}
        \BubbleTree(\Pfr) =
        \Corolla(\Qfr) \circ
        \left[\BubbleTree(\Rfr_1), \dots, \BubbleTree(\Rfr_k)\right]
    \end{math}
    where $\Pfr$ uniquely decomposes as $\Pfr = \Qfr \circ [\Rfr_1, \dots, \Rfr_k]$ under
    the conditions stated by  Proposition~\ref{prop:unique_decomposition_NC_M}. By
    definition of area and of the map $\BubbleTree$, $\Qfr$ is an $\Mca$-bubble. Moreover,
    by induction hypothesis, any $\BubbleTree(\Rfr_i)$, $i \in [k]$, is a treelike
    expression on $\Bubbles_\Mca$ of $\Rfr_i$. Hence, $\BubbleTree(\Pfr)$ is a treelike
    expression on $\Bubbles_\Mca$ of~$\Pfr$.
\end{proof}

\begin{Proposition} \label{prop:map_NC_M_bubble_tree}
    Let $\Mca$ be a unitary magma. Then the map $\BubbleTree$ is injective and the image of
    $\BubbleTree$ is the linear span of all syntax trees $\Tfr$ on $\Bubbles_\Mca$ such that
    \begin{enumerate}[label={(\it\roman*)}]
        \item \label{item:map_NC_M_bubble_tree_1}
        the root of $\Tfr$ is labeled by an $\Mca$-bubble;
        \item \label{item:map_NC_M_bubble_tree_2}
        the internal nodes of $\Tfr$ different from the root are labeled by $\Mca$-bubbles
        whose bases are labeled by~$\Unit_\Mca$;
        \item \label{item:map_NC_M_bubble_tree_3}
        if $x$ and $y$ are two internal nodes of $\Tfr$ such that $y$ is the $i$th child of
        $x$, the $i$th edge of the bubble labeling $x$ is solid.
    \end{enumerate}
\end{Proposition}
\begin{proof}
    First of all, since by definition $\BubbleTree$ sends a basis element of $\NC\Mca$ to a
    basis element of $\Free\left(\K \Angle{\Bubbles_\Mca}\right)$, it is sufficient to show
    that $\BubbleTree$ is injective as a map from $\Cliques_\Mca$ to the set of syntax trees
    on $\Bubbles_\Mca$ to establish that it is an injective linear map. For this, we proceed
    by induction on the arity $n$. If $n = 1$, since $\BubbleTree(\UnitClique) = \Leaf$ and
    $\NC\Mca(1)$ is of dimension $1$, $\BubbleTree$ is injective. Assume now that $\Pfr$ and
    $\Pfr'$ are two noncrossing $\Mca$-cliques of arity $n$ such that $\BubbleTree(\Pfr) =
    \BubbleTree(\Pfr')$. Hence, $\Pfr$ (resp.\ $\Pfr'$) uniquely decomposes as $\Pfr = \Qfr
    \circ [\Rfr_1, \dots, \Rfr_k]$ (resp.\ $\Pfr' = \Qfr' \circ [\Rfr'_1, \dots, \Rfr'_k]$)
    as stated by Proposition~\ref{prop:unique_decomposition_NC_M} and
    \begin{equation}
        \BubbleTree(\Pfr)
        = \Corolla(\Qfr) \circ
        \left[\BubbleTree(\Rfr_1), \dots, \BubbleTree(\Rfr_k)\right]
        = \Corolla(\Qfr') \circ
        \left[\BubbleTree(\Rfr'_1), \dots, \BubbleTree(\Rfr'_k)\right]
        = \BubbleTree(\Pfr').
    \end{equation}
    Now, because by definition of area, all bases of the $\Rfr_i$ and $\Rfr'_i$, $i \in
    [k]$, are labeled by $\Unit_\Mca$, this implies that $\Qfr = \Qfr'$. Therefore, we have
    $\BubbleTree(\Rfr_i) = \BubbleTree(\Rfr'_i)$ for all $i \in [k]$, so that, by induction
    hypothesis, $\Rfr_i = \Rfr'_i$ for all $i \in [k]$. Hence, $\BubbleTree$ is injective.

    The definition of $\BubbleTree$ together with
    Proposition~\ref{prop:unique_decomposition_NC_M} leads to the fact that, for a
    noncrossing $\Mca$-clique $\Pfr$, the syntax tree $\BubbleTree(\Pfr)$
    satisfies~\ref{item:map_NC_M_bubble_tree_1}, \ref{item:map_NC_M_bubble_tree_2},
    and~\ref{item:map_NC_M_bubble_tree_3}. Conversely, let $\Tfr$ be a syntax tree
    satisfying~\ref{item:map_NC_M_bubble_tree_1}, \ref{item:map_NC_M_bubble_tree_2},
    and~\ref{item:map_NC_M_bubble_tree_3}. We show by structural induction on $\Tfr$
    that there is a noncrossing $\Mca$-clique $\Pfr$ such that $\BubbleTree(\Pfr) = \Tfr$.
    If $\Tfr = \Leaf$, the property holds because $\BubbleTree(\UnitClique) = \Leaf$.
    Otherwise, we have $\Tfr = \Sfr \circ [\Ufr_1, \dots, \Ufr_k]$ where $\Sfr$ is a syntax
    tree of degree $1$ and the $\Ufr_i$, $i \in [k]$, are syntax trees. Since $\Tfr$
    satisfies~\ref{item:map_NC_M_bubble_tree_1}, \ref{item:map_NC_M_bubble_tree_2},
    and~\ref{item:map_NC_M_bubble_tree_3}, the trees $\Sfr$ and $\Ufr_i$, $i \in [k]$,
    satisfy the same three properties. Therefore, by induction hypothesis, there are
    noncrossing $\Mca$-cliques $\Qfr$ and $\Rfr_i$, $i \in [k]$, such that
    $\BubbleTree(\Qfr) = \Sfr$ and $\BubbleTree(\Rfr_i) = \Ufr_i$. Now define $\Pfr$ as the
    noncrossing $\Mca$-clique $\Qfr \circ [\Rfr_1, \dots, \Rfr_k]$. By definition of the map
    $\BubbleTree$ and the unique decomposition stated in
    Proposition~\ref{prop:unique_decomposition_NC_M} for $\Pfr$, one obtains that
    $\BubbleTree(\Pfr) = \Tfr$.
\end{proof}

Observe that $\BubbleTree$ is not an operad morphism. Indeed,
\begin{equation} \label{equ:bubble_tree_not_morphism}
    \BubbleTree\left(
        \begin{tikzpicture}[scale=0.3,Centering]
        \node[CliquePoint](1)at(-0.87,-0.50){};
        \node[CliquePoint](2)at(-0.00,1.00){};
        \node[CliquePoint](3)at(0.87,-0.50){};
        \draw[CliqueEmptyEdge](1)edge[]node[]{}(2);
        \draw[CliqueEmptyEdge](1)edge[]node[]{}(3);
        \draw[CliqueEmptyEdge](2)edge[]node[]{}(3);
    \end{tikzpicture}
    \circ_1
        \begin{tikzpicture}[scale=0.3,Centering]
        \node[CliquePoint](1)at(-0.87,-0.50){};
        \node[CliquePoint](2)at(-0.00,1.00){};
        \node[CliquePoint](3)at(0.87,-0.50){};
        \draw[CliqueEmptyEdge](1)edge[]node[]{}(2);
        \draw[CliqueEmptyEdge](1)edge[]node[]{}(3);
        \draw[CliqueEmptyEdge](2)edge[]node[]{}(3);
    \end{tikzpicture}
    \right)
    =
    \begin{tikzpicture}[scale=0.4,Centering]
        \node[](0)at(0.00,-2.00){};
        \node[](2)at(1.00,-2.00){};
        \node[](3)at(2.00,-2.00){};
        \node[](1)at(1.00,0.00){
            \begin{tikzpicture}[scale=0.4]
                \node[CliquePoint](1_1)at(-0.71,-0.71){};
                \node[CliquePoint](1_2)at(-0.71,0.71){};
                \node[CliquePoint](1_3)at(0.71,0.71){};
                \node[CliquePoint](1_4)at(0.71,-0.71){};
                \draw[CliqueEmptyEdge](1_1)edge[]node[]{}(1_2);
                \draw[CliqueEmptyEdge](1_1)edge[]node[]{}(1_4);
                \draw[CliqueEmptyEdge](1_2)edge[]node[]{}(1_3);
                \draw[CliqueEmptyEdge](1_3)edge[]node[]{}(1_4);
            \end{tikzpicture}};
        \draw[Edge](0)--(1);
        \draw[Edge](2)--(1);
        \draw[Edge](3)--(1);
        \node(r)at(1.00,1.50){};
        \draw[Edge](r)--(1);
    \end{tikzpicture}
    \enspace \ne
    \begin{tikzpicture}[scale=0.4,Centering]
        \node[](0)at(0.00,-3.33){};
        \node[](2)at(2.00,-3.33){};
        \node[](4)at(4.00,-1.67){};
        \node[](1)at(1.00,-1.67){
            \begin{tikzpicture}[scale=0.3]
                \node[CliquePoint](1_1)at(-0.87,-0.50){};
                \node[CliquePoint](1_2)at(-0.00,1.00){};
                \node[CliquePoint](1_3)at(0.87,-0.50){};
                \draw[CliqueEmptyEdge](1_1)edge[]node[]{}(1_2);
                \draw[CliqueEmptyEdge](1_1)edge[]node[]{}(1_3);
                \draw[CliqueEmptyEdge](1_2)edge[]node[]{}(1_3);
            \end{tikzpicture}};
        \node[](3)at(3.00,0.00){
            \begin{tikzpicture}[scale=0.3]
                \node[CliquePoint](3_1)at(-0.87,-0.50){};
                \node[CliquePoint](3_2)at(-0.00,1.00){};
                \node[CliquePoint](3_3)at(0.87,-0.50){};
                \draw[CliqueEmptyEdge](3_1)edge[]node[]{}(3_2);
                \draw[CliqueEmptyEdge](3_1)edge[]node[]{}(3_3);
                \draw[CliqueEmptyEdge](3_2)edge[]node[]{}(3_3);
            \end{tikzpicture}};
        \draw[Edge](0)--(1);
        \draw[Edge](1)--(3);
        \draw[Edge](2)--(1);
        \draw[Edge](4)--(3);
        \node(r)at(3.00,1.25){};
        \draw[Edge](r)--(3);
    \end{tikzpicture}
    =
    \BubbleTree\left(
    \begin{tikzpicture}[scale=0.3,Centering]
        \node[CliquePoint](1)at(-0.87,-0.50){};
        \node[CliquePoint](2)at(-0.00,1.00){};
        \node[CliquePoint](3)at(0.87,-0.50){};
        \draw[CliqueEmptyEdge](1)edge[]node[]{}(2);
        \draw[CliqueEmptyEdge](1)edge[]node[]{}(3);
        \draw[CliqueEmptyEdge](2)edge[]node[]{}(3);
    \end{tikzpicture}
    \right)
    \circ_1
    \BubbleTree\left(
    \begin{tikzpicture}[scale=0.3,Centering]
        \node[CliquePoint](1)at(-0.87,-0.50){};
        \node[CliquePoint](2)at(-0.00,1.00){};
        \node[CliquePoint](3)at(0.87,-0.50){};
        \draw[CliqueEmptyEdge](1)edge[]node[]{}(2);
        \draw[CliqueEmptyEdge](1)edge[]node[]{}(3);
        \draw[CliqueEmptyEdge](2)edge[]node[]{}(3);
    \end{tikzpicture}
    \right).
\end{equation}
Observe that~\eqref{equ:bubble_tree_not_morphism} holds for all unitary magmas $\Mca$ since
$\Unit_\Mca$ is always idempotent.

\subsubsection[Realization]{Realization in terms of decorated Schröder trees}%
\label{subsubsec:M_Schroder_trees}
Recall that a \Def{Schröder tree} is a tree such that all internal nodes have at least two
children. An \Def{$\Mca$-Schröder tree} $\Tfr$ is a Schröder tree such that each edge
connecting two internal nodes is labeled by $\bar{\Mca}$, each edge connecting an internal
node and a leaf is labeled by $\Mca$, and the outgoing edge from the root of $\Tfr$ is
labeled by $\Mca$ (see~\eqref{equ:example_M_Schroder_tree} for an example of a $\Z$-Schröder
tree).

From the description of the image of the map $\BubbleTree$ provided by
Proposition~\ref{prop:map_NC_M_bubble_tree}, any bubble tree $\Tfr$ of a noncrossing
$\Mca$-clique $\Pfr$ of arity $n$ can be encoded by an $\Mca$-Schröder tree $\Sfr$ with $n$
leaves. Indeed, this $\Mca$-Schröder tree is obtained by considering each internal node $x$
of $\Tfr$ and by labeling the edge connecting $x$ and its $i$th child by the label of the
$i$th edge of the $\Mca$-bubble labeling $x$. The outgoing edge from the root of $\Sfr$ is
labeled by the label of the base of the $\Mca$-bubble labeling the root of $\Tfr$. For
instance, the bubble tree of~\eqref{equ:bubble_tree_example} is encoded by the $\Z$-Schröder
tree
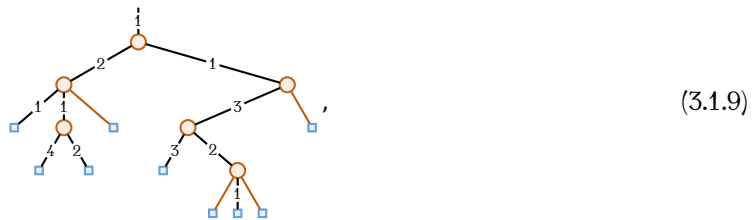
\begin{equation} \label{equ:example_M_Schroder_tree}
    \begin{tikzpicture}[xscale=.33,yscale=.19,Centering]
        \node[Leaf](0)at(0.00,-6.00){};
        \node[Leaf](11)at(9.00,-12.00){};
        \node[Leaf](12)at(10.00,-12.00){};
        \node[Leaf](14)at(12.00,-6.00){};
        \node[Leaf](2)at(1.00,-9.00){};
        \node[Leaf](4)at(3.00,-9.00){};
        \node[Leaf](5)at(4.00,-6.00){};
        \node[Leaf](7)at(6.00,-9.00){};
        \node[Leaf](9)at(8.00,-12.00){};
        \node[Node](1)at(2.00,-3.00){};
        \node[Node](10)at(9.00,-9.00){};
        \node[Node](13)at(11.00,-3.00){};
        \node[Node](3)at(2.00,-6.00){};
        \node[Node](6)at(5.00,0.00){};
        \node[Node](8)at(7.00,-6.00){};
        \draw[Edge](0)edge[]node[EdgeLabel]{$1$}(1);
        \draw[Edge](1)edge[]node[EdgeLabel]{$2$}(6);
        \draw[Edge](10)edge[]node[EdgeLabel]{$2$}(8);
        \draw[Edge](11)edge[]node[EdgeLabel]{$1$}(10);
        \draw[Edge](12)--(10);
        \draw[Edge](13)edge[]node[EdgeLabel]{$1$}(6);
        \draw[Edge](14)--(13);
        \draw[Edge](2)edge[]node[EdgeLabel]{$4$}(3);
        \draw[Edge](3)edge[]node[EdgeLabel]{$1$}(1);
        \draw[Edge](4)edge[]node[EdgeLabel]{$2$}(3);
        \draw[Edge](5)--(1);
        \draw[Edge](7)edge[]node[EdgeLabel]{$3$}(8);
        \draw[Edge](8)edge[]node[EdgeLabel]{$3$}(13);
        \draw[Edge](9)--(10);
        \node(r)at(5.00,3){};
        \draw[Edge](r)edge[]node[EdgeLabel]{$1$}(6);
    \end{tikzpicture}\,,
\end{equation}
where the labels of the edges are placed in their centers and where unlabeled edges are
implicitly labeled by $\Unit_\Mca$. We shall use these drawing conventions in the sequel. As
a side remark, observe that the $\Mca$-Schröder tree encoding a noncrossing $\Mca$-clique
$\Pfr$ and the dual tree of $\Pfr$ (in the usual meaning) have the same underlying unlabeled
tree.

This encoding of noncrossing $\Mca$-cliques by bubble trees is reversible, and hence one can
interpret $\NC\Mca$ as an operad on the linear span of all $\Mca$-Schröder trees. Hence,
through this interpretation, if $\Sfr$ and $\Tfr$ are two $\Mca$-Schröder trees and $i$ is a
valid integer, the tree $\Sfr \circ_i \Tfr$ is computed by grafting the root of $\Tfr$ to
the $i$th leaf of $\Sfr$. Then, by writing $b$ for the label of the edge adjacent to the
root of $\Tfr$ and $a$ for the label of the edge adjacent to the $i$th leaf of $\Sfr$, we
have two cases to consider, depending on the value of $c := a \Op b$. If $c \ne \Unit_\Mca$,
we label the edge connecting $\Sfr$ and $\Tfr$ by $c$. Otherwise, if $c = \Unit_\Mca$, we
contract the edge connecting $\Sfr$ and $\Tfr$ by merging the root of $\Tfr$ and the direct
ancestor of the $i$th leaf of $\Sfr$ (see Figure~\ref{fig:composition_NC_M_Schroder_trees}).
\begin{figure}[ht]
    \centering
    \captionsetup[subfigure]{width=5cm}
    \subfloat[]
    [The expression $\Sfr \circ_i \Tfr$ to compute. The displayed
    leaf is the $i$th one of~$\Sfr$.]{
        \begin{tikzpicture}[xscale=1,yscale=.55,Centering]
            \node[Subtree](1)at(0,0){$\Sfr'$};
            \node[Node](2)at(0,-1){};
            \node[Leaf](3)at(0,-2){};
            \draw[Edge](1)--(2);
            \draw[Edge](2)edge[]node[EdgeLabel]{$a$}(3);
            \node[below of=2,font=\footnotesize]{$i$};
        \end{tikzpicture}
        \quad $\circ_i$ \quad
        \begin{tikzpicture}[xscale=.9,yscale=.8,Centering]
            \node[Node](1)at(0,0){};
            \node(r)at(0,1){};
            \node[Subtree](2)at(-1,-1){$\Tfr_1$};
            \node[Subtree](3)at(1,-1){$\Tfr_k$};
            \draw[Edge](1)edge[]node[EdgeLabel]{$b$}(r);
            \draw[Edge](1)--(2);
            \draw[Edge](1)--(3);
            \node at(0,-1){$\dots$};
        \end{tikzpicture}
    \label{subfig:composition_NC_M_Schroder_trees_1}}

    \captionsetup[subfigure]{width=3cm}
    \subfloat[]
    [The resulting tree if $a \Op b \ne \Unit_\Mca$.]{
        \begin{tikzpicture}[xscale=.9,yscale=.68,Centering]
            \node[Subtree](1)at(0,-.5){$\Sfr'$};
            \node[Node](2)at(0,-1.5){};
            \draw[Edge](1)--(2);
            \node[Node](4)at(0,-3){};
            \node[Subtree](6)at(-1,-4){$\Tfr_1$};
            \node[Subtree](7)at(1,-4){$\Tfr_k$};
            \draw[Edge](4)--(6);
            \draw[Edge](4)--(7);
            \node at(0,-4){$\dots$};
            \draw[Edge](4)edge[]node[EdgeLabel]{$a \Op b$}(2);
        \end{tikzpicture}
    \label{subfig:composition_NC_M_Schroder_trees_2}}
    \qquad
    \qquad
    \qquad
    \subfloat[]
    [The resulting tree if $a \Op b = \Unit_\Mca$.]{
        \begin{tikzpicture}[xscale=.9,yscale=.7,Centering]
            \node[Subtree](1)at(0,-2){$\Sfr'$};
            \node[Node](4)at(0,-3){};
            \draw[Edge](1)--(4);
            \node[Subtree](6)at(-1,-4){$\Tfr_1$};
            \node[Subtree](7)at(1,-4){$\Tfr_k$};
            \draw[Edge](4)--(6);
            \draw[Edge](4)--(7);
            \node at(0,-4){$\dots$};
        \end{tikzpicture}
    \label{subfig:composition_NC_M_Schroder_trees_3}}
    \caption{
    The partial composition of $\NC\Mca$ realized on $\Mca$-Schröder trees. Here, the two
    cases~\protect\subref{subfig:composition_NC_M_Schroder_trees_2}
    and~\protect\subref{subfig:composition_NC_M_Schroder_trees_3} for the computation of
    $\Sfr \circ_i \Tfr$ are shown, where $\Sfr$ and $\Tfr$ are two $\Mca$-Schröder trees. In
    these drawings, the triangles denote subtrees.}
    \label{fig:composition_NC_M_Schroder_trees}
\end{figure}
For instance, in $\NC\N_3$, we have the two partial compositions
\begin{subequations}
\begin{equation}
    \begin{tikzpicture}[xscale=.28,yscale=.2,Centering]
        \node[Leaf](0)at(0.00,-5.33){};
        \node[Leaf](2)at(2.00,-5.33){};
        \node[Leaf](4)at(3.00,-5.33){};
        \node[Leaf](6)at(5.00,-5.33){};
        \node[Leaf](7)at(6.00,-2.67){};
        \node[Node](1)at(1.00,-2.67){};
        \node[Node](3)at(4.00,0.00){};
        \node[Node](5)at(4.00,-2.67){};
        \draw[Edge](0)--(1);
        \draw[Edge](1)edge[]node[EdgeLabel]{$1$}(3);
        \draw[Edge](2)edge[]node[EdgeLabel]{$1$}(1);
        \draw[Edge](4)edge[]node[EdgeLabel]{$2$}(5);
        \draw[Edge](5)edge[]node[EdgeLabel]{$1$}(3);
        \draw[Edge](6)--(5);
        \draw[Edge](7)--(3);
        \node(r)at(4.00,2.75){};
        \draw[Edge](r)edge[]node[EdgeLabel]{$2$}(3);
    \end{tikzpicture}
    \circ_2
    \begin{tikzpicture}[xscale=.3,yscale=.27,Centering]
        \node[Leaf](0)at(0.00,-1.67){};
        \node[Leaf](2)at(2.00,-3.33){};
        \node[Leaf](4)at(4.00,-3.33){};
        \node[Node](1)at(1.00,0.00){};
        \node[Node](3)at(3.00,-1.67){};
        \draw[Edge](0)--(1);
        \draw[Edge](2)--(3);
        \draw[Edge](3)edge[]node[EdgeLabel]{$1$}(1);
        \draw[Edge](4)edge[]node[EdgeLabel]{$2$}(3);
        \node(r)at(1.00,2){};
        \draw[Edge](r)edge[]node[EdgeLabel]{$1$}(1);
    \end{tikzpicture}
    =
    \begin{tikzpicture}[xscale=.28,yscale=.20,Centering]
        \node[Leaf](0)at(0.00,-4.80){};
        \node[Leaf](10)at(9.00,-4.80){};
        \node[Leaf](11)at(10.00,-2.40){};
        \node[Leaf](2)at(2.00,-7.20){};
        \node[Leaf](4)at(4.00,-9.60){};
        \node[Leaf](6)at(6.00,-9.60){};
        \node[Leaf](8)at(7.00,-4.80){};
        \node[Node](1)at(1.00,-2.40){};
        \node[Node](3)at(3.00,-4.80){};
        \node[Node](5)at(5.00,-7.20){};
        \node[Node](7)at(8.00,0.00){};
        \node[Node](9)at(8.00,-2.40){};
        \draw[Edge](0)--(1);
        \draw[Edge](1)edge[]node[EdgeLabel]{$1$}(7);
        \draw[Edge](10)--(9);
        \draw[Edge](11)--(7);
        \draw[Edge](2)--(3);
        \draw[Edge](3)edge[]node[EdgeLabel]{$2$}(1);
        \draw[Edge](4)--(5);
        \draw[Edge](5)edge[]node[EdgeLabel]{$1$}(3);
        \draw[Edge](6)edge[]node[EdgeLabel]{$2$}(5);
        \draw[Edge](8)edge[]node[EdgeLabel]{$2$}(9);
        \draw[Edge](9)edge[]node[EdgeLabel]{$1$}(7);
        \node(r)at(8.00,2.5){};
        \draw[Edge](r)edge[]node[EdgeLabel]{$2$}(7);
    \end{tikzpicture}\,,
\end{equation}
\begin{equation}
    \begin{tikzpicture}[xscale=.28,yscale=.2,Centering]
        \node[Leaf](0)at(0.00,-5.33){};
        \node[Leaf](2)at(2.00,-5.33){};
        \node[Leaf](4)at(3.00,-5.33){};
        \node[Leaf](6)at(5.00,-5.33){};
        \node[Leaf](7)at(6.00,-2.67){};
        \node[Node](1)at(1.00,-2.67){};
        \node[Node](3)at(4.00,0.00){};
        \node[Node](5)at(4.00,-2.67){};
        \draw[Edge](0)--(1);
        \draw[Edge](1)edge[]node[EdgeLabel]{$1$}(3);
        \draw[Edge](2)edge[]node[EdgeLabel]{$1$}(1);
        \draw[Edge](4)edge[]node[EdgeLabel]{$2$}(5);
        \draw[Edge](5)edge[]node[EdgeLabel]{$1$}(3);
        \draw[Edge](6)--(5);
        \draw[Edge](7)--(3);
        \node(r)at(4.00,2.75){};
        \draw[Edge](r)edge[]node[EdgeLabel]{$2$}(3);
    \end{tikzpicture}
    \circ_3
    \begin{tikzpicture}[xscale=.3,yscale=.27,Centering]
        \node[Leaf](0)at(0.00,-1.67){};
        \node[Leaf](2)at(2.00,-3.33){};
        \node[Leaf](4)at(4.00,-3.33){};
        \node[Node](1)at(1.00,0.00){};
        \node[Node](3)at(3.00,-1.67){};
        \draw[Edge](0)--(1);
        \draw[Edge](2)--(3);
        \draw[Edge](3)edge[]node[EdgeLabel]{$1$}(1);
        \draw[Edge](4)edge[]node[EdgeLabel]{$2$}(3);
        \node(r)at(1.00,2){};
        \draw[Edge](r)edge[]node[EdgeLabel]{$1$}(1);
    \end{tikzpicture}
    =
    \begin{tikzpicture}[xscale=.28,yscale=.20,Centering]
        \node[Leaf](0)at(0.00,-5.50){};
        \node[Leaf](10)at(8.00,-2.75){};
        \node[Leaf](2)at(2.00,-5.50){};
        \node[Leaf](4)at(3.00,-5.50){};
        \node[Leaf](6)at(4.00,-8.25){};
        \node[Leaf](8)at(6.00,-8.25){};
        \node[Leaf](9)at(7.00,-5.50){};
        \node[Node](1)at(1.00,-2.75){};
        \node[Node](3)at(5.00,0.00){};
        \node[Node](5)at(5.00,-2.75){};
        \node[Node](7)at(5.00,-5.50){};
        \draw[Edge](0)--(1);
        \draw[Edge](1)edge[]node[EdgeLabel]{$1$}(3);
        \draw[Edge](10)--(3);
        \draw[Edge](2)edge[]node[EdgeLabel]{$1$}(1);
        \draw[Edge](4)--(5);
        \draw[Edge](5)edge[]node[EdgeLabel]{$1$}(3);
        \draw[Edge](6)--(7);
        \draw[Edge](7)edge[]node[EdgeLabel]{$1$}(5);
        \draw[Edge](8)edge[]node[EdgeLabel]{$2$}(7);
        \draw[Edge](9)--(5);
        \node(r)at(5.00,2.5){};
        \draw[Edge](r)edge[]node[EdgeLabel]{$2$}(3);
    \end{tikzpicture}\,.
\end{equation}
\end{subequations}

In the sequel, we shall indifferently view $\NC\Mca$ as an operad on noncrossing
$\Mca$-cliques or on $\Mca$-Schröder trees.

\subsubsection{Minimal generating set}

\begin{Proposition} \label{prop:generating_set_NC_M}
    Let $\Mca$ be a unitary magma. The set $\Triangles_\Mca$ of all $\Mca$-triangles is a
    minimal generating set of~$\NC\Mca$.
\end{Proposition}
\begin{proof}
    We start by showing by induction on the arity that the suboperad
    $(\NC\Mca)^{\Triangles_\Mca}$ of $\NC\Mca$ generated by $\Triangles_\Mca$ is $\NC\Mca$.
    This is immediately true in arity $1$. Let $\Pfr$ be a noncrossing $\Mca$-clique of
    arity $n \geq 2$. Proposition~\ref{prop:unique_decomposition_NC_M} says in particular
    that we can express $\Pfr$ as $\Pfr = \Qfr \circ [\Rfr_1, \dots, \Rfr_k]$ where $\Qfr$
    is an $\Mca$-bubble of arity $k \geq 2$ and the $\Rfr_i$, $i \in [k]$, are noncrossing
    $\Mca$-cliques. Since $\Qfr$ is an $\Mca$-bubble, it can be expressed as
    \begin{equation}\label{equ:generating_set_NC_M}
        \Qfr =
        \TriangleXEX{\Qfr_0}{\Unit_\Mca}{\Qfr_k}
        \circ_1
        \TriangleEEX{\Unit_\Mca}{\Unit_\Mca}{\Qfr_{k - 1}}
        \circ_1 \dots \circ_1
        \TriangleEEX{\Unit_\Mca}{\Unit_\Mca}{\Qfr_3}
        \circ_1
        \TriangleEXX{\Unit_\Mca}{\Qfr_1}{\Qfr_2}\,.
    \end{equation}
    Observe that, in~\eqref{equ:generating_set_NC_M}, brackets are not necessary since
    $\circ_1$ is associative. Since $k \geq 2$, the arities of each $\Rfr_i$, $i \in [k]$,
    are smaller than the one of $\Pfr$. For this reason, by induction hypothesis, each
    $\Rfr_i$ belongs to $(\NC\Mca)^{\Triangles_\Mca}$. Moreover,
    since~\eqref{equ:generating_set_NC_M} shows an expression of $\Qfr$ by partial
    compositions of $\Mca$-triangles, $\Qfr$ also belongs to $(\NC\Mca)^{\Triangles_\Mca}$.
    This implies that this is also the case for $\Pfr$. Hence, $\NC\Mca$ is generated
    by~$\Triangles_\Mca$.

    Finally, due to the fact that the partial composition of two $\Mca$-triangles is an
    $\Mca$-clique of arity $3$, if $\Pfr$ is an $\Mca$-triangle, $\Pfr$ cannot be expressed
    as a partial composition of $\Mca$-triangles. Moreover, since the space $\NC\Mca(1)$ is
    trivial, these arguments imply that $\Triangles_\Mca$ is a minimal generating set
    of~$\NC\Mca$.
\end{proof}

Proposition~\ref{prop:generating_set_NC_M} also says that $\NC\Mca$ is the smallest
suboperad of $\Cli\Mca$ that contains all $\Mca$-triangles and that $\NC\Mca$ is the biggest
binary suboperad of~$\Cli\Mca$.

\subsubsection{Dimensions}
We now use the notion of bubble trees introduced in Section~\ref{subsubsec:treelike_bubbles}
to compute the dimensions of~$\NC\Mca$.

\begin{Proposition} \label{prop:Hilbert_series_NC_M}
    Let $\Mca$ be a finite unitary magma. The Hilbert series $\Hilbert_{\NC\Mca}(t)$ of
    $\NC\Mca$ satisfies
    \begin{equation} \label{equ:Hilbert_series_NC_M}
        t + \left(m^3 - 2m^2 + 2m - 1\right)t^2
        + \left(2m^2t - 3mt + 2t - 1\right) \Hilbert_{\NC\Mca}(t)
        + \left(m - 1\right) \Hilbert_{\NC\Mca}(t)^2
        = 0,
    \end{equation}
    where $m := \# \Mca$.
\end{Proposition}
\begin{proof}
    By Proposition~\ref{prop:map_NC_M_bubble_tree}, the set of noncrossing $\Mca$-cliques is
    in one-to-one correspondence with the set of syntax trees on $\Bubbles_\Mca$ that
    satisfy~\ref{item:map_NC_M_bubble_tree_1}, \ref{item:map_NC_M_bubble_tree_2},
    and~\ref{item:map_NC_M_bubble_tree_3}. We call $T(t)$ the generating series of these
    trees and $S(t)$ the generating series of these trees with the extra condition that the
    roots are labeled by $\Mca$-bubbles whose bases are labeled by~$\Unit_\Mca$. Immediately
    from its description, $S(t)$ satisfies
    \begin{equation} \label{equ:Hilbert_series_NC_M_1}
        S(t) = t + \sum_{n \geq 2} \left((m - 1) S(t) + t\right)^n,
    \end{equation}
    and $T(t)$ satisfies
    \begin{equation} \label{equ:Hilbert_series_NC_M_2}
        T(t) = t + m (S(t) - t).
    \end{equation}
    As the set of noncrossing $\Mca$-cliques forms the fundamental basis of $\NC\Mca$, we
    have $\Hilbert_{\NC\Mca}(t) = T(t)$. We eventually
    obtain~\eqref{equ:Hilbert_series_NC_M} from~\eqref{equ:Hilbert_series_NC_M_1}
    and~\eqref{equ:Hilbert_series_NC_M_2} by a direct computation.
\end{proof}

From Proposition~\ref{prop:Hilbert_series_NC_M}, we deduce that the Hilbert series of
$\NC\Mca$ satisfies
\begin{equation} \label{equ:Hilbert_series_NC_M_function}
    \Hilbert_{\NC\Mca}(t) =
    \frac{1 - (2m^2 - 3m + 2)t - \sqrt{1 - 2(2m^2 - m)t + m^2t^2}}{2(m - 1)},
\end{equation}
where $m := \# \Mca \ne 1$.

By using Narayana numbers, whose definition is recalled in
Section~\ref{subsubsec:quotient_Cli_M_Inf} of~\cite{Cliques1}, we can state the following
result.

\begin{Proposition} \label{prop:dimensions_NC_M}
    Let $\Mca$ be a finite unitary magma. For all $n \geq 2$,
    \begin{equation} \label{equ:dimensions_NC_M}
        \dim \NC\Mca(n) =
        \sum_{0 \leq k \leq n - 2}
            m^{n + k + 1} (m - 1)^{n - k - 2} \;
            \Nar(n, k),
    \end{equation}
    where $m := \# \Mca$.
\end{Proposition}
\begin{proof}
    As shown by Proposition~\ref{prop:map_NC_M_bubble_tree}, each noncrossing $\Mca$-clique
    $\Pfr$ of $\NC\Mca(n)$ can be encoded by a unique syntax tree $\BubbleTree(\Pfr)$ on
    $\Bubbles_\Mca$ satisfying some conditions. Moreover,
    Proposition~\ref{prop:generating_set_NC_M} shows that a noncrossing $\Mca$-clique can
    be expressed (not necessarily in a unique way) as partial compositions of several
    $\Mca$-triangles. By combining these two results, we obtain that a noncrossing
    $\Mca$-clique $\Pfr$ can be encoded by a syntax tree on $\Triangles_\Mca$ obtained from
    $\BubbleTree(\Pfr)$ by replacing each of its nodes $\Sfr$ of arity $\ell \geq 3$ by left
    comb binary syntax trees $\Sfr'$ on $\Triangles_\Mca$ satisfying
    \begin{equation} \label{equ:dimensions_NC_M_demo}
        \Sfr' :=
        \Corolla\left(\Qfr^1\right) \circ_1
        \Corolla\left(\Qfr^2\right) \circ_1 \dots \circ_1
        \Corolla\left(\Qfr^{\ell - 1}\right),
    \end{equation}
    where the $\Qfr^i$, $i \in [\ell - 1]$, are the unique $\Mca$-triangles such that, for
    every $i \in [2, \ell - 1]$, the base of $\Qfr^i$ is labeled by $\Unit_\Mca$, for every
    $j \in [\ell - 2]$, the first edge of $\Qfr^j$ is labeled by $\Unit_\Mca$, and
    $\Eval(\Sfr') = \Eval(\Sfr)$. Observe that, in~\eqref{equ:dimensions_NC_M_demo},
    brackets are not necessary since $\circ_1$ is associative. Therefore, $\Pfr$ can be
    encoded in a unique way as a binary syntax tree $\Tfr$ on $\Triangles_\Mca$ satisfying
    the following restrictions:
    \begin{enumerate}[fullwidth,label={(\it\roman*)}]
        \item the $\Mca$-triangles labeling the internal nodes of $\Tfr$ which are not the
        root have bases labeled by $\Unit_\Mca$;
        \item if $x$ and $y$ are two internal nodes of $\Tfr$ such that $y$ is the right
        child of $x$, the second edge of the bubble labeling $x$ is solid.
    \end{enumerate}

    To establish~\eqref{equ:dimensions_NC_M}, since the set of noncrossing $\Mca$-cliques
    forms the fundamental basis of $\NC\Mca$, we now have to count these binary trees.
    Consider a binary tree $\Tfr$ of arity $n \geq 2$ with exactly $k \in [0, n - 2]$
    internal nodes having an internal node as a left child. There are $m$ ways to label the
    base of the $\Mca$-triangle labeling the root of $\Tfr$, $m^k$ ways to label the first
    edges of the $\Mca$-triangles labeling the internal nodes of $\Tfr$ that have an
    internal node as left child, $m^n$ ways to label the first (resp.\ second) edges of the
    $\Mca$-triangles labeling the internal nodes of $\Tfr$ having a leaf as left (resp.\
    right) child, and, since there are exactly $n - k - 2$ internal nodes of $\Tfr$ having
    an internal node as a right child, there are $(m - 1)^{n - k - 2}$ ways to label the
    second edges of the $\Mca$-triangles labeling these internal nodes. Now, since $\Nar(n,
    k)$ counts the binary trees with $n$ leaves and exactly $k$ internal nodes having an
    internal node as a left child, and a binary tree with $n$ leaves can have at most $n -
    2$ internal nodes having an internal node as left child, \eqref{equ:dimensions_NC_M}
    follows.
\end{proof}

We can use Proposition~\ref{prop:dimensions_NC_M} to compute the first dimensions of
$\NC\Mca$. For instance, depending on $m := \# \Mca$, we have the following sequences of
dimensions:
\begin{subequations}
\begin{equation}
    1, 1, 1, 1, 1, 1, 1, 1,
    \qquad m = 1,
\end{equation}
\begin{equation}
    1, 8, 48, 352, 2880, 25216, 231168, 2190848,
    \qquad m = 2,
\end{equation}
\begin{equation}
    1, 27, 405, 7533, 156735, 349263, 81520425, 1967414265,
    \qquad m = 3,
\end{equation}
\begin{equation}
    1, 64, 1792, 62464, 2437120, 101859328, 4459528192, 201889939456.
    \qquad m = 4,
\end{equation}
\end{subequations}
The second one forms, except for the first terms, Sequence~\OEIS{A054726} of~\cite{Slo}. The
last two sequences are not listed in~\cite{Slo} at this time.

\subsection{Presentation and Koszulity}
The aim of this section is to establish a presentation by generators and relations of
$\NC\Mca$. For this, we will define an adequate rewrite rule on the set of syntax trees on
$\Triangles_\Mca$ and prove that it admits the required properties.

\subsubsection{Space of relations} \label{subsubsec:space_of_relations_NC_M}
Let $\Rel_{\NC\Mca}$ be the subspace of $\Free\left(\K \Angle{\Triangles_\Mca}\right)(3)$
generated by the elements
\begin{subequations}
\begin{equation} \label{equ:relation_1_NC_M}
    \Corolla\left(\Triangle{\Pfr_0}{\Pfr_1}{\Pfr_2}\right)
    \circ_1
    \Corolla\left(\Triangle{\Qfr_0}{\Qfr_1}{\Qfr_2}\right)
    -
    \Corolla\left(\Triangle{\Pfr_0}{\Rfr_1}{\Pfr_2}\right)
    \circ_1
    \Corolla\left(\Triangle{\Rfr_0}{\Qfr_1}{\Qfr_2}\right),
    \qquad
    \mbox{if } \Pfr_1 \Op \Qfr_0 = \Rfr_1 \Op \Rfr_0 \ne \Unit_\Mca,
\end{equation}
\begin{equation} \label{equ:relation_2_NC_M}
    \Corolla\left(\Triangle{\Pfr_0}{\Pfr_1}{\Pfr_2}\right)
    \circ_1
    \Corolla\left(\Triangle{\Qfr_0}{\Qfr_1}{\Qfr_2}\right)
    -
    \Corolla\left(\Triangle{\Pfr_0}{\Qfr_1}{\Rfr_2}\right)
    \circ_2
    \Corolla\left(\Triangle{\Rfr_0}{\Qfr_2}{\Pfr_2}\right),
    \qquad
    \mbox{if } \Pfr_1 \Op \Qfr_0 = \Rfr_2 \Op \Rfr_0 = \Unit_\Mca,
\end{equation}
\begin{equation} \label{equ:relation_3_NC_M}
    \Corolla\left(\Triangle{\Pfr_0}{\Pfr_1}{\Pfr_2}\right)
    \circ_2
    \Corolla\left(\Triangle{\Qfr_0}{\Qfr_1}{\Qfr_2}\right)
    -
    \Corolla\left(\Triangle{\Pfr_0}{\Pfr_1}{\Rfr_2}\right)
    \circ_2
    \Corolla\left(\Triangle{\Rfr_0}{\Qfr_1}{\Qfr_2}\right),
    \qquad
    \mbox{if } \Pfr_2 \Op \Qfr_0 = \Rfr_2 \Op \Rfr_0 \ne \Unit_\Mca,
\end{equation}
\end{subequations}
where $\Pfr$, $\Qfr$, and $\Rfr$ are $\Mca$-triangles.

\begin{Lemma} \label{lem:quadratic_relations_NC_M}
    Let $\Mca$ be a unitary magma, and $\Sfr$ and $\Tfr$ be two syntax trees of arity $3$ on
    $\Triangles_\Mca$. Then $\Sfr - \Tfr$ belongs to $\Rel_{\NC\Mca}$ if and only if
    $\Eval(\Sfr) = \Eval(\Tfr)$.
\end{Lemma}
\begin{proof}
    Assume first that $\Sfr - \Tfr$ belongs to $\Rel_{\NC\Mca}$. Then $\Sfr - \Tfr$ is a
    linear combination of elements of the form~\eqref{equ:relation_1_NC_M},
    \eqref{equ:relation_2_NC_M}, and~\eqref{equ:relation_3_NC_M}. Now, observe that, if
    $\Pfr$, $\Qfr$, and $\Rfr$ are three $\Mca$-triangles,
    \begin{enumerate}[fullwidth,label=(\alph*)]
        \item if $\delta := \Pfr_1 \Op \Qfr_0 = \Rfr_1 \Op \Rfr_0 \ne \Unit_\Mca$, we have
        \begin{equation}
            \Eval\left(
                \Corolla\left(\Triangle{\Pfr_0}{\Pfr_1}{\Pfr_2}\right)
                \circ_1
                \Corolla\left(\Triangle{\Qfr_0}{\Qfr_1}{\Qfr_2}\right)
            \right)
            =
            \SquareRight{\Qfr_1}{\Qfr_2}{\Pfr_2}{\Pfr_0}{\delta}
            =
            \Eval\left(
                \Corolla\left(\Triangle{\Pfr_0}{\Rfr_1}{\Pfr_2}\right)
                \circ_1
                \Corolla\left(\Triangle{\Rfr_0}{\Qfr_1}{\Qfr_2}\right)
            \right),
        \end{equation}
        \item if $\Pfr_1 \Op \Qfr_0 = \Rfr_2 \Op \Rfr_0 = \Unit_\Mca$, we have
        \begin{equation}
            \Eval\left(
                \Corolla\left(\Triangle{\Pfr_0}{\Pfr_1}{\Pfr_2}\right)
                \circ_1
                \Corolla\left(\Triangle{\Qfr_0}{\Qfr_1}{\Qfr_2}\right)
            \right)
            =
            \SquareN{\Qfr_1}{\Qfr_2}{\Pfr_2}{\Pfr_0}
            =
            \Eval\left(
                \Corolla\left(\Triangle{\Pfr_0}{\Qfr_1}{\Rfr_2}\right)
                \circ_2
                \Corolla\left(\Triangle{\Rfr_0}{\Qfr_2}{\Pfr_2}\right)
            \right),
        \end{equation}
        \item if $\delta := \Pfr_2 \Op \Qfr_0 = \Rfr_2 \Op \Rfr_0 \ne \Unit_\Mca$, we have
        \begin{equation}
            \Eval\left(
                \Corolla\left(\Triangle{\Pfr_0}{\Pfr_1}{\Pfr_2}\right)
                \circ_2
                \Corolla\left(\Triangle{\Qfr_0}{\Qfr_1}{\Qfr_2}\right)
            \right)
            =
            \SquareLeft{\Pfr_1}{\Qfr_1}{\Qfr_2}{\Pfr_0}{\delta}
            =
            \Eval\left(
                \Corolla\left(\Triangle{\Pfr_0}{\Pfr_1}{\Rfr_2}\right)
                \circ_2
                \Corolla\left(\Triangle{\Rfr_0}{\Qfr_1}{\Qfr_2}\right)
            \right).
        \end{equation}
    \end{enumerate}
    This shows that all evaluations in $\NC\Mca$ of~\eqref{equ:relation_1_NC_M},
    \eqref{equ:relation_2_NC_M}, and~\eqref{equ:relation_3_NC_M} are equal to zero.
    Therefore, $\Eval(\Sfr - \Tfr) = 0$, and hence we have $\Eval(\Sfr) - \Eval(\Tfr) = 0$
    and, as expected, $\Eval(\Sfr) = \Eval(\Tfr)$.

    We now assume that $\Eval(\Sfr) = \Eval(\Tfr)$ and let $\Rfr := \Eval(\Sfr)$. As
    $\Sfr$ is of arity $3$, $\Rfr$ also is of arity~$3$ and thus,
    \begin{equation}\label{equ:quadratic_relations_NC_M_demo_1}
        \Rfr \in \left\{
        \SquareRight{\Qfr_1}{\Qfr_2}{\Pfr_2}{\Pfr_0}{\delta},
        \SquareN{\Qfr_1}{\Qfr_2}{\Pfr_2}{\Pfr_0},
        \SquareLeft{\Pfr_1}{\Qfr_1}{\Qfr_2}{\Pfr_0}{\delta} :
        \Pfr, \Qfr \in \Triangles_\Mca,
        \delta \in \bar{\Mca}
        \right\}.
    \end{equation}
    Now, by definition of the partial composition of $\NC\Mca$, if $\Rfr$ has the form of
    the  first (resp.\ second, third) noncrossing $\Mca$-clique appearing
    in~\eqref{equ:quadratic_relations_NC_M_demo_1}, $\Sfr$ and $\Tfr$ are of the form of the
    first and second syntax trees of~\eqref{equ:relation_1_NC_M} (resp.\
    \eqref{equ:relation_2_NC_M}, \eqref{equ:relation_3_NC_M}). Hence, in all cases, $\Sfr -
    \Tfr$ is in~$\Rel_{\NC\Mca}$.
\end{proof}

\begin{Proposition} \label{prop:dimensions_relations_NC_M}
    Let $\Mca$ be a finite unitary magma. Then the dimension of the space $\Rel_{\NC\Mca}$
    satisfies
    \begin{equation}
        \dim \Rel_{\NC\Mca} = 2m^6 - 2m^5 + m^4,
    \end{equation}
    where $m := \# \Mca$.
\end{Proposition}
\begin{proof}
    For $x \in \Mca$, let $f(x)$ be the number of ordered pairs $(y, z) \in \Mca^2$ such
    that $x = y \Op z$. Since $\Mca$ is finite, $f : \Mca \to \N$ is a well-defined map.

    Let $\RelEq$ be the equivalence relation on the set of syntax trees on $\Triangles_\Mca$
    of arity $3$ satisfying $\Sfr \RelEq \Tfr$ if $\Sfr$ and $\Tfr$ are two such syntax
    trees satisfying $\Eval(\Sfr) = \Eval(\Tfr)$. Let also $C$ be the set of noncrossing
    $\Mca$-cliques of arity $3$. For $\Rfr \in C$, we denote the set of syntax trees $\Sfr$
    satisfying $\Eval(\Sfr) = \Rfr$ by $[\Rfr]_\RelEq$.
    Proposition~\ref{prop:generating_set_NC_M} says in particular that any $\Rfr \in C$ can
    be obtained by a partial composition of two $\Mca$-triangles, and hence all
    $[\Rfr]_\RelEq$ are nonempty sets and thus, are $\RelEq$-equivalence classes.

    Moreover, by Lemma~\ref{lem:quadratic_relations_NC_M}, for syntax trees $\Sfr$ and
    $\Tfr$, we have $\Sfr \RelEq \Tfr$ if and only if $\Sfr - \Tfr$ is in $\Rel_{\NC\Mca}$.
    For this reason, the dimension of $\Rel_{\NC\Mca}$ is linked with the cardinalities of
    all $\RelEq$-equivalence classes by
    \begin{equation} \label{equ:dimensions_relations_NC_M_demo_1}
        \dim \Rel_{\NC\Mca} = \sum_{\Rfr \in C} \Par{\# [\Rfr]_\RelEq - 1}.
    \end{equation}
    We now compute~\eqref{equ:dimensions_relations_NC_M_demo_1} by enumerating each
    $\RelEq$-equivalence class $[\Rfr]_\RelEq$.

    Observe that, since $\Rfr$ is of arity $3$, it can be of three different forms according
    to the presence of a solid diagonal.
    \begin{enumerate}[fullwidth,label=(\alph*)]
        \item If
        \begin{equation}
            \Rfr = \SquareRight{\Qfr_1}{\Qfr_2}{\Pfr_2}{\Pfr_0}{\delta}
        \end{equation}
        for some $\Pfr_0, \Pfr_2, \Qfr_1, \Qfr_2 \in \Mca$ and
        $\delta \in \bar{\Mca}$, to have $\Sfr \in [\Rfr]_\RelEq$, we
        necessarily have
        \begin{equation}
            \Sfr =
            \Corolla\left(\Triangle{\Pfr_0}{\Pfr_1}{\Pfr_2}\right)
            \circ_1
            \Corolla\left(\Triangle{\Qfr_0}{\Qfr_1}{\Qfr_2}\right)
        \end{equation}
        where $\Pfr_1, \Qfr_0 \in \Mca$ and $\Pfr_1 \Op \Qfr_0 = \delta$. Hence, $\#
        [\Rfr]_\RelEq = f(\delta)$.
        \item If
        \begin{equation}
            \Rfr = \SquareN{\Qfr_1}{\Qfr_2}{\Pfr_2}{\Pfr_0}
        \end{equation}
        for some $\Pfr_0, \Pfr_2, \Qfr_1, \Qfr_2 \in \Mca$, to have $\Sfr \in
            [\Rfr]_\RelEq$, we necessarily have
        \begin{equation}
            \Sfr \in \left\{
            \Corolla\left(\Triangle{\Pfr_0}{\Pfr_1}{\Pfr_2}\right)
            \circ_1
            \Corolla\left(\Triangle{\Qfr_0}{\Qfr_1}{\Qfr_2}\right),
            \Corolla\left(\Triangle{\Pfr_0}{\Qfr_1}{\Rfr_2}\right)
            \circ_2
            \Corolla\left(\Triangle{\Rfr_0}{\Qfr_2}{\Pfr_2}\right)
            \right\}
        \end{equation}
        where $\Pfr_1, \Qfr_0, \Rfr_0, \Rfr_2 \in \Mca$, $\Pfr_1 \Op \Qfr_0 = \Unit_\Mca$,
        and $\Rfr_2 \Op \Rfr_0 = \Unit_\Mca$. Hence, $\# [\Rfr]_\RelEq = 2f(\Unit_\Mca)$.
        \item Otherwise,
        \begin{equation}
            \Rfr = \SquareLeft{\Pfr_1}{\Qfr_1}{\Qfr_2}{\Pfr_0}{\delta}
        \end{equation}
        for some $\Pfr_0, \Pfr_1, \Qfr_1, \Qfr_2 \in \Mca$ and $\delta \in \bar{\Mca}$, and
            to have $\Sfr \in [\Rfr]_\RelEq$, we necessarily have
        \begin{equation}
            \Sfr =
            \Corolla\left(\Triangle{\Pfr_0}{\Pfr_1}{\Pfr_2}\right)
            \circ_2
            \Corolla\left(\Triangle{\Qfr_0}{\Qfr_1}{\Qfr_2}\right)
        \end{equation}
        where $\Pfr_2, \Qfr_0 \in \Mca$ and $\Pfr_2 \Op \Qfr_0 = \delta$. Hence, $\#
        [\Rfr]_{\RelEq} = f(\delta)$.
    \end{enumerate}
    Therefore, by using the fact that
    \begin{equation}
        \sum_{\delta \in \Mca} f(\delta) = m^2,
    \end{equation}
    from~\eqref{equ:dimensions_relations_NC_M_demo_1} we obtain
    \begin{equation}\begin{split}
        \dim \Rel_{\NC\Mca} & =
            \sum_{\substack{
                \Pfr_0, \Pfr_2, \Qfr_1, \Qfr_2 \in \Mca \\
                \delta \in \bar{\Mca}
            }}
            \Par{f(\delta) - 1}
        +
            \sum_{\Pfr_0, \Pfr_2, \Qfr_1, \Qfr_2 \in \Mca}
            \Par{2f(\Unit_\Mca) - 1}
        +
            \sum_{\substack{
                \Pfr_0, \Pfr_1, \Qfr_1, \Qfr_2 \in \Mca \\
                \delta \in \bar{\Mca}
            }}
            \Par{f(\delta) - 1}
        \\
        & = m^4
        \left(
            2
                \sum_{\delta \in \bar{\Mca}}
                \Par{f(\delta) - 1}
            + 2 f(\Unit_\Mca) - 1
        \right) \\
        & = 2m^6 - 2m^5 + m^4,
    \end{split}\end{equation}
    establishing the statement of the proposition.
\end{proof}

Observe that, by Proposition~\ref{prop:dimensions_relations_NC_M}, the dimension of
$\Rel_{\NC\Mca}$ only depends on the cardinality of $\Mca$ and not on its operation~$\Op$.

\subsubsection{Rewrite rule} \label{subsubsec:rewrite_rule_NC_M}
Let $\Rew$ be the rewrite rule on the set of syntax trees on $\Triangles_\Mca$ satisfying
\begin{subequations}
\begin{equation} \label{equ:rewrite_1_NC_M}
    \Corolla\left(\Triangle{\Pfr_0}{\Pfr_1}{\Pfr_2}\right)
    \circ_1
    \Corolla\left(\Triangle{\Qfr_0}{\Qfr_1}{\Qfr_2}\right)
    \Rew
    \Corolla\left(\Triangle{\Pfr_0}{\delta}{\Pfr_2}\right)
    \circ_1
    \Corolla\left(\TriangleEXX{\Unit_\Mca}{\Qfr_1}{\Qfr_2}\right),
    \qquad
    \mbox{if } \Qfr_0 \ne \Unit_\Mca,
    \mbox{ where } \delta := \Pfr_1 \Op \Qfr_0,
\end{equation}
\begin{equation} \label{equ:rewrite_2_NC_M}
    \Corolla\left(\Triangle{\Pfr_0}{\Pfr_1}{\Pfr_2}\right)
    \circ_1
    \Corolla\left(\Triangle{\Qfr_0}{\Qfr_1}{\Qfr_2}\right)
    \Rew
    \Corolla\left(\TriangleXXE{\Pfr_0}{\Qfr_1}{\Unit_\Mca}\right)
    \circ_2
    \Corolla\left(\TriangleEXX{\Unit_\Mca}{\Qfr_2}{\Pfr_2}\right),
    \qquad
    \mbox{if } \Pfr_1 \Op \Qfr_0 = \Unit_\Mca,
\end{equation}
\begin{equation} \label{equ:rewrite_3_NC_M}
    \Corolla\left(\Triangle{\Pfr_0}{\Pfr_1}{\Pfr_2}\right)
    \circ_2
    \Corolla\left(\Triangle{\Qfr_0}{\Qfr_1}{\Qfr_2}\right)
    \Rew
    \Corolla\left(\Triangle{\Pfr_0}{\Pfr_1}{\delta}\right)
    \circ_2
    \Corolla\left(\TriangleEXX{\Unit_\Mca}{\Qfr_1}{\Qfr_2}\right),
    \qquad
    \mbox{if } \Qfr_0 \ne \Unit_\Mca,
    \mbox{ where } \delta := \Pfr_2 \Op \Qfr_0,
\end{equation}
\end{subequations}
where $\Pfr$ and $\Qfr$ are $\Mca$-triangles.

\begin{Lemma} \label{lem:equivalence_relation_rewrite_rule_NC_M}
    Let $\Mca$ be a unitary magma. Then the vector space induced by the rewrite rule $\Rew$
    is $\Rel_{\NC\Mca}$.
\end{Lemma}
\begin{proof}
    Let $\Sfr$ and $\Tfr$ be two syntax trees on $\Triangles_\Mca$ such that $\Sfr \Rew
    \Tfr$. We have three cases to consider depending on the form of $\Sfr$ and~$\Tfr$.
    \begin{enumerate}[fullwidth,label=(\alph*)]
        \item if $\Sfr$ (resp.\ $\Tfr$) is of the form described by the left (resp.\ right)
        side of~\eqref{equ:rewrite_1_NC_M}, we have
        \begin{equation}
            \Eval(\Sfr)
            = \SquareRight{\Qfr_1}{\Qfr_2}{\Pfr_2}{\Pfr_0}{\delta}
            = \Eval(\Tfr),
        \end{equation}
        where $\delta := \Pfr_1 \Op \Qfr_0$.
        \item If $\Sfr$ (resp.\ $\Tfr$) is of the form described by the left (resp.\ right)
        side of~\eqref{equ:rewrite_2_NC_M}, we have
        \begin{equation}
            \Eval(\Sfr)
            = \SquareN{\Qfr_1}{\Qfr_2}{\Pfr_2}{\Pfr_0}
            = \Eval(\Tfr).
        \end{equation}
        \item Otherwise, $\Sfr$ (resp.\ $\Tfr$) is of the form described by the left (resp.\
        right) side of~\eqref{equ:rewrite_3_NC_M}. We have
        \begin{equation}
            \Eval(\Sfr)
            = \SquareLeft{\Pfr_1}{\Qfr_1}{\Qfr_2}{\Pfr_0}{\delta}
            = \Eval(\Tfr),
        \end{equation}
        where $\delta := \Pfr_2 \Op \Qfr_0$.
    \end{enumerate}
    Therefore, by Lemma~\ref{lem:quadratic_relations_NC_M} we have $\Sfr - \Tfr \in
    \Rel_{\NC\Mca}$ for each case. This leads to the fact that $\Sfr \RewTransSym \Tfr$
    implies $\Sfr - \Tfr \in \Rel_{\NC\Mca}$, and shows that the space induced by $\Rew$ is
    a subspace of~$\Rel_{\NC\Mca}$.

    We now assume that $\Sfr$ and $\Tfr$ are two syntax trees on $\Triangles_\Mca$ such
    that $\Sfr - \Tfr$ is a generator of $\Rel_{\NC\Mca}$ among~\eqref{equ:relation_1_NC_M},
    \eqref{equ:relation_2_NC_M}, and~\eqref{equ:relation_3_NC_M}.
    \begin{enumerate}[fullwidth,label=(\alph*)]
        \item If $\Sfr$ (resp.\ $\Tfr$) is of the form described by the left (resp.\ right)
        side of~\eqref{equ:relation_1_NC_M}, we have by~\eqref{equ:rewrite_1_NC_M},
        \begin{equation}
            \Sfr \RewTrans
            \Corolla\left(\Triangle{\Pfr_0}{\delta}{\Pfr_2}\right)
            \circ_1
            \Corolla\left(\TriangleEXX{\Unit_\Mca}{\Qfr_1}{\Qfr_2}\right)
            \quad \mbox{and} \quad
            \Tfr \RewTrans
            \Corolla\left(\Triangle{\Pfr_0}{\delta'}{\Pfr_2}\right)
            \circ_1
            \Corolla\left(\TriangleEXX{\Unit_\Mca}{\Qfr_1}{\Qfr_2}\right),
        \end{equation}
        where $\delta := \Pfr_1 \Op \Qfr_0$ and $\delta' := \Rfr_1 \Op \Rfr_0$. Since
        by~\eqref{equ:relation_1_NC_M}, we have $\delta = \delta'$, we obtain that $\Sfr
        \RewTransSym \Tfr$.
        \item If $\Sfr$ (resp.\ $\Tfr$) is of the form described by the left (resp.\ right)
        side of~\eqref{equ:relation_2_NC_M}, we have by~\eqref{equ:rewrite_2_NC_M} and
        by~\eqref{equ:rewrite_3_NC_M},
        \begin{equation}
            \Sfr \Rew
            \Corolla\left(\TriangleXXE{\Pfr_0}{\Qfr_1}{\Unit_\Mca}\right)
            \circ_2
            \Corolla\left(\TriangleEXX{\Unit_\Mca}{\Qfr_2}{\Pfr_2}\right)
            \quad \mbox{and} \quad
            \Tfr \RewTrans
            \Corolla\left(\TriangleXXE{\Pfr_0}{\Qfr_1}{\Unit_\Mca}\right)
            \circ_2
            \Corolla\left(\TriangleEXX{\Unit_\Mca}{\Qfr_2}{\Pfr_2}\right).
        \end{equation}
        We obtain that $\Sfr \RewTransSym \Tfr$.
        \item Otherwise, $\Sfr$ (resp.\ $\Tfr$) is of the form described by the left (resp.\
        right) side of~\eqref{equ:relation_3_NC_M}. We have by~\eqref{equ:rewrite_3_NC_M},
        \begin{equation}
            \Sfr \RewTrans
            \Corolla\left(\Triangle{\Pfr_0}{\Pfr_1}{\delta}\right)
            \circ_2
            \Corolla\left(
                \TriangleEXX{\Unit_\Mca}{\Qfr_1}{\Qfr_2}\right)
            \quad \mbox{and} \quad
            \Tfr \RewTrans
            \Corolla\left(\Triangle{\Pfr_0}{\Pfr_1}{\delta'}\right)
            \circ_2
            \Corolla\left(
                \TriangleEXX{\Unit_\Mca}{\Qfr_1}{\Qfr_2}\right),
        \end{equation}
        where $\delta := \Pfr_2 \Op \Qfr_0$ and $\delta' := \Rfr_2 \Op \Rfr_0$. Since
        by~\eqref{equ:relation_3_NC_M}, $\delta = \delta'$, we obtain that $\Sfr
        \RewTransSym \Tfr$.
    \end{enumerate}
    Hence, for each case, we have $\Sfr \RewTransSym \Tfr$. This shows that $\Rel_{\NC\Mca}$
    is a subspace of the space induced by $\Rew$. The statement of the lemma follows.
\end{proof}

\begin{Lemma} \label{lem:rewrite_rule_NC_M_terminating}
    For a unitary magma $\Mca$, the rewrite rule $\Rew$ is terminating.
\end{Lemma}
\begin{proof}
    Writing $T_n$ for the set of syntax trees on $\Triangles_\Mca$ of arity $n$, let
    $\phi : T_n \to \N^2$ be the map defined in the following way. For a syntax tree
    $\Tfr$ of $T_n$, $\phi(\Tfr) := (\alpha, \beta)$, where $\alpha$ is the sum of the number
    of internal nodes in the left subtree of $x$ taken over all internal nodes $x$ of
    $\Tfr$, and $\beta$ is the number of internal nodes of $\Tfr$ labeled by an
    $\Mca$-triangle whose base is not labeled by $\Unit_\Mca$. Let $\Sfr$ and $\Tfr$ be two
    syntax trees of $T_3$ such that $\Sfr \Rew \Tfr$. Due to the definition of $\Rew$, we
    have three configurations to explore. In what follows, $\eta : \Mca \to \N$ is the map
    satisfying
    $\eta(a) := 0$ if $a = \Unit_\Mca$ and $\eta(a) := 1$ otherwise.
    \begin{enumerate}[fullwidth,label=(\alph*)]
        \item If $\Sfr$ (resp.\ $\Tfr$) is of the form described by the left (resp.\ right)
        side of~\eqref{equ:rewrite_1_NC_M}, writing $\leq$ for the lexicographic order on
        $\N^2$, we have
        \begin{equation}\begin{split}
            \phi\left(
            \Corolla\left(\Triangle{\Pfr_0}{\Pfr_1}{\Pfr_2}\right)
            \circ_1
            \Corolla\left(\Triangle{\Qfr_0}{\Qfr_1}{\Qfr_2}\right)
            \right)
            & =
            \left(1, \eta(\Pfr_0) + 1\right) \\
            & >
            \left(1, \eta(\Pfr_0)\right)
            =
            \phi\left(
            \Corolla\left(\Triangle{\Pfr_0}{\delta}{\Pfr_2}\right)
            \circ_1
            \Corolla\left(\TriangleEXX{\Unit_\Mca}{\Qfr_1}{\Qfr_2}\right)
            \right),
        \end{split}\end{equation}
        where $\delta := \Pfr_1 \Op \Qfr_0$.
        \item If $\Sfr$ (resp.\ $\Tfr$) is of the form described by the left (resp.\ right)
        side of~\eqref{equ:rewrite_2_NC_M}, we have
        \begin{equation}\begin{split}
            \phi\left(
            \Corolla\left(\Triangle{\Pfr_0}{\Pfr_1}{\Pfr_2}\right)
            \circ_1
            \Corolla\left(\Triangle{\Qfr_0}{\Qfr_1}{\Qfr_2}\right)
            \right)
            & =
            \left(1, \eta(\Pfr_0) + \eta(\Qfr_0)\right) \\
            & >
            \left(0, \eta(\Pfr_0)\right)
            =
            \phi\left(
            \Corolla\left(\TriangleXXE{\Pfr_0}{\Qfr_1}{\Unit_\Mca}\right)
            \circ_2
            \Corolla\left(\TriangleEXX{\Unit_\Mca}{\Qfr_2}{\Pfr_2}\right)
            \right).
        \end{split}\end{equation}
        \item Otherwise, $\Sfr$ (resp.\ $\Tfr$) is of the form described by the left (resp.\
        right) side of~\eqref{equ:rewrite_3_NC_M}. We have
        \begin{equation}\begin{split}
            \phi\left(
            \Corolla\left(\Triangle{\Pfr_0}{\Pfr_1}{\Pfr_2}\right)
            \circ_2
            \Corolla\left(\Triangle{\Qfr_0}{\Qfr_1}{\Qfr_2}\right)
            \right)
            & =
            \left(0, \eta(\Pfr_0) + 1\right) \\
            & >
            \left(0, \eta(\Pfr_0)\right)
            =
            \phi\left(
            \Corolla\left(\Triangle{\Pfr_0}{\Pfr_1}{\delta}\right)
            \circ_2
            \Corolla\left(\TriangleEXX{\Unit_\Mca}{\Qfr_1}{\Qfr_2}\right)
            \right),
        \end{split}\end{equation}
        where $\delta := \Pfr_2 \Op \Qfr_0$.
    \end{enumerate}
    Therefore, for all syntax trees $\Sfr$ and $\Tfr$ such that $\Sfr \Rew \Tfr$,
    $\phi(\Sfr) > \phi(\Tfr)$. This implies that, for all syntax trees $\Sfr$ and $\Tfr$ such
    that $\Sfr \ne \Tfr$ and $\Sfr \RewTrans \Tfr$, we have $\phi(\Sfr) > \phi(\Tfr)$. Since
    $(0, 0)$ is the smallest element of $\N^2$ with respect to the lexicographic order
    $\leq$, the statement of the lemma follows.
\end{proof}

\begin{Lemma} \label{lem:rewrite_rule_NC_M_normal_forms}
    Let $\Mca$ be a unitary magma. The set of normal forms of the rewrite rule $\Rew$ is the
    set of syntax trees $\Tfr$ on $\Triangles_\Mca$ such that, for any internal nodes $x$
    and $y$ of $\Tfr$ where $y$ is a child of $x$,
    \begin{enumerate}[label={(\it\roman*)}]
        \item \label{item:rewrite_rule_NC_M_normal_forms_1} the base of the $\Mca$-triangle
        labeling $y$ is labeled by~$\Unit_\Mca$;
        \item \label{item:rewrite_rule_NC_M_normal_forms_2} if $y$ is a left child of $x$,
        the first edge of the $\Mca$-triangle labeling $x$ is not labeled by~$\Unit_\Mca$.
    \end{enumerate}
\end{Lemma}
\begin{proof}
    By Lemma~\ref{lem:rewrite_rule_NC_M_terminating}, $\Rew$ is terminating. Therefore,
    $\Rew$ admits normal forms, which are by definition the syntax trees on
    $\Triangles_\Mca$ that cannot be rewritten by~$\Rew$.

    Let $\Tfr$ be a normal form of $\Rew$. The fact that $\Tfr$
    satisfies~\ref{item:rewrite_rule_NC_M_normal_forms_1} is an immediate consequence of the
    fact that $\Tfr$ avoids the patterns appearing as left sides
    of~\eqref{equ:rewrite_1_NC_M} and~\eqref{equ:rewrite_3_NC_M}. Moreover, since $\Tfr$
    avoids the patterns appearing as left sides of~\eqref{equ:rewrite_2_NC_M}, one cannot
    have $\Pfr_1 \Op \Qfr_0 = \Unit_\Mca$, where $\Pfr$ (resp.\ $\Qfr$) is the label of $x$
    (resp.\ $y$). Since by~\ref{item:rewrite_rule_NC_M_normal_forms_1}, we have $\Qfr_0 =
    \Unit_\Mca$, we necessarily get $\Pfr_1 \ne \Unit_\Mca$. Hence, $\Tfr$
    satisfies~\ref{item:rewrite_rule_NC_M_normal_forms_2}.

    Conversely, if $\Tfr$ is a syntax tree on $\Triangles_\Mca$
    satisfying~\ref{item:rewrite_rule_NC_M_normal_forms_1}
    and~\ref{item:rewrite_rule_NC_M_normal_forms_2}, a direct inspection shows that one
    cannot rewrite $\Tfr$ by $\Rew$. Therefore, $\Tfr$ is a normal form of~$\Rew$.
\end{proof}

\begin{Lemma} \label{lem:generating_series_normal_forms_rewrite_rule_NC_M}
    Let $\Mca$ be a finite unitary magma. The generating series of the normal forms of the
    rewrite rule $\Rew$ is the Hilbert series $\Hilbert_{\NC\Mca(t)}$ of~$\NC\Mca$.
\end{Lemma}
\begin{proof}
    First, since by Lemma~\ref{lem:rewrite_rule_NC_M_terminating}, $\Rew$ is terminating,
    and since for $n \geq 1$, due to the finiteness of $\Mca$, there are finitely many
    syntax trees on $\Triangles_\Mca$ of arity $n$, the generating series $T(t)$ of the
    normal forms of $\Rew$ is well-defined.

    Let $S(t)$ be the generating series of the normal forms of $\Rew$ such that the bases of
    the $\Mca$-triangles labeling the roots are labeled by $\Unit_\Mca$. From the
    description of the normal forms of $\Rew$ provided by
    Lemma~\ref{lem:rewrite_rule_NC_M_normal_forms}, we obtain that $S(t)$ satisfies
    \begin{equation}
        S(t) = t + mtS(t) + (m - 1)m S(t)^2.
    \end{equation}
    Again by Lemma~\ref{lem:rewrite_rule_NC_M_normal_forms}, we have
    \begin{equation}
        T(t) = t + m(S(t) - t).
    \end{equation}
    A direct computation shows that $T(t)$ satisfies the algebraic equation
    \begin{equation}
        t + \left(m^3 - 2m^2 + 2m - 1\right)t^2
        + \left(2m^2t - 3mt + 2t - 1\right) T(t)
        + \left(m - 1\right) T(t)^2
        = 0.
    \end{equation}
    Hence, by Proposition~\ref{prop:Hilbert_series_NC_M}, we observe that $T(t) =
    \Hilbert_{\NC\Mca(t)}$.
\end{proof}

\begin{Lemma} \label{lem:rewrite_rule_NC_M_confluent}
    For a finite unitary magma $\Mca$, the rewrite rule $\Rew$ is confluent.
\end{Lemma}
\begin{proof}
    Arguing by contradiction, assume that $\Rew$ is not confluent. Since by
    Lemma~\ref{lem:rewrite_rule_NC_M_terminating}, $\Rew$ is terminating, there is an
    integer $n \geq 1$ and two normal forms $\Tfr$ and $\Tfr'$ of $\Rew$ of arity $n$ such
    that $\Tfr \ne \Tfr'$ and $\Tfr \CRewTransSym \Tfr'$. Now,
    Lemma~\ref{lem:quadratic_relations_NC_M} together with
    Lemma~\ref{lem:equivalence_relation_rewrite_rule_NC_M} implies that $\Eval(\Tfr) =
    \Eval(\Tfr')$. By Proposition~\ref{prop:generating_set_NC_M}, the map
    \begin{math}
        \Eval : \Free\left(\K \Angle{\Triangles_\Mca}\right) \to \NC\Mca
    \end{math}
    is surjective, leading to the fact that the number of normal forms of $\Rew$ of arity
    $n$ is greater than the number of noncrossing $\Mca$-cliques of arity $n$. However, by
    Lemma~\ref{lem:generating_series_normal_forms_rewrite_rule_NC_M}, there are as many
    normal forms of $\Rew$ of arity $n$ as noncrossing $\Mca$-cliques of arity $n$. This
    leads to a contradiction and proves the statement of the lemma.
\end{proof}

\subsubsection{Presentation and Koszulity}
The results of Sections~\ref{subsubsec:space_of_relations_NC_M}
and~\ref{subsubsec:rewrite_rule_NC_M} are finally used here to provide a presentation of
$\NC\Mca$ and prove that $\NC\Mca$ is a Koszul operad.

\begin{Theorem} \label{thm:presentation_NC_M}
    Let $\Mca$ be a finite unitary magma. Then $\NC\Mca$ admits the
    presentation~$\left(\Triangles_\Mca, \Rel_{\NC\Mca}\right)$.
\end{Theorem}
\begin{proof}
    First, since by Lemmas~\ref{lem:rewrite_rule_NC_M_terminating}
    and~\ref{lem:rewrite_rule_NC_M_confluent}, $\Rew$ is a convergent rewrite rule, and
    since by Lemma~\ref{lem:equivalence_relation_rewrite_rule_NC_M}, the space induced by
    $\Rew$ is $\Rel_{\NC\Mca}$, we can regard the underlying space of the quotient operad
    \begin{equation}
        \Oca :=
        \Free\left(\K \Angle{\Triangles_\Mca}\right) /_{\langle\Rel_{\NC\Mca}\rangle}
    \end{equation}
    as the linear span of all normal forms of $\Rew$. Moreover, as a consequence of
    Lemma~\ref{lem:quadratic_relations_NC_M}, the linear map
    \begin{math}
        \phi : \Oca \to \NC\Mca
    \end{math}
    defined for any normal form $\Tfr$ of $\Rew$ by $\phi(\Tfr) := \Eval(\Tfr)$ is an operad
    morphism. Now, by Proposition~\ref{prop:generating_set_NC_M}, $\phi$ is surjective.
    Moreover, by Lemma~\ref{lem:generating_series_normal_forms_rewrite_rule_NC_M}, the
    dimensions of the spaces $\Oca(n)$, $n \geq 1$ are the ones of $\NC\Mca(n)$. Hence,
    $\phi$ is an operad isomorphism and the statement of the theorem follows.
\end{proof}

We use Theorem~\ref{thm:presentation_NC_M} to express the presentations of the operads
$\NC\N_2$ and $\NC\Dbb_0$, where $\N_2$ is the cyclic monoid $\Z / 2 \Z$ and $\Dbb_0$ is the
multiplicative monoid on $\{0, 1\}$. The operad $\NC\N_2$ is generated by
\begin{equation}
    \Triangles_{\N_2} =
    \left\{
    \TriangleEEE{}{}{},
    \TriangleEXE{}{1}{},
    \TriangleEEX{}{}{1},
    \TriangleEXX{}{1}{1},
    \TriangleXEE{1}{}{},
    \TriangleXXE{1}{1}{},
    \TriangleXEX{1}{}{1},
    \Triangle{1}{1}{1}
    \right\},
\end{equation}
and these generators satisfies only the nontrivial relations
\begin{subequations}
\begin{equation}
    \TriangleXEX{a}{}{b_3} \circ_1 \Triangle{1}{b_1}{b_2}
    = \Triangle{a}{1}{b_3} \circ_1 \TriangleEXX{}{b_1}{b_2},
    \qquad a, b_1, b_2, b_3 \in \N_2,
\end{equation}
\begin{equation}
    \Triangle{a}{1}{b_3} \circ_1 \Triangle{1}{b_1}{b_2}
    = \TriangleXEX{a}{}{b_3} \circ_1 \TriangleEXX{}{b_1}{b_2}
    = \TriangleXXE{a}{b_1}{} \circ_2 \TriangleEXX{}{b_2}{b_3}
    = \Triangle{a}{b_1}{1} \circ_2 \Triangle{1}{b_2}{b_3},
    \qquad a, b_1, b_2, b_3 \in \N_2,
\end{equation}
\begin{equation}
    \TriangleXXE{a}{b_1}{} \circ_2 \Triangle{1}{b_2}{b_3}
    = \Triangle{a}{b_1}{1} \circ_2 \TriangleEXX{}{b_2}{b_3},
    \qquad a, b_1, b_2, b_3 \in \N_2.
\end{equation}
\end{subequations}
On the other hand, the operad $\NC\Dbb_0$ is generated by
\begin{equation}
    \Triangles_{\Dbb_0} =
    \left\{
    \TriangleEEE{}{}{},
    \TriangleEXE{}{0}{},
    \TriangleEEX{}{}{0},
    \TriangleEXX{}{0}{0},
    \TriangleXEE{0}{}{},
    \TriangleXXE{0}{0}{},
    \TriangleXEX{0}{}{0},
    \Triangle{0}{0}{0}
    \right\},
\end{equation}
and these generators satisfies only the nontrivial relations
\begin{subequations}
\begin{equation}
    \TriangleXEX{a}{}{b_3}
    \circ_1
    \Triangle{0}{b_1}{b_2}
    =
    \Triangle{a}{0}{b_3}
    \circ_1
    \Triangle{0}{b_1}{b_2}
    =
    \Triangle{a}{0}{b_3}
    \circ_1
    \TriangleEXX{}{b_1}{b_2},
    \qquad
    a, b_1, b_2, b_3 \in \Dbb_0,
\end{equation}
\begin{equation}
    \TriangleXEX{a}{}{b_3}
    \circ_1
    \TriangleEXX{}{b_1}{b_2}
    =
    \TriangleXXE{a}{b_1}{}
    \circ_2
    \TriangleEXX{}{b_2}{b_3},
    \qquad
    a, b_1, b_2, b_3 \in \Dbb_0,
\end{equation}
\begin{equation}
    \TriangleXXE{a}{b_1}{}
    \circ_2
    \Triangle{0}{b_2}{b_3}
    =
    \Triangle{a}{b_1}{0}
    \circ_2
    \Triangle{0}{b_2}{b_3}
    =
    \Triangle{a}{b_1}{0}
    \circ_2
    \TriangleEXX{}{b_2}{b_3},
    \qquad
    a, b_1, b_2, b_3 \in \Dbb_0.
\end{equation}
\end{subequations}

\begin{Theorem} \label{thm:Koszul_NC_M}
    For a finite unitary magma $\Mca$, $\NC\Mca$ is Koszul and the set of normal forms
    of $\Rew$ forms a Poincaré--Birkhoff--Witt basis of~$\NC\Mca$.
\end{Theorem}
\begin{proof}
    By Lemma~\ref{lem:equivalence_relation_rewrite_rule_NC_M} and
    Theorem~\ref{thm:presentation_NC_M}, the rewrite rule $\Rew$ is an orientation of the
    space of relations $\Rel_{\NC\Mca}$ of $\NC\Mca$. Moreover, by
    Lemmas~\ref{lem:rewrite_rule_NC_M_terminating}
    and~\ref{lem:rewrite_rule_NC_M_confluent}, this rewrite rule is convergent. Therefore,
    by Lemma~\ref{lem:koszulity_criterion_pbw}, $\NC\Mca$ is Koszul. Finally, the set of
    normal forms of $\Rew$ described by Lemma~\ref{lem:rewrite_rule_NC_M_normal_forms} is,
    by definition, a Poincaré--Birkhoff--Witt basis of~$\NC\Mca$.
\end{proof}

\subsection{Suboperads generated by bubbles}
In this section, we consider suboperads of $\NC\Mca$ generated by finite sets of
$\Mca$-bubbles. We assume here that $\Mca$ is endowed with an arbitrary total order so that
$\Mca = \{x_0, x_1, \dots\}$ with $x_0 = \Unit_\Mca$.

If $\Pfr$ is an $\Mca$-clique, the \Def{border} of $\Pfr$ is the word $\Border(\Pfr)$ of
length $n$ such that, for $i \in [n]$, we have $\Border(\Pfr)_i = \Pfr_i$.

\subsubsection{Treelike expressions on bubbles}
Let $B$ and $E$ be two subsets of $\Mca$. We write $\Bubbles_\Mca^{B, E}$ for the set of
$\Mca$-bubbles $\Pfr$ such that the bases of $\Pfr$ are labeled by $B$ and all edges of
$\Pfr$ are labeled by $E$. Moreover, we say that $\Mca$ is \Def{$(E, B)$-quasi-injective} if
for all $x, x' \in E$ and $y, y' \in B$, $x \Op y = x' \Op y' \ne \Unit_\Mca$ implies $x =
x'$ and $y = y'$.

\begin{Lemma} \label{lem:treelike_expression_suboperad_bubbles}
    Let $\Mca$ be a unitary magma, and $B$ and $E$ be two subsets of $\Mca$. If $\Mca$ is
    $(E, B)$-quasi-injective, then any $\Mca$-clique admits at most one treelike expression
    on $\Bubbles_\Mca^{B, E}$ of minimal degree.
\end{Lemma}
\begin{proof}
    Assume that $\Pfr$ is an $\Mca$-clique admitting a treelike expression on
    $\Bubbles_\Mca^{B, E}$. This implies that the base of $\Pfr$ is labeled by $B$, all
    solid diagonals of $\Pfr$ are labeled by $B \Op E$, and all edges of $\Pfr$ are labeled
    by $E$. By Proposition~\ref{prop:unique_decomposition_NC_M} and
    Lemma~\ref{lem:map_NC_M_bubble_tree_treelike_expression}, the tree $\Tfr :=
    \BubbleTree(\Pfr)$ is a treelike expression of $\Pfr$ on $\Bubbles_\Mca$ of minimal
    degree. Now, observe that $\Tfr$ is not necessarily a syntax tree on $\Bubbles_\Mca^{B,
    E}$ as required since some of its internal nodes may be labeled by bubbles that do not
    belong to $\Bubbles_\Mca^{B, E}$. Since $\Mca$ is $(E, B)$-quasi-injective, there is a
    unique way to relabel the internal nodes of $\Tfr$ by bubbles of $\Bubbles_\Mca^{B, E}$
    to obtain a syntax tree on $\Bubbles_\Mca^{B, E}$ such that $\Eval(\Tfr') =
    \Eval(\Tfr)$. By construction, $\Tfr'$ satisfies the properties of the statement of the
    lemma.
\end{proof}

\subsubsection{Dimensions} \label{subsubsec:dimensions_suboperads_triangles}
Let $G$ be a set of $\Mca$-bubbles and $\Xi := \left\{\xi_{x_0}, \xi_{x_1}, \dots\right\}$
be a set of noncommutative variables. Given $x_i \in \Mca$, let $\SeriesBubbles_{x_i}$ be
the series of $\N \langle\langle \Xi \rangle\rangle$ defined by
\begin{equation} \label{equ:series_bubbles}
    \SeriesBubbles_{x_i}\left(\xi_{x_0}, \xi_{x_1}, \dots\right) :=
    \sum_{\substack{
        \Pfr \in \Bubbles_\Mca^G \\
        \Pfr \ne \UnitClique}}
    \enspace \prod_{i \in [|\Pfr|]} \xi_{\Pfr_i},
\end{equation}
where $\Bubbles_\Mca^G$ is the set of $\Mca$-bubbles that can be obtained by partial
compositions of elements of $G$. Observe from~\eqref{equ:series_bubbles} that a
noncommutative monomial $u \in \Xi^{\geq 2}$ appears in $\SeriesBubbles_{x_i}$ with $1$ as
coefficient if and only if there is an $\Mca$-bubble with a base labeled by $x_i$ and with
$u$ as border in the suboperad of $\NC\Mca$ generated by $G$.

Moreover, for $x_i \in \Mca$, let the series $\SeriesElements_{x_i}$ of $\N \langle\langle t
\rangle\rangle$ defined by
\begin{equation}
    \SeriesElements_{x_i}(t) :=
    \SeriesBubbles_{x_i}\left(t + \bar \SeriesElements_{x_0}(t),
    t + \bar \SeriesElements_{x_1}(t), \dots\right),
\end{equation}
where, for $x_i \in \Mca$,
\begin{equation} \label{equ:bar_series_elements_based}
    \bar \SeriesElements_{x_i}(t) :=
    \sum_{\substack{
        x_j \in \Mca \\
        x_i \Op x_j \ne \Unit_\Mca
    }}
    \SeriesElements_{x_j}(t).
\end{equation}

\begin{Proposition} \label{prop:suboperads_NC_M_triangles_dimensions}
    Let $\Mca$ be a unitary magma and $G$ be a finite set of $\Mca$-bubbles such that,
    writing $B$ (resp.\ $E$) for the set of labels of the bases (resp.\ edges) of the
    elements of $G$, $\Mca$ is $(E, B)$-quasi-injective. Then, the Hilbert series
    $\Hilbert_{(\NC\Mca)^G}(t)$ of the suboperad of $\NC\Mca$ generated by $G$ satisfies
    \begin{equation} \label{equ:suboperads_NC_M_triangles_dimensions}
        \Hilbert_{(\NC\Mca)^G}(t) = t + \sum_{x_i \in \Mca} \SeriesElements_{x_i}(t).
    \end{equation}
\end{Proposition}
\begin{proof}
    By Lemma~\ref{lem:treelike_expression_suboperad_bubbles}, an $\Mca$-clique of
    $(\NC\Mca)^G$ admits exactly one treelike expression on $\Mca$-bubbles of $(\NC\Mca)^G$
    of minimal degree. For this reason, and as a consequence of the
    definition~\eqref{equ:bar_series_elements_based} of the series $\bar
    \SeriesElements_{x_i}(t)$, $x_i \in \Mca$, the series $\SeriesElements_{x_i}(t)$ is the
    generating series of all $\Mca$-cliques of $(\NC\Mca)^G$ different from $\UnitClique$
    and with a base labeled by $x_i \in \Mca$. Therefore, the
    expression~\eqref{equ:suboperads_NC_M_triangles_dimensions} for the Hilbert series of
    $(\NC\Mca)^G$ follows.
\end{proof}

As a side remark, Proposition~\ref{prop:suboperads_NC_M_triangles_dimensions} can be proved
by using the notion of bubble decompositions of operads developed in~\cite{CG14}. This
result provides a practical method to compute the dimensions of some suboperads
$(\NC\Mca)^G$ of $\NC\Mca$ by describing the series~\eqref{equ:series_bubbles} of the
bubbles of $\Bubbles_\Mca^G$. If $G$ satisfies the requirement of
Proposition~\ref{prop:suboperads_NC_M_triangles_dimensions}, this result also implies that
the Hilbert series of $(\NC\Mca)^G$ is algebraic.

\subsubsection[A cubic suboperad]{First example~: a cubic suboperad}
Consider the suboperad of $\NC\Ebb_2$ generated by
\begin{equation}
    G := \left\{
        \TriangleXEX{\Ett_1}{}{\Ett_1},
        \TriangleXEX{\Ett_2}{}{\Ett_2}
    \right\}.
\end{equation}
Computer experiments show that the generators of $(\NC\Ebb_2)^G$ do not satisfy any
nontrivial quadratic relation and that they satisfy the only four nontrivial cubic
relations
\begin{subequations}
\begin{equation} \label{equ:rel_suboperad_1_1}
    \TriangleXEX{\Ett_1}{}{\Ett_1}
    \circ_2
    \left(
    \TriangleXEX{\Ett_1}{}{\Ett_1}
    \circ_2
    \TriangleXEX{\Ett_1}{}{\Ett_1}
    \right)
    =
    \TriangleXEX{\Ett_1}{}{\Ett_1}
    \circ_2
    \left(
    \TriangleXEX{\Ett_2}{}{\Ett_2}
    \circ_2
    \TriangleXEX{\Ett_1}{}{\Ett_1}
    \right),
\end{equation}
\begin{equation} \label{equ:rel_suboperad_1_2}
    \TriangleXEX{\Ett_1}{}{\Ett_1}
    \circ_2
    \left(
    \TriangleXEX{\Ett_1}{}{\Ett_1}
    \circ_2
    \TriangleXEX{\Ett_2}{}{\Ett_2}
    \right)
    =
    \TriangleXEX{\Ett_1}{}{\Ett_1}
    \circ_2
    \left(
    \TriangleXEX{\Ett_2}{}{\Ett_2}
    \circ_2
    \TriangleXEX{\Ett_2}{}{\Ett_2}
    \right),
\end{equation}
\begin{equation} \label{equ:rel_suboperad_1_3}
    \TriangleXEX{\Ett_2}{}{\Ett_2}
    \circ_2
    \left(
    \TriangleXEX{\Ett_1}{}{\Ett_1}
    \circ_2
    \TriangleXEX{\Ett_1}{}{\Ett_1}
    \right)
    =
    \TriangleXEX{\Ett_2}{}{\Ett_2}
    \circ_2
    \left(
    \TriangleXEX{\Ett_2}{}{\Ett_2}
    \circ_2
    \TriangleXEX{\Ett_1}{}{\Ett_1}
    \right),
\end{equation}
\begin{equation} \label{equ:rel_suboperad_1_4}
    \TriangleXEX{\Ett_2}{}{\Ett_2}
    \circ_2
    \left(
    \TriangleXEX{\Ett_1}{}{\Ett_1}
    \circ_2
    \TriangleXEX{\Ett_2}{}{\Ett_2}
    \right)
    =
    \TriangleXEX{\Ett_2}{}{\Ett_2}
    \circ_2
    \left(
    \TriangleXEX{\Ett_2}{}{\Ett_2}
    \circ_2
    \TriangleXEX{\Ett_2}{}{\Ett_2}
    \right).
\end{equation}
\end{subequations}
Hence, $(\NC\Ebb_2)^G$ is not a quadratic operad. Moreover, it is possible to prove that
this operad does not admit any other nontrivial relations between its generators. This can
be performed by defining a rewrite rule on the syntax trees on $G$, consisting in rewriting
the left patterns of~\eqref{equ:rel_suboperad_1_1}, \eqref{equ:rel_suboperad_1_2},
\eqref{equ:rel_suboperad_1_3}, and~\eqref{equ:rel_suboperad_1_4} into their respective right
patterns, and by checking that this rewrite rule admits the required properties (like the
ones establishing the presentation of $\NC\Mca$ by Theorem~\ref{thm:presentation_NC_M}). The
existence of this nonquadratic operad shows that $\NC\Mca$ contains nonquadratic suboperads
even if it is quadratic.

One can prove by induction on the arity that the set of bubbles of $(\NC\Ebb_2)^G$ is the
set $B_1 \sqcup B_2$ where $B_1$ (resp.\ $B_2$) is the set of bubbles whose bases are
labeled by $\Ett_1$ (resp.\ $\Ett_2$) and the border is $\Unit \Ett_1$ (resp.\ $\Unit
\Ett_2$), or $\Unit \Unit \Unit^* \Ett_1$, or $\Unit \Unit \Unit^* \Ett_2$. Hence, we obtain
\begin{subequations}
\begin{equation}
    \SeriesBubbles_\Unit
    \left(\xi_\Unit, \xi_{\Ett_1}, \xi_{\Ett_2}\right)
    = 0,
\end{equation}
\begin{equation}
    \SeriesBubbles_{\Ett_1}
    \left(\xi_\Unit, \xi_{\Ett_1}, \xi_{\Ett_2}\right) =
    \frac{\xi_\Unit}{1 - \xi_\Unit}
    \left(\xi_{\Ett_1} + \xi_\Unit \xi_{\Ett_2}\right)
    =
    \SeriesBubbles_{\Ett_2}
    \left(\xi_\Unit, \xi_{\Ett_2}, \xi_{\Ett_1}\right),
\end{equation}
\end{subequations}
\noindent Moreover, one can check that $G$ satisfies the conditions required by
Proposition~\ref{prop:suboperads_NC_M_triangles_dimensions}. We hence have
\vspace{-2ex}
\begin{multicols}{2}
\begin{subequations}
\begin{equation}
    \bar\SeriesElements_\Unit(t)
    = \SeriesElements_{\Ett_1}(t) + \SeriesElements_{\Ett_2}(t),
\end{equation}

\begin{equation}
    \bar\SeriesElements_{\Ett_1}(t)
    = \SeriesElements_\Unit(t)
    = \bar\SeriesElements_{\Ett_2}(t),
\end{equation}
\end{subequations}
\end{multicols}
\noindent and
\begin{subequations}
\begin{equation}
    \SeriesElements_\Unit(t) = 0,
\end{equation}
\begin{equation}
    \SeriesElements_{\Ett_1}(t) =
    \SeriesBubbles_{\Ett_1}(
        t + \SeriesElements_{\Ett_1}(t) + \SeriesElements_{\Ett_2}(t),
        t, t)
    =
    \SeriesBubbles_{\Ett_2}(
        t + \SeriesElements_{\Ett_1}(t) + \SeriesElements_{\Ett_2}(t),
        t, t)
    =
    \SeriesElements_{\Ett_2}(t),
\end{equation}
\end{subequations}
By Proposition~\ref{prop:suboperads_NC_M_triangles_dimensions}, the Hilbert series of
$(\NC\Ebb_2)^G$ satisfies
\begin{equation}
    \Hilbert_{(\NC\Ebb_2)^G}(t) = t + \SeriesElements_\Unit(t) +
    \SeriesElements_{\Ett_1}(t) +
    \SeriesElements_{\Ett_2}(t) = t + 2 \SeriesElements_{\Ett_1}(t),
\end{equation}
and, by a straightforward computation, we obtain that this series satisfies the algebraic
equation
\begin{equation}
    t + (t - 1) \Hilbert_{(\NC\Ebb_2)^G}(t) +
    (2t + 1)\Hilbert_{(\NC\Ebb_2)^G}(t)^2 = 0.
\end{equation}
The first dimensions of $(\NC\Ebb_2)^G$ are
\begin{equation}
    1, 2, 8, 36, 180, 956, 5300, 30316,
\end{equation}
and form Sequence~\OEIS{A129148} of~\cite{Slo}.

\subsubsection[A suboperad of Motzkin paths]{Second example~: a suboperad of Motzkin paths}
Consider the suboperad of $\NC\Dbb_0$ generated by
\begin{equation}
    G := \left\{
        \TriangleEEE{}{}{},
        \SquareMotz
    \right\}.
\end{equation}
Computer experiments show that the generators of $(\NC\Dbb_0)^G$ satisfy the only four
nontrivial quadratic relations
\begin{subequations}
\begin{equation} \label{equ:rel_suboperad_2_1}
    \TriangleEEE{}{}{} \circ_1 \TriangleEEE{}{}{}
    =
    \TriangleEEE{}{}{} \circ_2 \TriangleEEE{}{}{}\,,
\end{equation}
\begin{equation} \label{equ:rel_suboperad_2_2}
    \SquareMotz \circ_1 \TriangleEEE{}{}{}
    =
    \TriangleEEE{}{}{} \circ_2 \SquareMotz\,,
\end{equation}
\begin{equation} \label{equ:rel_suboperad_2_3}
    \TriangleEEE{}{}{} \circ_1 \SquareMotz
    =
    \SquareMotz \circ_2 \TriangleEEE{}{}{}\,,
\end{equation}
\begin{equation} \label{equ:rel_suboperad_2_4}
    \SquareMotz \circ_1 \SquareMotz
    =
    \SquareMotz \circ_3 \SquareMotz\,.
\end{equation}
\end{subequations}
It is possible to prove that this operad does not admit any other nontrivial relations
between its generators. This can be performed by defining a rewrite rule on the syntax trees
on $G$, consisting in rewriting the left patterns of~\eqref{equ:rel_suboperad_2_1},
\eqref{equ:rel_suboperad_2_2}, \eqref{equ:rel_suboperad_2_3},
and~\eqref{equ:rel_suboperad_2_4} into their respective right patterns, and by checking that
this rewrite rule admits the required properties (like the ones establishing the
presentation of $\NC\Mca$ by Theorem~\ref{thm:presentation_NC_M}).

One can prove by induction on the arity that the set of bubbles of $(\NC\Dbb_0)^G$ is the
set of bubbles whose bases are labeled by $\Unit$ and borders are words of $\{\Unit,
0\}^{\geq 2}$ such that each occurrence of $0$ has a $\Unit$ immediately to its left and a
$\Unit$ immediately to its right. Hence, we obtain
\vspace{-2ex}
\begin{multicols}{2}
\begin{subequations}
\begin{equation}
    \SeriesBubbles_\Unit\left(\xi_\Unit, \xi_0\right)
    = \frac{1}{1 - \xi_\Unit - \xi_\Unit \xi_0} \xi_\Unit - \xi_\Unit,
\end{equation}

\begin{equation}
    \SeriesBubbles_0\left(\xi_\Unit, \xi_0\right) = 0.
\end{equation}
\end{subequations}
\end{multicols}
\noindent Moreover, one can check that $G$ satisfies the conditions required by
Proposition~\ref{prop:suboperads_NC_M_triangles_dimensions}. We hence have
\vspace{-2ex}
\begin{multicols}{2}
\begin{subequations}
\begin{equation}
    \bar\SeriesElements_\Unit(t) = \SeriesElements_0(t),
\end{equation}

\begin{equation}
    \bar\SeriesElements_0(t)
    = \SeriesElements_\Unit(t) + \SeriesElements_0(t),
\end{equation}
\end{subequations}
\end{multicols}
\noindent and
\vspace{-2ex}
\begin{multicols}{2}
\begin{subequations}
\begin{equation}
    \SeriesElements_\Unit(t)
    = \SeriesBubbles_\Unit\left(t, t + \SeriesElements_\Unit(t)\right),
\end{equation}

\begin{equation}
    \SeriesElements_0(t) = 0.
\end{equation}
\end{subequations}
\end{multicols}
\noindent By Proposition~\ref{prop:suboperads_NC_M_triangles_dimensions}, the Hilbert series
of $(\NC\Dbb_0)^G$ satisfies
\begin{equation}
    \Hilbert_{(\NC\Dbb_0)^G}(t) = t + \SeriesElements_\Unit(t),
\end{equation}
and, by a straightforward computation, we obtain that this series satisfies the algebraic
equation
\begin{equation}
    t + (t - 1) \Hilbert_{(\NC\Dbb_0)^G}(t) + t\Hilbert_{(\NC\Dbb_0)^G}(t)^2 = 0.
\end{equation}
The first dimensions of $(\NC\Dbb_0)^G$ are
\begin{equation}
    1, 1, 2, 4, 9, 21, 51, 127,
\end{equation}
and form Sequence~\OEIS{A001006} of~\cite{Slo}. The operad $(\NC\Dbb_0)^G$ has the same
presentation by generators and relations (and thus, the same Hilbert series) as the operad
$\Motz$ defined in~\cite{Gir15}, involving Motzkin paths. Hence, $(\NC\Dbb_0)^G$ and $\Motz$
are two isomorphic operads. Note in passing that these two operads are not isomorphic to the
operad $\Motzkin\Dbb_0$ constructed in Section~\ref{subsubsec:Motzkin_configurations}
of~\cite{Cliques1} and involving Motzkin configurations. Indeed, the sequence of the
dimensions of this last operad is a shifted version of the one of $(\NC\Dbb_0)^G$
and~$\Motz$.

\subsection{Algebras over the noncrossing clique operads}
We begin by briefly describing $\NC\Mca$-algebras in terms of relations between their
operations and the free $\NC\Mca$-algebras over one generator. We continue this section by
providing two ways to construct (not necessarily free) $\NC\Mca$-algebras. The first one
takes as input an associative algebra endowed with endofunctions satisfying some conditions,
and the second one takes as input a monoid.

\subsubsection{Relations}
From the presentation of $\NC\Mca$ established by Theorem~\ref{thm:presentation_NC_M}, any
$\NC\Mca$-algebra is a vector space $\Alg$ endowed with binary linear operations
\begin{equation}
    \TriangleOp{\Pfr_0}{\Pfr_1}{\Pfr_2} :
    \Alg \otimes \Alg \to \Alg,
    \qquad
    \Pfr \in \Triangles_\Mca,
\end{equation}
satisfying, for all $a_1, a_2, a_3 \in \Alg$, the relations
\begin{subequations}
\begin{equation} \label{equ:relation_NC_M_algebras_1}
    \left(a_1 \TriangleOp{\Qfr_0}{\Qfr_1}{\Qfr_2} a_2\right)
    \TriangleOp{\Pfr_0}{\Pfr_1}{\Pfr_2} a_3
    =
    \left(a_1 \TriangleOp{\Rfr_0}{\Qfr_1}{\Qfr_2} a_2\right)
    \TriangleOp{\Pfr_0}{\Rfr_1}{\Pfr_2} a_3,
    \qquad
    \mbox{if } \Pfr_1 \Op \Qfr_0 = \Rfr_1 \Op \Rfr_0 \ne \Unit_\Mca,
\end{equation}
\begin{equation} \label{equ:relation_NC_M_algebras_2}
    \left(a_1 \TriangleOp{\Qfr_0}{\Qfr_1}{\Qfr_2} a_2\right)
    \TriangleOp{\Pfr_0}{\Pfr_1}{\Pfr_2} a_3
    =
    a_1 \TriangleOp{\Pfr_0}{\Qfr_1}{\Rfr_2}
    \left(a_2 \TriangleOp{\Rfr_0}{\Qfr_2}{\Pfr_2} a_3 \right),
    \qquad
    \mbox{if } \Pfr_1 \Op \Qfr_0 = \Rfr_2 \Op \Rfr_0 = \Unit_\Mca,
\end{equation}
\begin{equation} \label{equ:relation_NC_M_algebras_3}
    a_1 \TriangleOp{\Pfr_0}{\Pfr_1}{\Pfr_2}
    \left( a_2 \TriangleOp{\Qfr_0}{\Qfr_1}{\Qfr_2} a_3 \right)
    =
    a_1 \TriangleOp{\Pfr_0}{\Pfr_1}{\Rfr_2}
    \left( a_2 \TriangleOp{\Rfr_0}{\Qfr_1}{\Qfr_2} a_3 \right),
    \qquad
    \mbox{if } \Pfr_2 \Op \Qfr_0 = \Rfr_2 \Op \Rfr_0 \ne \Unit_\Mca,
\end{equation}
\end{subequations}
where $\Pfr$, $\Qfr$, and $\Rfr$ are $\Mca$-triangles. Observe that $\Mca$ has to be finite
because Theorem~\ref{thm:presentation_NC_M} requires this property as premise.

\subsubsection{Free algebras over one generator}
From the realization of $\NC\Mca$ coming from its definition as a suboperad of $\Cli\Mca$,
the free $\NC\Mca$-algebra over one generator is the linear span $\NC\Mca$ of all
noncrossing $\Mca$-cliques endowed with the linear operations
\begin{equation}
    \TriangleOp{\Pfr_0}{\Pfr_1}{\Pfr_2} :
    \NC\Mca(n) \otimes \NC\Mca(m) \to \NC\Mca(n + m),
    \qquad
    \Pfr \in \Triangles_\Mca,
    n, m \geq 1,
\end{equation}
defined, for noncrossing $\Mca$-cliques $\Qfr$ and $\Rfr$, by
\begin{equation} \label{equ:free_algebra_MT_M_product}
    \Qfr \TriangleOp{\Pfr_0}{\Pfr_1}{\Pfr_2} \Rfr
    :=
    \left(\Triangle{\Pfr_0}{\Pfr_1}{\Pfr_2}
    \circ_2 \Rfr\right) \circ_1 \Qfr.
\end{equation}
In terms of $\Mca$-Schröder trees (see Section~\ref{subsubsec:M_Schroder_trees}),
\eqref{equ:free_algebra_MT_M_product} is the $\Mca$-Schröder tree obtained by grafting the
$\Mca$-Schröder trees of $\Qfr$ and $\Rfr$ respectively as left and right children of a
binary corolla having its edge adjacent to the root labeled by $\Pfr_0$, its first edge
labeled by $\Pfr_1 \Op \Qfr_0$, and second edge labeled by $\Pfr_2 \Op \Rfr_0$, and by
contracting each of these two edges when labeled by~$\Unit_\Mca$. For instance, in the free
$\NC\N_3$-algebra, we have
\begin{subequations}
\begin{equation}
    \begin{tikzpicture}[xscale=.34,yscale=.28,Centering]
        \node[Leaf](0)at(0.00,-2.00){};
        \node[Leaf](2)at(1.00,-4.00){};
        \node[Leaf](4)at(3.00,-4.00){};
        \node[Leaf](5)at(4.00,-2.00){};
        \node[Node](1)at(2.00,0.00){};
        \node[Node](3)at(2.00,-2.00){};
        \draw[Edge](0)edge[]node[EdgeLabel]{$1$}(1);
        \draw[Edge](2)--(3);
        \draw[Edge](3)edge[]node[EdgeLabel]{$2$}(1);
        \draw[Edge](4)edge[]node[EdgeLabel]{$1$}(3);
        \draw[Edge](5)edge[]node[EdgeLabel]{$1$}(1);
        \node(r)at(2.00,2){};
        \draw[Edge](r)edge[]node[EdgeLabel]{$2$}(1);
    \end{tikzpicture}
    \enspace \TriangleOp{1}{2}{1} \enspace
    \begin{tikzpicture}[xscale=.32,yscale=.25,Centering]
        \node[Leaf](0)at(0.00,-4.00){};
        \node[Leaf](2)at(1.00,-4.00){};
        \node[Leaf](3)at(2.00,-4.00){};
        \node[Leaf](5)at(4.00,-2.00){};
        \node[Node](1)at(1.00,-2.00){};
        \node[Node](4)at(3.00,0.00){};
        \draw[Edge](0)--(1);
        \draw[Edge](1)edge[]node[EdgeLabel]{$1$}(4);
        \draw[Edge](2)--(1);
        \draw[Edge](3)edge[]node[EdgeLabel]{$2$}(1);
        \draw[Edge](5)edge[]node[EdgeLabel]{$2$}(4);
        \node(r)at(3.00,1.25){};
        \draw[Edge](r)--(4);
    \end{tikzpicture}
    \enspace = \enspace
    \begin{tikzpicture}[xscale=.32,yscale=.16,Centering]
        \node[Leaf](0)at(0.00,-6.50){};
        \node[Leaf](10)at(8.00,-9.75){};
        \node[Leaf](12)at(10.00,-6.50){};
        \node[Leaf](2)at(1.00,-9.75){};
        \node[Leaf](4)at(3.00,-9.75){};
        \node[Leaf](5)at(4.00,-6.50){};
        \node[Leaf](7)at(6.00,-9.75){};
        \node[Leaf](9)at(7.00,-9.75){};
        \node[Node](1)at(2.00,-3.25){};
        \node[Node](11)at(9.00,-3.25){};
        \node[Node](3)at(2.00,-6.50){};
        \node[Node](6)at(5.00,0.00){};
        \node[Node](8)at(7.00,-6.50){};
        \draw[Edge](0)edge[]node[EdgeLabel]{$1$}(1);
        \draw[Edge](1)edge[]node[EdgeLabel]{$1$}(6);
        \draw[Edge](10)edge[]node[EdgeLabel]{$2$}(8);
        \draw[Edge](11)edge[]node[EdgeLabel]{$1$}(6);
        \draw[Edge](12)edge[]node[EdgeLabel]{$2$}(11);
        \draw[Edge](2)--(3);
        \draw[Edge](3)edge[]node[EdgeLabel]{$2$}(1);
        \draw[Edge](4)edge[]node[EdgeLabel]{$1$}(3);
        \draw[Edge](5)edge[]node[EdgeLabel]{$1$}(1);
        \draw[Edge](7)--(8);
        \draw[Edge](8)edge[]node[EdgeLabel]{$1$}(11);
        \draw[Edge](9)--(8);
        \node(r)at(5.00,3){};
        \draw[Edge](r)edge[]node[EdgeLabel]{$1$}(6);
    \end{tikzpicture}\,,
\end{equation}
\begin{equation}
    \begin{tikzpicture}[xscale=.34,yscale=.28,Centering]
        \node[Leaf](0)at(0.00,-2.00){};
        \node[Leaf](2)at(1.00,-4.00){};
        \node[Leaf](4)at(3.00,-4.00){};
        \node[Leaf](5)at(4.00,-2.00){};
        \node[Node](1)at(2.00,0.00){};
        \node[Node](3)at(2.00,-2.00){};
        \draw[Edge](0)edge[]node[EdgeLabel]{$1$}(1);
        \draw[Edge](2)--(3);
        \draw[Edge](3)edge[]node[EdgeLabel]{$2$}(1);
        \draw[Edge](4)edge[]node[EdgeLabel]{$1$}(3);
        \draw[Edge](5)edge[]node[EdgeLabel]{$1$}(1);
        \node(r)at(2.00,2){};
        \draw[Edge](r)edge[]node[EdgeLabel]{$2$}(1);
    \end{tikzpicture}
    \enspace \TriangleOp{1}{1}{1} \enspace
    \begin{tikzpicture}[xscale=.32,yscale=.25,Centering]
        \node[Leaf](0)at(0.00,-4.00){};
        \node[Leaf](2)at(1.00,-4.00){};
        \node[Leaf](3)at(2.00,-4.00){};
        \node[Leaf](5)at(4.00,-2.00){};
        \node[Node](1)at(1.00,-2.00){};
        \node[Node](4)at(3.00,0.00){};
        \draw[Edge](0)--(1);
        \draw[Edge](1)edge[]node[EdgeLabel]{$1$}(4);
        \draw[Edge](2)--(1);
        \draw[Edge](3)edge[]node[EdgeLabel]{$2$}(1);
        \draw[Edge](5)edge[]node[EdgeLabel]{$2$}(4);
        \node(r)at(3.00,1.25){};
        \draw[Edge](r)--(4);
    \end{tikzpicture}
    \enspace = \enspace
    \begin{tikzpicture}[xscale=.33,yscale=.19,Centering]
        \node[Leaf](0)at(0.00,-3.00){};
        \node[Leaf](1)at(1.00,-6.00){};
        \node[Leaf](11)at(10.00,-6.00){};
        \node[Leaf](3)at(3.00,-6.00){};
        \node[Leaf](5)at(5.00,-3.00){};
        \node[Leaf](6)at(6.00,-9.00){};
        \node[Leaf](8)at(7.00,-9.00){};
        \node[Leaf](9)at(8.00,-9.00){};
        \node[Node](10)at(9.00,-3.00){};
        \node[Node](2)at(2.00,-3.00){};
        \node[Node](4)at(4.00,0.00){};
        \node[Node](7)at(7.00,-6.00){};
        \draw[Edge](0)edge[]node[EdgeLabel]{$1$}(4);
        \draw[Edge](1)--(2);
        \draw[Edge](10)edge[]node[EdgeLabel]{$1$}(4);
        \draw[Edge](11)edge[]node[EdgeLabel]{$2$}(10);
        \draw[Edge](2)edge[]node[EdgeLabel]{$2$}(4);
        \draw[Edge](3)edge[]node[EdgeLabel]{$1$}(2);
        \draw[Edge](5)edge[]node[EdgeLabel]{$1$}(4);
        \draw[Edge](6)--(7);
        \draw[Edge](7)edge[]node[EdgeLabel]{$1$}(10);
        \draw[Edge](8)--(7);
        \draw[Edge](9)edge[]node[EdgeLabel]{$2$}(7);
        \node(r)at(4.00,3){};
        \draw[Edge](r)edge[]node[EdgeLabel]{$1$}(4);
    \end{tikzpicture}\,,
\end{equation}
\begin{equation}
    \begin{tikzpicture}[xscale=.34,yscale=.28,Centering]
        \node[Leaf](0)at(0.00,-2.00){};
        \node[Leaf](2)at(1.00,-4.00){};
        \node[Leaf](4)at(3.00,-4.00){};
        \node[Leaf](5)at(4.00,-2.00){};
        \node[Node](1)at(2.00,0.00){};
        \node[Node](3)at(2.00,-2.00){};
        \draw[Edge](0)edge[]node[EdgeLabel]{$1$}(1);
        \draw[Edge](2)--(3);
        \draw[Edge](3)edge[]node[EdgeLabel]{$2$}(1);
        \draw[Edge](4)edge[]node[EdgeLabel]{$1$}(3);
        \draw[Edge](5)edge[]node[EdgeLabel]{$1$}(1);
        \node(r)at(2.00,2){};
        \draw[Edge](r)edge[]node[EdgeLabel]{$2$}(1);
    \end{tikzpicture}
    \enspace \TriangleOp{1}{2}{0} \enspace
    \begin{tikzpicture}[xscale=.32,yscale=.25,Centering]
        \node[Leaf](0)at(0.00,-4.00){};
        \node[Leaf](2)at(1.00,-4.00){};
        \node[Leaf](3)at(2.00,-4.00){};
        \node[Leaf](5)at(4.00,-2.00){};
        \node[Node](1)at(1.00,-2.00){};
        \node[Node](4)at(3.00,0.00){};
        \draw[Edge](0)--(1);
        \draw[Edge](1)edge[]node[EdgeLabel]{$1$}(4);
        \draw[Edge](2)--(1);
        \draw[Edge](3)edge[]node[EdgeLabel]{$2$}(1);
        \draw[Edge](5)edge[]node[EdgeLabel]{$2$}(4);
        \node(r)at(3.00,1.25){};
        \draw[Edge](r)--(4);
    \end{tikzpicture}
    \enspace = \enspace
    \begin{tikzpicture}[xscale=.32,yscale=.16,Centering]
        \node[Leaf](0)at(0.00,-6.00){};
        \node[Leaf](10)at(7.00,-6.00){};
        \node[Leaf](11)at(8.00,-3.00){};
        \node[Leaf](2)at(1.00,-9.00){};
        \node[Leaf](4)at(3.00,-9.00){};
        \node[Leaf](5)at(4.00,-6.00){};
        \node[Leaf](7)at(5.00,-6.00){};
        \node[Leaf](9)at(6.00,-6.00){};
        \node[Node](1)at(2.00,-3.00){};
        \node[Node](3)at(2.00,-6.00){};
        \node[Node](6)at(6.00,0.00){};
        \node[Node](8)at(6.00,-3.00){};
        \draw[Edge](0)edge[]node[EdgeLabel]{$1$}(1);
        \draw[Edge](1)edge[]node[EdgeLabel]{$1$}(6);
        \draw[Edge](10)edge[]node[EdgeLabel]{$2$}(8);
        \draw[Edge](11)edge[]node[EdgeLabel]{$2$}(6);
        \draw[Edge](2)--(3);
        \draw[Edge](3)edge[]node[EdgeLabel]{$2$}(1);
        \draw[Edge](4)edge[]node[EdgeLabel]{$1$}(3);
        \draw[Edge](5)edge[]node[EdgeLabel]{$1$}(1);
        \draw[Edge](7)--(8);
        \draw[Edge](8)edge[]node[EdgeLabel]{$1$}(6);
        \draw[Edge](9)--(8);
        \node(r)at(6.00,3){};
        \draw[Edge](r)edge[]node[EdgeLabel]{$1$}(6);
    \end{tikzpicture}\,,
\end{equation}
\begin{equation}
    \begin{tikzpicture}[xscale=.34,yscale=.28,Centering]
        \node[Leaf](0)at(0.00,-2.00){};
        \node[Leaf](2)at(1.00,-4.00){};
        \node[Leaf](4)at(3.00,-4.00){};
        \node[Leaf](5)at(4.00,-2.00){};
        \node[Node](1)at(2.00,0.00){};
        \node[Node](3)at(2.00,-2.00){};
        \draw[Edge](0)edge[]node[EdgeLabel]{$1$}(1);
        \draw[Edge](2)--(3);
        \draw[Edge](3)edge[]node[EdgeLabel]{$2$}(1);
        \draw[Edge](4)edge[]node[EdgeLabel]{$1$}(3);
        \draw[Edge](5)edge[]node[EdgeLabel]{$1$}(1);
        \node(r)at(2.00,2){};
        \draw[Edge](r)edge[]node[EdgeLabel]{$2$}(1);
    \end{tikzpicture}
    \enspace \TriangleOp{1}{1}{0} \enspace
    \begin{tikzpicture}[xscale=.32,yscale=.25,Centering]
        \node[Leaf](0)at(0.00,-4.00){};
        \node[Leaf](2)at(1.00,-4.00){};
        \node[Leaf](3)at(2.00,-4.00){};
        \node[Leaf](5)at(4.00,-2.00){};
        \node[Node](1)at(1.00,-2.00){};
        \node[Node](4)at(3.00,0.00){};
        \draw[Edge](0)--(1);
        \draw[Edge](1)edge[]node[EdgeLabel]{$1$}(4);
        \draw[Edge](2)--(1);
        \draw[Edge](3)edge[]node[EdgeLabel]{$2$}(1);
        \draw[Edge](5)edge[]node[EdgeLabel]{$2$}(4);
        \node(r)at(3.00,1.25){};
        \draw[Edge](r)--(4);
    \end{tikzpicture}
    \enspace = \enspace
    \begin{tikzpicture}[xscale=.32,yscale=.18,Centering]
        \node[Leaf](0)at(0.00,-3.67){};
        \node[Leaf](1)at(1.00,-7.33){};
        \node[Leaf](10)at(8.00,-3.67){};
        \node[Leaf](3)at(3.00,-7.33){};
        \node[Leaf](5)at(4.00,-3.67){};
        \node[Leaf](6)at(5.00,-7.33){};
        \node[Leaf](8)at(6.00,-7.33){};
        \node[Leaf](9)at(7.00,-7.33){};
        \node[Node](2)at(2.00,-3.67){};
        \node[Node](4)at(4.00,0.00){};
        \node[Node](7)at(6.00,-3.67){};
        \draw[Edge](0)edge[]node[EdgeLabel]{$1$}(4);
        \draw[Edge](1)--(2);
        \draw[Edge](10)edge[]node[EdgeLabel]{$2$}(4);
        \draw[Edge](2)edge[]node[EdgeLabel]{$2$}(4);
        \draw[Edge](3)edge[]node[EdgeLabel]{$1$}(2);
        \draw[Edge](5)edge[]node[EdgeLabel]{$1$}(4);
        \draw[Edge](6)--(7);
        \draw[Edge](7)edge[]node[EdgeLabel]{$1$}(4);
        \draw[Edge](8)--(7);
        \draw[Edge](9)edge[]node[EdgeLabel]{$2$}(7);
        \node(r)at(4.00,3){};
        \draw[Edge](r)edge[]node[EdgeLabel]{$1$}(4);
    \end{tikzpicture}\,.
\end{equation}
\end{subequations}

\subsubsection[Associative algebras]{$\NC \Mca$-algebras from associative algebras}
Let $\Alg$ be an associative algebra with associative product denoted by $\OpAssoc$, and
\begin{equation} \label{equ:compatible_set}
    \omega_x : \Alg \to \Alg, \qquad x \in \Mca,
\end{equation}
be a family of linear maps, not necessarily associative algebra morphisms, indexed by the
elements of $\Mca$. We say that $\Alg$ together with this family~\eqref{equ:compatible_set}
of maps is
\Def{$\Mca$-compatible} if
\begin{equation} \label{equ:compatible_magma_on_algebra_1}
    \omega_{\Unit_\Mca} = \Id_\Alg
\end{equation}
where $\Id_\Alg$ is the identity map on $\Alg$, and
\begin{equation} \label{equ:compatible_magma_on_algebra_2}
    \omega_x \circ \omega_y = \omega_{x \Op y},
\end{equation}
for all $x, y \in \Mca$. We now use $\Mca$-compatible associative algebras to construct
$\NC\Mca$-algebras.

\begin{Theorem} \label{thm:NC_M_algebras}
    Let $\Mca$ be a finite unitary magma and $\Alg$ be an $\Mca$-compatible associative
    algebra. The vector space $\Alg$ endowed with the binary linear operations
    \begin{equation} \label{equ:NC_M_algebras}
        \TriangleOp{\Pfr_0}{\Pfr_1}{\Pfr_2} :
        \Alg \otimes \Alg \to \Alg,
        \qquad
        \Pfr \in \Triangles_\Mca,
    \end{equation}
    defined for each $\Mca$-triangle $\Pfr$ and $a_1, a_2 \in \Alg$ by
    \begin{equation} \label{equ:NC_M_algebras_def}
        a_1 \TriangleOp{\Pfr_0}{\Pfr_1}{\Pfr_2} a_2
        :=
        \omega_{\Pfr_0}\left(\omega_{\Pfr_1}\left(a_1\right)
        \OpAssoc \omega_{\Pfr_2}\left(a_2\right)\right),
    \end{equation}
    is an $\NC\Mca$-algebra.
\end{Theorem}
\begin{proof}
    We prove that the operations~\eqref{equ:NC_M_algebras} satisfy
    Relations~\eqref{equ:relation_NC_M_algebras_1}, \eqref{equ:relation_NC_M_algebras_2},
    and~\eqref{equ:relation_NC_M_algebras_3} of $\NC\Mca$-algebras. Since $\Mca$ is finite,
    this amounts to showing that these operations endow $\Alg$ with an $\NC\Mca$-algebra
    structure. For this, let $a_1$, $a_2$, and $a_3$ be three elements of $\Alg$, and
    $\Pfr$, $\Qfr$, and $\Rfr$ be three $\Mca$-triangles.
    \begin{enumerate}[fullwidth,label=(\alph*)]
        \item If $\Pfr_1 \Op \Qfr_0 = \Rfr_1 \Op \Rfr_0 \ne \Unit_\Mca$, then, since
        by~\eqref{equ:compatible_magma_on_algebra_2} we have $\omega_{\Pfr_1} \circ
        \omega_{\Qfr_0} = \omega_{\Rfr_1} \circ \omega_{\Rfr_0}$, we obtain
        \begin{equation}\begin{split}
            \left(a_1 \TriangleOp{\Qfr_0}{\Qfr_1}{\Qfr_2} a_2\right)
            \TriangleOp{\Pfr_0}{\Pfr_1}{\Pfr_2} a_3
            & =
            \omega_{\Qfr_0}(
            \omega_{\Qfr_1}(a_1) \OpAssoc
            \omega_{\Qfr_2}(a_2))
            \TriangleOp{\Pfr_0}{\Pfr_1}{\Pfr_2} a_3 \\
            & =
            \omega_{\Pfr_0}(
            \omega_{\Pfr_1}(\omega_{\Qfr_0}(
            \omega_{\Qfr_1}(a_1) \OpAssoc
            \omega_{\Qfr_2}(a_2)))
            \OpAssoc
            \omega_{\Pfr_2}(a_3)) \\
            & =
            \omega_{\Pfr_0}(
            \omega_{\Rfr_1}(\omega_{\Rfr_0}(
            \omega_{\Qfr_1}(a_1) \OpAssoc
            \omega_{\Qfr_2}(a_2)))
            \OpAssoc
            \omega_{\Pfr_2}(a_3)) \\
            & =
            \left(a_1 \TriangleOp{\Rfr_0}{\Qfr_1}{\Qfr_2} a_2\right)
            \TriangleOp{\Pfr_0}{\Rfr_1}{\Pfr_2} a_3.
        \end{split}\end{equation}
        Hence~\eqref{equ:relation_NC_M_algebras_1} holds.
        \item If $\Pfr_1 \Op \Qfr_0 = \Rfr_2 \Op \Rfr_0 = \Unit_\Mca$, then, since
        by~\eqref{equ:compatible_magma_on_algebra_1} we have $\omega_{\Pfr_1} \circ
        \omega_{\Qfr_0} = \omega_{\Rfr_2} \circ \omega_{\Rfr_0} = \Id_\Alg$ and since
        $\OpAssoc$ is associative, we get
        \begin{equation}\begin{split}
            \left(a_1 \TriangleOp{\Qfr_0}{\Qfr_1}{\Qfr_2} a_2\right)
            \TriangleOp{\Pfr_0}{\Pfr_1}{\Pfr_2} a_3
            & =
            \omega_{\Qfr_0}(\omega_{\Qfr_1}(a_1)
            \OpAssoc \omega_{\Qfr_2}(a_2))
            \TriangleOp{\Pfr_0}{\Pfr_1}{\Pfr_2} a_3 \\
            & =
            \omega_{\Pfr_0}(
            (\omega_{\Pfr_1}(
            \omega_{\Qfr_0}(\omega_{\Qfr_1}(a_1)
            \OpAssoc \omega_{\Qfr_2}(a_2))))
            \OpAssoc
            \omega_{\Pfr_2}(a_3)) \\
            & =
            \omega_{\Pfr_0}(
            \omega_{\Qfr_1}(a_1) \OpAssoc \omega_{\Qfr_2}(a_2)
            \OpAssoc
            \omega_{\Pfr_2}(a_3)) \\
            & =
            \omega_{\Pfr_0}(\omega_{\Qfr_1}(a_1) \OpAssoc
            (\omega_{\Qfr_2}(a_2) \OpAssoc
            \omega_{\Pfr_2}(a_3))) \\
            & =
            \omega_{\Pfr_0}(\omega_{\Qfr_1}(a_1) \OpAssoc
            \omega_{\Rfr_2}(\omega_{\Rfr_0}(\omega_{\Qfr_2}(a_2) \OpAssoc
            \omega_{\Pfr_2}(a_3)))) \\
            & = a_1 \TriangleOp{\Pfr_0}{\Qfr_1}{\Rfr_2}
            \left(a_2 \TriangleOp{\Rfr_0}{\Qfr_2}{\Pfr_2} a_3 \right).
        \end{split}\end{equation}
        Hence~\eqref{equ:relation_NC_M_algebras_2} holds.
        \item If $\Pfr_2 \Op \Qfr_0 = \Rfr_2 \Op \Rfr_0 \ne \Unit_\Mca$, then, since
        by~\eqref{equ:compatible_magma_on_algebra_2} we have $\omega_{\Pfr_2} \circ
        \omega_{\Qfr_0} = \omega_{\Rfr_2} \circ \omega_{\Rfr_0}$, we obtain
        \begin{equation}\begin{split}
            a_1 \TriangleOp{\Pfr_0}{\Pfr_1}{\Pfr_2}
            \left( a_2 \TriangleOp{\Qfr_0}{\Qfr_1}{\Qfr_2} a_3 \right)
            & =
            a_1 \TriangleOp{\Pfr_0}{\Pfr_1}{\Pfr_2}
            \omega_{\Qfr_0}
            (\omega_{\Qfr_1}(a_2) \OpAssoc \omega_{\Qfr_2}(a_3)) \\
            & =
            \omega_{\Pfr_0}(
            \omega_{\Pfr_1}(a_1)
            \OpAssoc
            \omega_{\Pfr_2}(
            \omega_{\Qfr_0}
            (\omega_{\Qfr_1}(a_2) \OpAssoc \omega_{\Qfr_2}(a_3)))) \\
            & =
            \omega_{\Pfr_0}(
            \omega_{\Pfr_1}(a_1)
            \OpAssoc
            \omega_{\Rfr_2}(
            \omega_{\Rfr_0}
            (\omega_{\Qfr_1}(a_2) \OpAssoc \omega_{\Qfr_2}(a_3)))) \\
            & =
            a_1 \TriangleOp{\Pfr_0}{\Pfr_1}{\Rfr_2}
            \left( a_2 \TriangleOp{\Rfr_0}{\Qfr_1}{\Qfr_2} a_3 \right).
        \end{split}\end{equation}
        Hence~\eqref{equ:relation_NC_M_algebras_3} holds.
    \end{enumerate}
    Consequently, $\Alg$ is an $\NC\Mca$-algebra.
\end{proof}

By Theorem~\ref{thm:NC_M_algebras}, $\Alg$ has the structure of an $\NC\Mca$-algebra. Hence,
there is a left action~$\cdot$ of the operad $\NC\Mca$ on the tensor algebra of $\Alg$ of
the form
\begin{equation}
    \cdot : \NC\Mca(n) \otimes \Alg^{\otimes n} \to \Alg,
    \qquad n \geq 1,
\end{equation}
whose definition comes from the ones of the operations~\eqref{equ:NC_M_algebras} and
Relation~\eqref{equ:algebra_over_operad}. We describe here an algorithm to compute the
action of an element of $\NC\Mca$ of arity $n$ on tensors $a_1 \otimes \dots \otimes a_n$
of $\Alg^{\otimes n}$. First, if $\Bfr$ is an $\Mca$-bubble of arity $n$,
\begin{equation} \label{equ:NC_M_algebras_action_bubbles}
    \Bfr \cdot \left(a_1 \otimes \dots \otimes a_n\right)
    = \omega_{\Bfr_0}\left(\prod_{i \in [n]}
    \omega_{\Bfr_i}\left(a_i\right) \right),
\end{equation}
where the product of~\eqref{equ:NC_M_algebras_action_bubbles} denotes the iterated version
of the associative product $\OpAssoc$ of~$\Alg$. If $\Pfr$ is a noncrossing $\Mca$-clique
of arity $n$, then $\Pfr$ acts recursively on $a_1 \otimes \dots \otimes a_n$ as follows.
We have
\begin{equation}
    \Pfr \cdot a_1 = a_1
\end{equation}
if $\Pfr = \UnitClique$, and
\begin{equation} \label{equ:NC_M_algebras_action_cliques}
    \Pfr \cdot \left(a_1 \otimes \dots \otimes a_n\right) =
    \Bfr \cdot \left(
        \left(\Rfr_1 \cdot \left(a_1 \otimes \dots
                \otimes a_{|\Rfr_1|}\right)\right)
        \otimes \dots \otimes
        \left(\Rfr_k \cdot
            \left(a_{|\Rfr_1| + \dots + |\Rfr_{k - 1}| + 1}
            \otimes \dots  \otimes a_n\right)\right)
    \right),
\end{equation}
where, by setting $\Tfr$ as the bubble tree $\BubbleTree(\Pfr)$ of $\Pfr$ (see
Section~\ref{subsubsec:treelike_bubbles}), $\Bfr$ and $\Rfr_1$, \dots, $\Rfr_k$ are the
unique $\Mca$-bubble and noncrossing $\Mca$-cliques such that
\begin{math}
    \Tfr = \Corolla(\Bfr)
    \circ [\BubbleTree(\Rfr_1), \dots, \BubbleTree(\Rfr_k)].
\end{math}

Here are a few examples of the construction provided by Theorem~\ref{thm:NC_M_algebras}.
\begin{description}[fullwidth]
    \item[Noncommutative polynomials and selected concatenation]
    Let us consider the unitary magma $\Sbb_\ell$ of the subsets of $[\ell]$ with the union
    as product. Let $A := \{a_j : j \in [\ell]\}$ be an alphabet of noncommutative letters.
    On the associative algebra $\K \langle A \rangle$ of polynomials on $A$, we define the
    linear maps
    \begin{equation}
        \omega_S : \K \langle A \rangle \to \K \langle A \rangle,
        \qquad S \in \Sbb_\ell,
    \end{equation}
    as follows. For $u \in A^*$ and $S \in \Sbb_\ell$, we set
    \begin{equation}
        \omega_S(u) :=
        \begin{cases}
            u, & \mbox{if } |u|_{a_j} \geq 1 \mbox{ for all } j \in S, \\
            0, & \mbox{otherwise}.
        \end{cases}
    \end{equation}
    Since, for all $u \in A^*$, $\omega_{\emptyset}(u) = u$ and $(\omega_S \circ
    \omega_{S'})(u) = \omega_{S \cup S'}(u)$, and $\emptyset$ is the unit of $\Sbb_\ell$, we
    obtain from Theorem~\ref{thm:NC_M_algebras} that the
    operations~\eqref{equ:NC_M_algebras} endow $ \K \langle A \rangle$ with an
    $\NC\Sbb_\ell$-algebra structure. For instance, with $\ell := 3$ we have
    \begin{subequations}
    \begin{equation}
        \left(a_1 + a_1 a_3 + a_2 a_2\right)
        \;
        \begin{tikzpicture}[scale=.52,Centering]
            \node[shape=coordinate](1)at(0,1){};
            \node[shape=coordinate](2)at(0.87,-0.5){};
            \node[shape=coordinate](3)at(-0.87,-0.5){};
            \draw[draw=ColB!90](1)edge[]node[CliqueLabel,font=\tiny]{$\{2\}$}(2);
            \draw[draw=ColB!90](1)edge[]node[CliqueLabel,font=\tiny]{$\{1\}$}(3);
            \draw[draw=ColB!90](2)edge[]node[CliqueLabel,font=\tiny]{$\{2, 3\}$}(3);
        \end{tikzpicture}
        \;
        \left(1 + a_3 + a_2 a_1\right)
        =
        a_1 a_3 a_2 a_1,
    \end{equation}
    \begin{equation}
        \left(a_1 + a_1 a_3 + a_2 a_2\right)
        \;
        \begin{tikzpicture}[scale=.52,Centering]
            \node[shape=coordinate](1)at(0,1){};
            \node[shape=coordinate](2)at(0.87,-0.5){};
            \node[shape=coordinate](3)at(-0.87,-0.5){};
            \draw[draw=ColB!90](1)edge[]node[CliqueLabel,font=\tiny]{$\emptyset$}(2);
            \draw[draw=ColB!90](1)edge[]node[CliqueLabel,font=\tiny]{$\{1\}$}(3);
            \draw[draw=ColB!90](2)edge[]node[CliqueLabel,font=\tiny]{$\{1, 3\}$}(3);
        \end{tikzpicture}
        \;
        \left(1 + a_3 + a_2 a_1\right)
        =
        2 \; a_1 a_3 + a_1 a_3 a_3 +
        a_1 a_3 a_2 a_1.
    \end{equation}
    \end{subequations}
    Moreover, to compute the action
    \begin{equation} \label{equ:example_action_S_clique}
        \begin{tikzpicture}[scale=1.3,Centering]
            \node[CliquePoint](1)at(-0.34,-0.94){};
            \node[CliquePoint](2)at(-0.87,-0.50){};
            \node[CliquePoint](3)at(-0.98,0.17){};
            \node[CliquePoint](4)at(-0.64,0.77){};
            \node[CliquePoint](5)at(-0.00,1.00){};
            \node[CliquePoint](6)at(0.64,0.77){};
            \node[CliquePoint](7)at(0.98,0.17){};
            \node[CliquePoint](8)at(0.87,-0.50){};
            \node[CliquePoint](9)at(0.34,-0.94){};
            \draw[CliqueEmptyEdge](1)edge[]node[CliqueLabel]{}(2);
            \draw[CliqueEdge](1)edge[bend left=30]node[CliqueLabel,near start]{$\{1\}$}(7);
            \draw[CliqueEmptyEdge](1)edge[]node[CliqueLabel]{}(9);
            \draw[CliqueEdge](2)edge[]node[CliqueLabel]{$\{1\}$}(3);
            \draw[CliqueEdge](2)edge[bend right=30]node[CliqueLabel]{$\{1\}$}(4);
            \draw[CliqueEdge](3)edge[]node[CliqueLabel]{$\{2\}$}(4);
            \draw[CliqueEdge](4)edge[bend right=30]node[CliqueLabel]{$\{1, 2\}$}(6);
            \draw[CliqueEdge](4)edge[]node[CliqueLabel]{$\{2\}$}(5);
            \draw[CliqueEmptyEdge](5)edge[]node[CliqueLabel]{}(6);
            \draw[CliqueEdge](6)edge[]node[CliqueLabel]{$\{3\}$}(7);
            \draw[CliqueEmptyEdge](7)edge[]node[CliqueLabel]{}(8);
            \draw[CliqueEdge](7)edge[]node[CliqueLabel]{$\{1, 2\}$}(9);
            \draw[CliqueEmptyEdge](8)edge[]node[CliqueLabel]{}(9);
        \end{tikzpicture}
        \cdot
        (f \otimes f \otimes f \otimes f \otimes f \otimes f \otimes f
        \otimes f),
    \end{equation}
    where $f := a_1 + a_2 + a_3$, we use the above algorithm
    and~\eqref{equ:NC_M_algebras_action_bubbles}
    and~\eqref{equ:NC_M_algebras_action_cliques}. By presenting the computation on the
    bubble tree of the noncrossing $\Sbb_3$-clique of~\eqref{equ:example_action_S_clique},
    we obtain
    
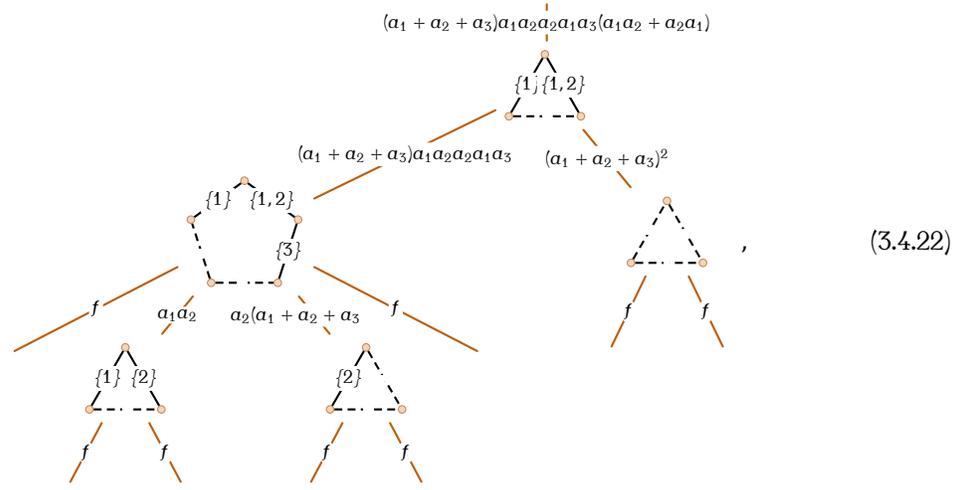
\begin{equation}
        \begin{tikzpicture}[xscale=.8,yscale=.6,Centering]
            \node(0)at(0.00,-6){};
            \node(1)at(1.00,-9){};
            \node(10)at(10.00,-6){};
            \node(12)at(12.00,-6){};
            \node(3)at(3.00,-9){};
            \node(5)at(5.00,-9){};
            \node(7)at(7.00,-9){};
            \node(8)at(8.00,-6){};
            \node(11)at(11.00,-3.25){
                \begin{tikzpicture}[scale=.55]
                    \node[CliquePoint](A1)at(-0.87,-0.50){};
                    \node[CliquePoint](A2)at(-0.00,1.00){};
                    \node[CliquePoint](A3)at(0.87,-0.50){};
                    \draw[CliqueEmptyEdge](A1)edge[]node[CliqueLabel]{}(A2);
                    \draw[CliqueEmptyEdge](A1)edge[]node[CliqueLabel]{}(A3);
                    \draw[CliqueEmptyEdge](A2)edge[]node[CliqueLabel]{}(A3);
                \end{tikzpicture}
            };
            \node(2)at(2.00,-6.50){
                \begin{tikzpicture}[scale=.55]
                    \node[CliquePoint](B1)at(-0.87,-0.50){};
                    \node[CliquePoint](B2)at(-0.00,1.00){};
                    \node[CliquePoint](B3)at(0.87,-0.50){};
                    \draw[CliqueEdge](B1)edge[]node[CliqueLabel]{$\{1\}$}(B2);
                    \draw[CliqueEmptyEdge](B1)edge[]node[CliqueLabel]{}(B3);
                    \draw[CliqueEdge](B2)edge[]node[CliqueLabel]{$\{2\}$}(B3);
                \end{tikzpicture}
            };
            \node(4)at(4.00,-3.25){
                \begin{tikzpicture}[scale=.75]
                \node[CliquePoint](C1)at(-0.59,-0.81){};
                \node[CliquePoint](C2)at(-0.95,0.31){};
                \node[CliquePoint](C3)at(-0.00,1.00){};
                \node[CliquePoint](C4)at(0.95,0.31){};
                \node[CliquePoint](C5)at(0.59,-0.81){};
                \draw[CliqueEmptyEdge](C1)edge[]node[CliqueLabel]{}(C2);
                \draw[CliqueEmptyEdge](C1)edge[]node[CliqueLabel]{}(C5);
                \draw[CliqueEdge](C2)edge[]node[CliqueLabel]{$\{1\}$}(C3);
                \draw[CliqueEdge](C3)edge[]node[CliqueLabel]{$\{1, 2\}$}(C4);
                \draw[CliqueEdge](C4)edge[]node[CliqueLabel]{$\{3\}$}(C5);
                \end{tikzpicture}
            };
            \node(6)at(6.00,-6.50){
                \begin{tikzpicture}[scale=.55]
                    \node[CliquePoint](D1)at(-0.87,-0.50){};
                    \node[CliquePoint](D2)at(-0.00,1.00){};
                    \node[CliquePoint](D3)at(0.87,-0.50){};
                    \draw[CliqueEdge](D1)edge[]node[CliqueLabel]{$\{2\}$}(D2);
                    \draw[CliqueEmptyEdge](D1)edge[]node[CliqueLabel]{}(D3);
                    \draw[CliqueEmptyEdge](D2)edge[]node[CliqueLabel]{}(D3);
                \end{tikzpicture}
            };
            \node(9)at(9.00,0.00){
                \begin{tikzpicture}[scale=.55]
                    \node[CliquePoint](E1)at(-0.87,-0.50){};
                    \node[CliquePoint](E2)at(-0.00,1.00){};
                    \node[CliquePoint](E3)at(0.87,-0.50){};
                    \draw[CliqueEdge](E1)edge[]node[CliqueLabel]{$\{1\}$}(E2);
                    \draw[CliqueEmptyEdge](E1)edge[]node[CliqueLabel]{}(E3);
                    \draw[CliqueEdge](E2)edge[]node[CliqueLabel]{$\{1, 2\}$}(E3);
                \end{tikzpicture}
            };
            \draw(0)edge[Edge]node[CliqueLabel]{$f$}(4);
            \draw(1)edge[Edge]node[CliqueLabel]{$f$}(2);
            \draw(10)edge[Edge]node[CliqueLabel]{$f$}(11);
            \draw(11)edge[Edge]node[CliqueLabel]{$(a_1 + a_2 + a_3)^2$}(9);
            \draw(12)edge[Edge]node[CliqueLabel]{$f$}(11);
            \draw(2)edge[Edge]node[CliqueLabel]{$a_1 a_2$}(4);
            \draw(3)edge[Edge]node[CliqueLabel]{$f$}(2);
            \draw(4)edge[Edge]node[CliqueLabel]{$(a_1 + a_2 + a_3) a_1 a_2 a_2 a_1 a_3$}(9);
            \draw(5)edge[Edge]node[CliqueLabel]{$f$}(6);
            \draw(6)edge[Edge]node[CliqueLabel]{$a_2 (a_1 + a_2 + a_3 \qquad$}(4);
            \draw(7)edge[Edge]node[CliqueLabel]{$f$}(6);
            \draw(8)edge[Edge]node[CliqueLabel]{$f$}(4);
            \node(r)at(9.00,2){};
            \draw(r)edge[Edge]node[CliqueLabel]
                {$(a_1 + a_2 + a_3) a_1 a_2 a_2 a_1 a_3 (a_1 a_2 + a_2 a_1)$}(9);
        \end{tikzpicture}\,,
    \end{equation}
    so that~\eqref{equ:example_action_S_clique} is equal to the polynomial
    \begin{math}
        (a_1 + a_2 + a_3) a_1 a_2 a_2 a_1 a_3 (a_1 a_2 + a_2 a_1).
    \end{math}

    \item[Noncommutative polynomials and constant term product]
    Let us now consider the unitary magma $\Dbb_0$. Let $A := \{a_1, a_2, \dots\}$ be an
    infinite alphabet of noncommutative letters. On the associative algebra $\K
    \langle A \rangle$ of polynomials on $A$, we define the linear maps
    \begin{equation}
        \omega_\Unit, \omega_0 :
        \K \langle A \rangle \to \K \langle A \rangle,
    \end{equation}
    as follows. For $u \in A^*$, we set $\omega_\Unit(u) := u$, and
    \begin{equation}
        \omega_0(u) :=
        \begin{cases}
            1, & \mbox{if } u = \epsilon, \\
            0, & \mbox{otherwise}.
        \end{cases}
    \end{equation}
    In other terms, $\omega_0(f)$ is the constant term, denoted by $f(0)$, of the polynomial
    $f \in \K \langle A \rangle$. Since $\omega_\Unit$ is the identity map on $\K \langle A
    \rangle$ and, for all $u \in A^*$,
    \begin{equation}
        (\omega_0 \circ \omega_0)(f)
        = (f(0))(0)
        = f(0)
        = \omega_0(f),
    \end{equation}
    we obtain from Theorem~\ref{thm:NC_M_algebras} that the
    operations~\eqref{equ:NC_M_algebras} endow $\K \langle A \rangle$ with an
    $\NC\Dbb_0$-algebra structure. For instance, for all polynomials $f_1$ and $f_2$ of $\K
    \langle A \rangle$, we have
    \vspace{-2ex}
    \begin{multicols}{2}
    \begin{subequations}
    \begin{equation}
        f_1 \TriangleOp{\Unit}{\Unit}{\Unit} f_2 = f_1 f_2,
    \end{equation}
    \begin{equation}
        f_1 \TriangleOp{0}{\Unit}{\Unit} f_2 = (f_1 f_2)(0) = f_1(0) \; f_2(0),
    \end{equation}

    \begin{equation} \label{equ:constant_term_product_example_1}
        f_1 \TriangleOp{\Unit}{0}{\Unit} f_2 = f_1(0) \; f_2,
    \end{equation}
    \begin{equation} \label{equ:constant_term_product_example_2}
        f_1 \TriangleOp{\Unit}{\Unit}{0} f_2 = f_1 \; (f_2(0)).
    \end{equation}
    \end{subequations}
    \end{multicols}
    \noindent If $f_1(0) = 1 = f_2(0)$, we obtain
    from~\eqref{equ:constant_term_product_example_1}
    and~\eqref{equ:constant_term_product_example_2} that
    \begin{equation}
        f_1
        \; \left(
        \TriangleOp{\Unit}{0}{\Unit} + \TriangleOp{\Unit}{\Unit}{0}
        \right) \;
        f_2
        = f_1(0) \; f_2 + f_1 \; (f_2(0))
        = f_1 + f_2.
    \end{equation}
\end{description}

\subsubsection[Monoids]{$\NC \Mca$-algebras from monoids}
If $\Mca$ is a monoid, with binary associative operation $\Op$ and unit $\Unit_\Mca$, we
write $\K \langle \Mca^* \rangle$ for the space of all noncommutative polynomials on $\Mca$,
viewed as an alphabet, with coefficients in $\K$. This space can be endowed with an
$\NC\Mca$-algebra structure as follows.

For $x \in \Mca$ and a word $w = w_1 \dots w_{|w|} \in \Mca^*$, let
\begin{equation}
    x * w := (x \Op w_1) \dots (x \Op w_{|w|}).
\end{equation}
This operation $*$ is linearly extended on the right to $\K \langle \Mca^* \rangle$.

\begin{Proposition} \label{prop:NC_M_algebras_monoid_polynomials}
    Let $\Mca$ be a finite monoid. The vector space $\K \langle \Mca^* \rangle$ endowed with
    the binary linear operations
    \begin{equation} \label{equ:NC_M_algebras_monoid_polynomials_op}
        \TriangleOp{\Pfr_0}{\Pfr_1}{\Pfr_2} :
        \K \langle \Mca^* \rangle \otimes \K \langle \Mca^* \rangle
        \to \K \langle \Mca^* \rangle,
        \qquad \Pfr \in \Triangles_\Mca,
    \end{equation}
    defined for each $\Mca$-triangle $\Pfr$ and $f_1, f_2 \in \K \langle \Mca^* \rangle$ by
    \begin{equation} \label{equ:NC_M_algebras_monoid_polynomials}
        f_1 \TriangleOp{\Pfr_0}{\Pfr_1}{\Pfr_2} f_2 :=
        \Pfr_0 * \left(\left(\Pfr_1 * f_1\right) \;
        \left(\Pfr_2 * f_2\right)\right),
    \end{equation}
    is an $\NC\Mca$-algebra.
\end{Proposition}
\begin{proof}
    This follows from Theorem~\ref{thm:NC_M_algebras} as a particular case of the general
    construction it provides. Indeed, $\K \langle \Mca^* \rangle$ is an associative algebra
    for the concatenation product of words. Moreover, by defining linear maps $\omega_x : \K
    \langle \Mca^* \rangle \to \K \langle \Mca^* \rangle$, $x \in \Mca$, by $\omega_x(u)
    := x * u$ for a word $u \in \Mca^*$, we obtain, since $\Mca$ is a monoid, that
    this family of maps satisfies~\eqref{equ:compatible_magma_on_algebra_1}
    and~\eqref{equ:compatible_magma_on_algebra_2}. Now, since the
    definition~\eqref{equ:NC_M_algebras_monoid_polynomials} is the specialization of the
    definition~\eqref{equ:NC_M_algebras_def} in this particular case, the statement of
    the proposition follows.
\end{proof}

Here are a few examples of the construction provided by
Proposition~\ref{prop:NC_M_algebras_monoid_polynomials}.
\begin{description}[fullwidth]
    \item[Words and double shifted concatenation]
    Consider the monoid $\N_\ell := \Z/_{\ell\Z}$ for an $\ell \geq 1$. By
    Proposition~\ref{prop:NC_M_algebras_monoid_polynomials}, the
    operations~\eqref{equ:NC_M_algebras_monoid_polynomials_op} endow $\K \left\langle
    \N_\ell^* \right\rangle$ with the structure of an $\NC\N_\ell$-algebra. For instance, in
    $\K \left\langle \N_4^* \right\rangle$ we have
    \begin{equation}
        0211 \; \TriangleOp{1}{2}{0} \; 312 = 3100023.
    \end{equation}

    \item[Words and erasing concatenation]
    Consider here the monoid $\Dbb_\ell$ for an $\ell \geq 0$ defined in~\cite{Cliques1}. By
    Proposition~\ref{prop:NC_M_algebras_monoid_polynomials}, the
    operations~\eqref{equ:NC_M_algebras_monoid_polynomials_op} endow $\K \langle \Dbb_\ell^*
    \rangle$ with the structure of an $\NC\Dbb_\ell$-algebra. For instance, for all words
    $u$ and $v$ of $\Dbb_\ell^*$, we have
    \vspace{-2ex}
    \begin{multicols}{2}
    \begin{subequations}
    \begin{equation}
        u \TriangleOp{\Unit}{\Unit}{\Unit} v = u v,
    \end{equation}
    \begin{equation}
        u \TriangleOp{\Dtt_i}{\Unit}{\Unit} v = (u v)_{\Dtt_i},
    \end{equation}

    \begin{equation}
        u \TriangleOp{0}{\Unit}{\Unit} v = 0^{|u| + |v|},
    \end{equation}
    \begin{equation}
        u \TriangleOp{\Unit}{\Dtt_i}{\Dtt_j} v = u_{\Dtt_i} \; v_{\Dtt_j},
    \end{equation}
    \end{subequations}
    \end{multicols}
    \noindent where, for a word $w$ of $\Dbb_\ell^*$ and an element $\Dtt_j$ of $\Dbb_\ell$,
    $j \in [\ell]$, $w_{\Dtt_j}$ is the word obtained from $w$ by replacing each occurrence
    of $\Unit$ by $\Dtt_j$ and each occurrence of $\Dtt_i$, $i \in [\ell]$, by $0$.
\end{description}

\section{Koszul dual} \label{sec:dual_NC_M}
Since, by Theorem~\ref{thm:presentation_NC_M}, the operad $\NC\Mca$ is binary and quadratic,
this operad admits a Koszul dual $\NC\Mca^!$. We end the study of $\NC\Mca$ by collecting
the main properties of~$\NC\Mca^!$.

\subsection{Presentation}
Let  $\Rel_{\NC\Mca}^!$ be the subspace of $\Free\left(\K \Angle{\Triangles_\Mca}\right)(3)$
generated by the elements
\begin{subequations}
\begin{equation} \label{equ:relation_1_NC_M_dual}
    \sum_{\substack{
        \Pfr_1, \Qfr_0 \in \Mca \\
        \Pfr_1 \Op \Qfr_0 = \delta
    }}
    \Corolla\left(\Triangle{\Pfr_0}{\Pfr_1}{\Pfr_2}\right)
    \circ_1
    \Corolla\left(\Triangle{\Qfr_0}{\Qfr_1}{\Qfr_2}\right),
    \qquad
    \Pfr_0, \Pfr_2, \Qfr_1, \Qfr_2 \in \Mca,
    \delta \in \bar{\Mca},
\end{equation}
\begin{equation} \label{equ:relation_2_NC_M_dual}
    \sum_{\substack{
        \Pfr_1, \Qfr_0 \in \Mca \\
        \Pfr_1 \Op \Qfr_0 = \Unit_\Mca
    }}
    \Par{
    \Corolla\left(\Triangle{\Pfr_0}{\Pfr_1}{\Pfr_2}\right)
    \circ_1
    \Corolla\left(\Triangle{\Qfr_0}{\Qfr_1}{\Qfr_2}\right)
    -
    \Corolla\left(\Triangle{\Pfr_0}{\Qfr_1}{\Pfr_1}\right)
    \circ_2
    \Corolla\left(\Triangle{\Qfr_0}{\Qfr_2}{\Pfr_2}\right)},
    \qquad
    \Pfr_0, \Pfr_2, \Qfr_1, \Qfr_2 \in \Mca,
\end{equation}
\begin{equation} \label{equ:relation_3_NC_M_dual}
    \sum_{\substack{
        \Pfr_2, \Qfr_0 \in \Mca \\
        \Pfr_2 \Op \Qfr_0 = \delta
    }}
    \Corolla\left(\Triangle{\Pfr_0}{\Pfr_1}{\Pfr_2}\right)
    \circ_2
    \Corolla\left(\Triangle{\Qfr_0}{\Qfr_1}{\Qfr_2}\right),
    \qquad
    \Pfr_0, \Pfr_1, \Qfr_1, \Qfr_2 \in \Mca,
    \delta \in \bar{\Mca},
\end{equation}
\end{subequations}
where $\Pfr$ and $\Qfr$ are $\Mca$-triangles.

\begin{Proposition} \label{prop:presentation_dual_NC_M}
    Let $\Mca$ be a finite unitary magma. Then the Koszul dual $\NC\Mca^!$ of $\NC\Mca$
    admits the presentation $\left(\Triangles_\Mca, \Rel_{\NC\Mca}^!\right)$.
\end{Proposition}
\begin{proof}
    Let
    \begin{equation}
        f := \sum_{\Tfr \in T_3} \lambda_\Tfr \Tfr
    \end{equation}
    be a generic element of $\Rel_{\NC\Mca}^!$, where $T_3$ is the set of all syntax trees
    on $\Triangles_\Mca$ of arity $3$ and the $\lambda_\Tfr$ are coefficients of $\K$. By
    definition of Koszul duality of operads, $\langle r, f \rangle = 0$ for all $r \in
    \Rel_{\NC\Mca}$, where $\langle -, - \rangle$ is the scalar product defined
    in~\eqref{equ:scalar_product_koszul}. Then, since $\Rel_{\NC\Mca}$ is the subspace of
    $\Free\left(\K \Angle{\Triangles_\Mca}\right)(3)$ generated
    by~\eqref{equ:relation_1_NC_M}, \eqref{equ:relation_2_NC_M},
    and~\eqref{equ:relation_3_NC_M}, we have
    \begin{subequations}
    \begin{equation} \label{prop:presentation_dual_NC_M_proof_1}
        \lambda_{\Corolla\left(\Triangle{\Pfr_0}{\Pfr_1}{\Pfr_2}\right)
            \circ_1
            \Corolla\left(\Triangle{\Qfr_0}{\Qfr_1}{\Qfr_2}\right)}
        -
        \lambda_{\Corolla\left(\Triangle{\Pfr_0}{\Rfr_1}{\Pfr_2}\right)
            \circ_1
            \Corolla\left(\Triangle{\Rfr_0}{\Qfr_1}{\Qfr_2}\right)}
        = 0,
        \qquad
        \Pfr_1 \Op \Qfr_0 = \Rfr_1 \Op \Rfr_0 \ne \Unit_\Mca,
    \end{equation}
    \begin{equation} \label{prop:presentation_dual_NC_M_proof_2}
        \lambda_{\Corolla\left(\Triangle{\Pfr_0}{\Pfr_1}{\Pfr_2}\right)
            \circ_1
            \Corolla\left(\Triangle{\Qfr_0}{\Qfr_1}{\Qfr_2}\right)}
        +
        \lambda_{\Corolla\left(\Triangle{\Pfr_0}{\Qfr_1}{\Rfr_2}\right)
            \circ_2
            \Corolla\left(\Triangle{\Rfr_0}{\Qfr_2}{\Pfr_2}\right)}
        = 0,
        \qquad
        \Pfr_1 \Op \Qfr_0 = \Rfr_2 \Op \Rfr_0 = \Unit_\Mca,
    \end{equation}
    \begin{equation} \label{prop:presentation_dual_NC_M_proof_3}
        \lambda_{\Corolla\left(\Triangle{\Pfr_0}{\Pfr_1}{\Pfr_2}\right)
            \circ_2
            \Corolla\left(\Triangle{\Qfr_0}{\Qfr_1}{\Qfr_2}\right)}
        -
        \lambda_{\Corolla\left(\Triangle{\Pfr_0}{\Pfr_1}{\Rfr_2}\right)
            \circ_2
            \Corolla\left(\Triangle{\Rfr_0}{\Qfr_1}{\Qfr_2}\right)}
        = 0,
        \qquad
        \Pfr_2 \Op \Qfr_0 = \Rfr_2 \Op \Rfr_0 \ne \Unit_\Mca,
    \end{equation}
    \end{subequations}
    where $\Pfr$, $\Qfr$, and $\Rfr$ are $\Mca$-triangles. This implies that $f$ is of the
    form
    \begin{equation}\begin{split}
        f & =
        \sum_{\substack{
            \Pfr_0, \Pfr_2, \Qfr_1, \Qfr_2 \in \Mca \\
            \delta \in \bar{\Mca}
        }}
        \lambda^{(1)}_{\Pfr_0, \Pfr_2, \Qfr_1, \Qfr_2, \delta}
        \sum_{\substack{
            \Pfr_1, \Qfr_0 \in \Mca \\
            \Pfr_1 \Op \Qfr_0 = \delta
        }}
        \Corolla\left(\Triangle{\Pfr_0}{\Pfr_1}{\Pfr_2}\right)
        \circ_1
        \Corolla\left(\Triangle{\Qfr_0}{\Qfr_1}{\Qfr_2}\right) \\
        & +
        \sum_{\Pfr_0, \Pfr_2, \Qfr_1, \Qfr_2 \in \Mca}
        \lambda^{(2)}_{\Pfr_0, \Pfr_2, \Qfr_1, \Qfr_2}
        \sum_{\substack{
            \Pfr_1, \Qfr_0 \in \Mca \\
            \Pfr_1 \Op \Qfr_0 = \Unit_\Mca
        }}
        \Par{
        \Corolla\left(\Triangle{\Pfr_0}{\Pfr_1}{\Pfr_2}\right)
        \circ_1
        \Corolla\left(\Triangle{\Qfr_0}{\Qfr_1}{\Qfr_2}\right)
        -
        \Corolla\left(\Triangle{\Pfr_0}{\Qfr_1}{\Pfr_1}\right)
        \circ_2
        \Corolla\left(\Triangle{\Qfr_0}{\Qfr_2}{\Pfr_2}\right)}
        \\
        & +
        \sum_{\substack{
            \Pfr_0, \Pfr_1, \Qfr_1, \Qfr_2 \in \Mca \\
            \delta \in \bar{\Mca}
        }}
        \lambda^{(3)}_{\Pfr_0, \Pfr_1, \Qfr_1, \Qfr_2, \delta}
        \sum_{\substack{
            \Pfr_2, \Qfr_0 \in \Mca \\
            \Pfr_2 \Op \Qfr_0 = \delta
        }}
        \Corolla\left(\Triangle{\Pfr_0}{\Pfr_1}{\Pfr_2}\right)
        \circ_2
        \Corolla\left(\Triangle{\Qfr_0}{\Qfr_1}{\Qfr_2}\right),
    \end{split}\end{equation}
    where, for $\Mca$-triangles $\Pfr$ and $\Qfr$ and $\delta \in \bar{\Mca}$, the
    $\lambda^{(1)}_{\Pfr_0, \Pfr_2, \Qfr_1, \Qfr_0, \delta}$, $\lambda^{(2)}_{\Pfr_0,
    \Pfr_2, \Qfr_1, \Qfr_2}$, and $\lambda^{(3)}_{\Pfr_1, \Pfr_0, \Qfr_1, \Qfr_2, \delta}$
    are coefficients of $\K$. Therefore, $f$ belongs to the space generated
    by~\eqref{equ:relation_1_NC_M_dual}, \eqref{equ:relation_2_NC_M_dual},
    and~\eqref{equ:relation_3_NC_M_dual}. Finally, since the coefficients of each of these
    relations satisfy~\eqref{prop:presentation_dual_NC_M_proof_1},
    \eqref{prop:presentation_dual_NC_M_proof_2},
    and~\eqref{prop:presentation_dual_NC_M_proof_3}, the statement of the proposition
    follows.
\end{proof}

We use Proposition~\ref{prop:presentation_dual_NC_M} to express the presentations of the
operads $\NC\N_2^!$ and $\NC\Dbb_0^!$. The operad $\NC\N_2^!$ is generated by
\begin{equation}
    \Triangles_{\N_2} =
    \left\{
    \TriangleEEE{}{}{},
    \TriangleEXE{}{1}{},
    \TriangleEEX{}{}{1},
    \TriangleEXX{}{1}{1},
    \TriangleXEE{1}{}{},
    \TriangleXXE{1}{1}{},
    \TriangleXEX{1}{}{1},
    \Triangle{1}{1}{1}
    \right\},
\end{equation}
and these generators satisfy only the nontrivial relations
\begin{subequations}
\begin{equation}
    \TriangleXEX{a}{}{b_3} \circ_1 \Triangle{1}{b_1}{b_2}
    +
    \Triangle{a}{1}{b_3} \circ_1 \TriangleEXX{}{b_1}{b_2}
    = 0,
    \qquad
    a, b_1, b_2, b_3 \in \N_2,
\end{equation}
\begin{equation}
    \Triangle{a}{1}{b_3} \circ_1 \Triangle{1}{b_1}{b_2}
    +
    \TriangleXEX{a}{}{b_3} \circ_1 \TriangleEXX{}{b_1}{b_2}
    =
    \TriangleXXE{a}{b_1}{} \circ_2 \TriangleEXX{}{b_2}{b_3}
    +
    \Triangle{a}{b_1}{1} \circ_2 \Triangle{1}{b_2}{b_3},
    \qquad
    a, b_1, b_2, b_3 \in \N_2,
\end{equation}
\begin{equation}
    \TriangleXXE{a}{b_1}{} \circ_2 \Triangle{1}{b_2}{b_3}
    +
    \Triangle{a}{b_1}{1} \circ_2 \TriangleEXX{}{b_2}{b_3}
    = 0,
    \qquad
    a, b_1, b_2, b_3 \in \N_2.
\end{equation}
\end{subequations}
On the other hand, the operad $\NC\Dbb_0^!$ is generated by
\begin{equation}
    \Triangles_{\Dbb_0} =
    \left\{
    \TriangleEEE{}{}{},
    \TriangleEXE{}{0}{},
    \TriangleEEX{}{}{0},
    \TriangleEXX{}{0}{0},
    \TriangleXEE{0}{}{},
    \TriangleXXE{0}{0}{},
    \TriangleXEX{0}{}{0},
    \Triangle{0}{0}{0}
    \right\},
\end{equation}
and these generators satisfies only the nontrivial relations
\begin{subequations}
\begin{equation}
    \TriangleXEX{a}{}{b_3} \circ_1 \Triangle{0}{b_1}{b_2}
    +
    \Triangle{a}{0}{b_3} \circ_1 \Triangle{0}{b_1}{b_2}
    +
    \Triangle{a}{0}{b_3} \circ_1 \TriangleEXX{}{b_1}{b_2}
    = 0,
    \qquad
    a, b_1, b_2, b_3 \in \Dbb_0,
\end{equation}
\begin{equation}
    \TriangleXEX{a}{}{b_3} \circ_1 \TriangleEXX{}{b_1}{b_2}
    =
    \TriangleXXE{a}{b_1}{} \circ_2 \TriangleEXX{}{b_2}{b_3},
    \qquad
    a, b_1, b_2, b_3 \in \Dbb_0,
\end{equation}
\begin{equation}
    \TriangleXXE{a}{b_1}{} \circ_2 \Triangle{0}{b_2}{b_3}
    +
    \Triangle{a}{b_1}{0} \circ_2 \Triangle{0}{b_2}{b_3}
    +
    \Triangle{a}{b_1}{0} \circ_2 \TriangleEXX{}{b_2}{b_3}
    = 0,
    \qquad
    a, b_1, b_2, b_3 \in \Dbb_0.
\end{equation}
\end{subequations}

\begin{Proposition} \label{prop:dimensions_relations_NC_M_dual}
    Let $\Mca$ be a finite unitary magma. Then the dimension of the space
    $\Rel_{\NC\Mca}^!$ is given by
    \begin{equation} \label{equ:dimensions_relations_NC_M_dual}
        \dim \Rel_{\NC\Mca}^! = 2m^5 - m^4,
    \end{equation}
    where $m := \# \Mca$.
\end{Proposition}
\begin{proof}
    To compute the dimension of the space of relations $\Rel_{\NC\Mca}^!$ of $\NC\Mca^!$, we
    consider the presentation of $\NC\Mca^!$ provided by
    Proposition~\ref{prop:presentation_dual_NC_M}. Consider the space $\Rel_1$ generated by
    the family consisting in the elements~\eqref{equ:relation_1_NC_M_dual}. Since this
    family is linearly independent and each of its element is totally specified by a tuple
    \begin{math}
        (\Pfr_0, \Pfr_2, \Qfr_1, \Qfr_2, \delta)
        \in \Mca^4 \times \bar{\Mca},
    \end{math}
    we obtain
    \begin{equation}
        \dim \Rel_1 = m^4 (m - 1).
    \end{equation}
    For the same reason, the dimension of the space $\Rel_3$ generated by the
    elements~\eqref{equ:relation_3_NC_M_dual} satisfies $\dim \Rel_3 = \dim \Rel_1$. Now,
    let $\Rel_2$ be the space generated by the elements~\eqref{equ:relation_2_NC_M_dual}.
    Since this family is linearly independent and each of its elements is totally specified
    by a tuple
    \begin{math}
        (\Pfr_0, \Pfr_2, \Qfr_1, \Qfr_2) \in \Mca^4,
    \end{math}
    we obtain
    \begin{equation}
        \dim \Rel_2 = m^4.
    \end{equation}
    Therefore, since
    \begin{equation}
        \Rel_{\NC\Mca}^! =
        \Rel_1 \oplus \Rel_2 \oplus \Rel_3,
    \end{equation}
    we obtain the stated formula~\eqref{equ:dimensions_relations_NC_M_dual} by summing the
    dimensions of $\Rel_1$, $\Rel_2$, and $\Rel_3$.
\end{proof}

Observe that, by Propositions~\ref{prop:dimensions_relations_NC_M}
and~\ref{prop:dimensions_relations_NC_M_dual}, we have
\begin{equation}\begin{split}
    \dim \Rel_{\NC\Mca} + \dim \Rel_{\NC\Mca}^!
        & = 2m^6 - 2m^5 + m^4
          + 2m^5 - m^4 \\
        & = 2m^6 \\
        & = \dim \Free\left(\K \Angle{\Triangles_\Mca}\right)(3),
\end{split}\end{equation}
as expected by Koszul duality, where $m := \# \Mca$.

\subsection{Hilbert series and dimensions}
An algebraic equation for the Hilbert series of $\NC\Mca^!$ is described and a formula
involving Narayana numbers to compute its coefficients is provided.

\begin{Proposition} \label{prop:Hilbert_series_NC_M_dual}
    Let $\Mca$ be a finite unitary magma. The Hilbert series $\Hilbert_{\NC\Mca^!}(t)$ of
    $\NC\Mca^!$ satisfies
    \begin{equation} \label{equ:Hilbert_series_NC_M_dual}
        t + (m - 1)t^2
        + \left(2m^2t - 3mt + 2t -1\right)\Hilbert_{\NC\Mca^!}(t)
        + \left(m^3 - 2m^2 + 2m - 1\right)\Hilbert_{\NC\Mca^!}(t)^2 = 0,
    \end{equation}
    where $m := \# \Mca$.
\end{Proposition}
\begin{proof}
    Let $G(t)$ be the generating series such that $G(-t)$
    satisfies~\eqref{equ:Hilbert_series_NC_M_dual}. Therefore, $G(t)$ satisfies
    \begin{equation} \label{equ:Hilbert_series_NC_M_dual_1}
        -t + (m - 1)t^2
        + \left(-2m^2t + 3mt - 2t -1\right)G(t)
        + \left(m^3 - 2m^2 + 2m - 1\right)G(t)^2 = 0,
    \end{equation}
    and, by solving~\eqref{equ:Hilbert_series_NC_M_dual_1} as a quadratic equation where $t$
    is the unknown, we obtain
    \begin{equation} \label{equ:Hilbert_series_NC_M_dual_2}
        t =
        \frac{1 + (2m^2 - 3m + 2)G(t)
        - \sqrt{1 + 2(2m^2 - m)G(t)
        + m^2G(t)^2}}{2(m - 1)}.
    \end{equation}
    Moreover, by Proposition~\ref{prop:Hilbert_series_NC_M}
    and~\eqref{equ:Hilbert_series_NC_M_function}, by setting $F(t) :=
    \Hilbert_{\NC\Mca}(-t)$, we have
    \begin{equation} \label{equ:Hilbert_series_NC_M_dual_3}
       F(G(t)) =
        \frac{1 + (2m^2 - 3m + 2)G(t)
        - \sqrt{1 + 2(2m^2 - m)G(t)
        + m^2G(t)^2}}{2(m - 1)}
        = t,
    \end{equation}
    showing that $F(t)$ and $G(t)$ are the inverses of each other for series composition.

    Now, since by Theorem~\ref{thm:Koszul_NC_M}, $\NC\Mca$ is a Koszul operad, the Hilbert
    series of $\NC\Mca$ and $\NC\Mca^!$ satisfy~\eqref{equ:Hilbert_series_Koszul_operads}.
    Therefore, \eqref{equ:Hilbert_series_NC_M_dual_3} implies that the Hilbert series of
    $\NC\Mca^!$ is the series $\Hilbert_{\NC\Mca^!}(t)$, satisfying the stated
    relation~\eqref{equ:Hilbert_series_NC_M_dual}.
\end{proof}

From Proposition~\ref{prop:Hilbert_series_NC_M_dual} we deduce that the Hilbert series of
$\NC\Mca^!$ satisfies
\begin{equation} \label{equ:Hilbert_function_NC_M_dual}
    \Hilbert_{\NC\Mca^!}(t) =
    \frac{1 - (2m^2 - 3m + 2)t -\sqrt{1 - 2(2m^3 - 2m^2 + m)t +m^2t^2}}
    {2(m^3 - 2m^2 + 2m - 1)},
\end{equation}
where $m := \# \Mca \ne 1$.

\begin{Proposition} \label{prop:dimensions_NC_M_dual}
    Let $\Mca$ be a finite unitary magma. For all $n \geq 2$,
    \begin{equation} \label{equ:dimensions_NC_M_dual}
        \dim \NC\Mca^!(n)
        = \sum_{0 \leq k \leq n - 2}
        m^{n + 1} (m(m - 1) + 1)^k (m (m - 1))^{n - k - 2} \; \Nar(n, k).
    \end{equation}
\end{Proposition}
\begin{proof}
    The proof consists in enumerating dual $\Mca$-cliques, introduced in the upcoming
    Section~\ref{subsec:basis_Cli_M_dual}. Indeed, by
    Proposition~\ref{prop:elements_NC_M_dual}, $\dim \NC\Mca^!(n)$ is equal to the number of
    dual $\Mca$-cliques of arity $n$. The expression for $\dim \NC\Mca^!(n)$ claimed
    by~\eqref{equ:dimensions_NC_M_dual} can be proved by using similar arguments as the ones
    intervening in the proof of Proposition~\ref{prop:dimensions_NC_M} for the
    expression~\eqref{equ:dimensions_NC_M} of $\dim \NC\Mca(n)$.
\end{proof}

We can use Proposition~\ref{prop:dimensions_NC_M_dual} to compute the first dimensions of
$\NC\Mca^!$. For instance, depending on $m := \# \Mca$, we have the following sequences of
dimensions:
\begin{subequations}
\begin{equation}
    1, 1, 1, 1, 1, 1, 1, 1,
    \qquad m = 1,
\end{equation}
\begin{equation}
    1, 8, 80, 992, 13760, 204416, 3180800, 51176960,
    \qquad m = 2,
\end{equation}
\begin{equation}
    1, 27, 1053, 51273, 2795715, 163318599, 9994719033, 632496651597,
    \qquad m = 3,
\end{equation}
\begin{equation}
    1, 64, 6400, 799744, 111923200, 16782082048, 2636161024000,
    428208345579520,
    \qquad m = 4.
\end{equation}
\end{subequations}
The second one is Sequence~\OEIS{A234596} of~\cite{Slo}. The last two sequences are not
listed in~\cite{Slo} at this time. It is worthwhile to observe that the dimensions of
$\NC\Mca^!$ for $\# \Mca = 2$ are the ones of the operad~$\BNC$ of bicolored noncrossing
configurations (see Section~\ref{subsec:operad_bnc}).

\subsection{Combinatorial basis} \label{subsec:basis_Cli_M_dual}
To describe a basis of $\NC\Mca^!$, we introduce the following sort of $\Mca$-decorated
cliques. A \Def{dual $\Mca$-clique} is an $\Mca^2$-clique such that its base and its edges
are labeled by pairs $(a, a) \in \Mca^2$, and all solid diagonals are labeled by pairs $(a,
b) \in \Mca^2$ with $a \ne b$. Observe that a non-solid diagonal of a dual $\Mca$-clique is
labeled by $(\Unit_\Mca, \Unit_\Mca)$. All definitions about $\Mca$-cliques given
in~\cite{Cliques1} remain valid for dual $\Mca$-cliques. For example,
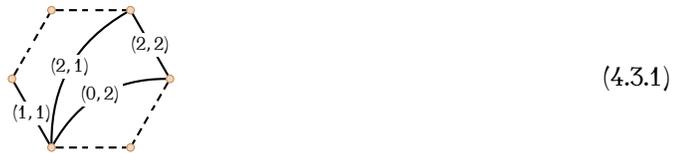
\begin{equation}
    \begin{tikzpicture}[scale=1.05,Centering]
        \node[CliquePoint](1)at(-0.50,-0.87){};
        \node[CliquePoint](2)at(-1.00,-0.00){};
        \node[CliquePoint](3)at(-0.50,0.87){};
        \node[CliquePoint](4)at(0.50,0.87){};
        \node[CliquePoint](5)at(1.00,0.00){};
        \node[CliquePoint](6)at(0.50,-0.87){};
        \draw[CliqueEdge](1)edge[]node[CliqueLabel]{$(1, 1)$}(2);
        \draw[CliqueEdge](1)edge[bend left=30]node[CliqueLabel]{$(0, 2)$}(5);
        \draw[CliqueEmptyEdge](1)edge[]node[CliqueLabel]{}(6);
        \draw[CliqueEdge](1)edge[bend left=30]node[CliqueLabel]{$(2, 1)$}(4);
        \draw[CliqueEmptyEdge](3)edge[]node[CliqueLabel]{}(4);
        \draw[CliqueEdge](4)edge[]node[CliqueLabel]{$(2, 2)$}(5);
        \draw[CliqueEmptyEdge](2)edge[]node[CliqueLabel]{}(3);
        \draw[CliqueEmptyEdge](5)edge[]node[CliqueLabel]{}(6);
    \end{tikzpicture}
\end{equation}
is a noncrossing dual $\N_3$-clique.

\begin{Proposition} \label{prop:elements_NC_M_dual}
    Let $\Mca$ be a finite unitary magma. The underlying graded vector space of $\NC\Mca^!$
    is the linear span of all noncrossing dual $\Mca$-cliques.
\end{Proposition}
\begin{proof}
    The statement of the proposition is equivalent to the fact that the generating series of
    noncrossing dual $\Mca$-cliques is the Hilbert series $\Hilbert_{\NC\Mca^!}(t)$ of
    $\NC\Mca^!$. From the definition of dual $\Mca$-cliques, we obtain that the set of dual
    $\Mca$-cliques of arity $n$, $n \geq 1$, is in bijection with the set of
    $\Mca^2$-Schröder trees of arity $n$ having the outgoing edges from the root and the
    edges connecting internal nodes with leaves labeled by pairs $(a, a) \in \Mca^2$, and
    the edges connecting two internal nodes labeled by pairs $(a, b) \in \Mca^2$ with $a \ne
    b$. The map $\BubbleTree$ defined in Section~\ref{subsubsec:treelike_bubbles} (see also
    Section~\ref{subsubsec:M_Schroder_trees}) realizes such a bijection. Let $T(t)$ be the
    generating series of these $\Mca^2$-Schröder trees, and let $S(t)$ be the generating
    series of the $\Mca^2$-Schröder trees of arities greater than $1$ and such that the
    outgoing edges from the roots and the edges connecting two internal nodes are labeled by
    pairs $(a, b) \in \Mca^2$ with $a \ne b$, and the edges connecting internal nodes with
    leaves are labeled by pairs $(a, a) \in \Mca^2$. From the description of these trees,
    we have
    \begin{equation}
        S(t) = m (m - 1) \frac{(mt + S(t))^2}{1 - mt - S},
    \end{equation}
    where $m := \# \Mca$. Moreover, for $m \ne 1$,  $T(t)$ satisfies
    \begin{equation}
        T(t) = t + \frac{S(t)}{m - 1},
    \end{equation}
    and we obtain that $T(t)$ admits~\eqref{equ:Hilbert_function_NC_M_dual} as solution.
    Then, by Proposition~\ref{prop:Hilbert_series_NC_M_dual}, for $m \ne 1$, this implies
    the statement of the proposition. If $m = 1$, it follows from
    Proposition~\ref{prop:presentation_dual_NC_M} that $\NC\Mca^!$ is isomorphic to the
    associative operad $\As$ (see for instance~\cite{LV12} for the definition of this
    operad). Hence, in this case, $\dim \NC\Mca^!(n) = 1$ for all $n \geq 1$. Since there is
    exactly one dual $\Mca$-clique of arity $n$ for $n \geq 1$, the statement of the
    proposition is satisfied.
\end{proof}

Proposition~\ref{prop:elements_NC_M_dual} gives a combinatorial description of the elements
of $\NC\Mca^!$. Nevertheless, for the time being we do not know a partial composition on the
linear span of these elements providing a realization of~$\NC\Mca^!$.

\section{Concrete constructions} \label{sec:concrete_constructions_nc}
The clique construction provides alternative definitions of known operads. We explore here
the cases of the operad $\NCP$ of based noncrossing trees, the operad $\FF_4$ of formal
fractions, and the operad $\BNC$ of bicolored noncrossing configurations.

\subsection{Rational functions and related operads}
We use here the noncrossing clique construction to interpret a few operads related to the
operad $\RatFct$ of rational functions of Loday~\cite{Lod10} (see also
Section~\ref{subsubsec:rational_functions} of~\cite{Cliques1}).

\subsubsection[Based noncrossing tree operad]{Dendriform and based noncrossing tree operads}
The \Def{operad of based noncrossing trees} $\NCP$ is an operad introduced in~\cite{Cha07}.
This operad is generated by two binary elements $\GDendr$ and $\DDendr$ satisfying one
nontrivial quadratic relation. The algebras over $\NCP$ are \Def{$\LOp$-algebras} and have
been studied in~\cite{Ler11}. We do not describe $\NCP$ in detail here because this is not
essential for the sequel. We just explain how to construct $\NCP$ through the clique
construction and interpret a known link between $\NCP$ and the dendriform operad through the
rational functions associated with $\Z$-cliques (see
Section~\ref{subsubsec:rational_functions} of~\cite{Cliques1}).

Let $\Oca_{\NCP}$ be the suboperad of $\Cli\Z$ generated by
\begin{equation}
    \left\{\TriangleEXE{}{-1}{},\TriangleEEX{}{}{-1}\right\}.
\end{equation}
By using Proposition~\ref{prop:suboperads_NC_M_triangles_dimensions}, we find that the
Hilbert series $\Hilbert_{\Oca_{\NCP}}(t)$ of
$\Oca_{\NCP}$ satisfies
\begin{equation} \label{equ:Hilbert_series_L_operad}
    t - \Hilbert_{\Oca_{\NCP}}(t) + 2 \Hilbert_{\Oca_{\NCP}}(t)^2
    - \Hilbert_{\Oca_{\NCP}}(t)^3 = 0.
\end{equation}
The first dimensions of $\Oca$ are
\begin{equation}
    1, 2, 7, 30, 143, 728, 3876, 21318,
\end{equation}
and form Sequence~\OEIS{A006013} of~\cite{Slo}. Moreover, one can see that
\begin{equation} \label{equ:relation_L_operad}
    \TriangleEEX{}{}{-1} \circ_1 \TriangleEXE{}{-1}{}
    = \TriangleEXE{}{-1}{} \circ_2 \TriangleEEX{}{}{-1},
\end{equation}
is the only nontrivial relation of degree $2$ between the generators of~$\Oca_{\NCP}$.

\begin{Proposition} \label{prop:construction_NCP}
    The operad $\Oca_{\NCP}$ is isomorphic to the operad $\NCP$.
\end{Proposition}
\begin{proof}
    Let $\phi_{\NCP} : \Oca_{\NCP}(2) \to \NCP(2)$ be the linear map satisfying
    \vspace{-2ex}
    \begin{multicols}{2}
    \begin{subequations}
    \begin{equation}
        \phi_{\NCP}\left(\TriangleEEX{}{}{-1}\right) = \; \GDendr,
    \end{equation}

    \begin{equation}
        \phi_{\NCP}\left(\TriangleEXE{}{-1}{}\right) = \; \DDendr,
    \end{equation}
    \end{subequations}
    \end{multicols}
    \noindent where $\GDendr$ and $\DDendr$ are the two binary generators of $\NCP$.
    In~\cite{Cha07}, a presentation of $\NCP$ is described wherein its generators satisfy
    one nontrivial relation of degree $2$. This relation can be obtained by replacing each
    $\Z$-clique appearing in~\eqref{equ:relation_L_operad} by its image by $\phi_{\NCP}$.
    For this reason, $\phi_{\NCP}$ uniquely extends to an operad morphism. Moreover, because
    the image of $\phi_{\NCP}$ contains all the generators of $\NCP$, this morphism is
    surjective. Finally, the Hilbert series of $\NCP$
    satisfies~\eqref{equ:Hilbert_series_L_operad}, so that $\Oca_{\NCP}$ and $\NCP$ have the
    same dimensions. Therefore, $\phi_{\NCP}$ is an operad isomorphism.
\end{proof}

Loday as shown in~\cite{Lod10} that the suboperad of $\RatFct$ generated by the rational
functions $f_1(u_1, u_2) := u_1^{-1}$ and $f_2(u_1, u_2) := u_2^{-1}$ is isomorphic to the
dendriform operad $\Dendr$~\cite{Lod01}. This operad is generated by two binary elements
$\GDendr$ and $\DDendr$ satisfying three nontrivial quadratic relations. An isomorphism
between $\Dendr$ and the suboperad of $\RatFct$ generated by $f_1$ and $f_2$ sends $\GDendr$
to $f_2$ and $\DDendr$ to $f_1$. The map $\Frac_\Id$ introduced in~\cite{Cliques1} is an
operad morphism from $\Cli\Z$ to $\RatFct$. Hence, the restriction of $\Frac_\Id$ on
$\Oca_{\NCP}$ is also an operad morphism from $\Oca_{\NCP}$ to $\RatFct$. Moreover, since
\vspace{-2ex}
\begin{multicols}{2}
\begin{subequations}
\begin{equation}
    \Frac_\Id\left(\TriangleEXE{}{-1}{}\right) = \frac{1}{u_1} = f_1,
\end{equation}

\begin{equation}
    \Frac_\Id\left(\TriangleEEX{}{}{-1}\right) = \frac{1}{u_2} = f_2,
\end{equation}
\end{subequations}
\end{multicols}
\noindent the map $\Frac_\Id$ is a surjective operad morphism from $\Oca_{\NCP}$
to~$\Dendr$.

\subsubsection{Operad of formal fractions} \label{subsubsec:operad_ff}
The \Def{operad of formal fractions} $\FF$ is an operad introduced in~\cite{CHN16}. Its
elements of arity $n \geq 1$ are fractions whose numerators and denominators are formal
products of subsets of $[n]$. For instance,
\begin{equation}
    \frac{\{1, 3, 4\} \{2\} \{4, 6\}}{\{2, 3, 5\} \{4\}}
\end{equation}
is an element of arity $6$ of $\FF$. We do not describe the partial composition of this
operad since its knowledge is not essential for the sequel. The operad $\FF$ admits a
suboperad $\FF_4$, defined as the binary suboperad of $\FF$ generated by
\begin{equation} \label{equ:generators_FF4}
    \left\{
        \frac{1}{\{1\} \{1, 2\}},
        \frac{1}{\{2\} \{1, 2\}},
        \frac{1}{\{1, 2\}},
        \frac{1}{\{1\} \{2\}}
    \right\}.
\end{equation}
We explain here how to construct $\FF_4$ through the clique construction.

Let $\Oca_{\FF_4}$ be the suboperad of $\Cli\Z$ generated by
\begin{equation}
    \left\{
        \Triangle{-1}{-1}{1}, \Triangle{-1}{1}{-1},
        \Triangle{-1}{1}{1}, \Triangle{1}{-1}{-1}
    \right\}.
\end{equation}
By using Proposition~\ref{prop:suboperads_NC_M_triangles_dimensions}, we find that the
Hilbert series $\Hilbert_{\Oca_{\FF_4}}(t)$ of $\Oca_{\FF_4}$ satisfies
\begin{equation} \label{equ:Hilbert_series_FF4}
    t + (2t - 1) \Hilbert_{\Oca_{\FF_4}}(t)
    + 2 \Hilbert_{\Oca_{\FF_4}}(t)^2 = 0.
\end{equation}
The first dimensions of $\Oca_{\FF_4}$ are
\begin{equation}
    1, 4, 24, 176, 1440, 12608, 115584, 1095424,
\end{equation}
and form Sequence~\OEIS{A156017} of~\cite{Slo}. Moreover, by computer exploration, we obtain
the list
\vspace{-2ex}
\begin{multicols}{2}
\begin{subequations}
\begin{equation} \label{equ:relation_FF4_1}
    \Triangle{-1}{1}{-1} \circ_1 \Triangle{-1}{-1}{1}
    =
    \Triangle{-1}{-1}{1} \circ_2 \Triangle{-1}{1}{-1},
\end{equation}
\begin{equation}
    \Triangle{-1}{1}{1} \circ_1 \Triangle{-1}{1}{1}
    =
    \Triangle{-1}{1}{1} \circ_2 \Triangle{-1}{1}{1},
\end{equation}
\begin{equation}
    \Triangle{-1}{1}{1} \circ_1 \Triangle{-1}{-1}{1}
    =
    \Triangle{-1}{-1}{1} \circ_2 \Triangle{-1}{1}{1},
\end{equation}
\begin{equation}
    \Triangle{-1}{1}{1} \circ_1 \Triangle{-1}{1}{-1}
    =
    \Triangle{-1}{1}{1} \circ_2 \Triangle{-1}{-1}{1},
\end{equation}

\begin{equation}
    \Triangle{-1}{1}{-1} \circ_1 \Triangle{-1}{1}{1}
    =
    \Triangle{-1}{1}{1} \circ_2 \Triangle{-1}{1}{-1},
\end{equation}
\begin{equation}
    \Triangle{1}{-1}{-1} \circ_1 \Triangle{1}{-1}{-1}
    =
    \Triangle{1}{-1}{-1} \circ_2 \Triangle{1}{-1}{-1},
\end{equation}
\begin{equation}
    \Triangle{-1}{1}{-1} \circ_1 \Triangle{-1}{1}{-1}
    =
    \Triangle{-1}{1}{-1} \circ_2 \Triangle{1}{-1}{-1},
\end{equation}
\begin{equation} \label{equ:relation_FF4_8}
    \Triangle{-1}{-1}{1} \circ_1 \Triangle{1}{-1}{-1}
    =
    \Triangle{-1}{-1}{1} \circ_2 \Triangle{-1}{-1}{1},
\end{equation}
\end{subequations}
\end{multicols}
\noindent of all nontrivial relations of degree $2$ between the generators
of~$\Oca_{\FF_4}$.

\begin{Proposition} \label{prop:construction_FF4}
    The operad $\Oca_{\FF_4}$ is isomorphic to the operad $\FF_4$.
\end{Proposition}
\begin{proof}
    Let $\phi_{\FF_4} : \Oca_{\FF_4}(2) \to \FF_4(2)$ be the linear map satisfying
    \vspace{-2ex}
    \begin{multicols}{2}
    \begin{subequations}
    \begin{equation}
        \phi_{\FF_4}\left(\Triangle{-1}{-1}{1}\right)
        = \frac{1}{\{1\} \{1, 2\}},
    \end{equation}
    \begin{equation}
        \phi_{\FF_4}\left(\Triangle{-1}{1}{-1}\right)
        = \frac{1}{\{2\} \{1, 2\}},
    \end{equation}

    \begin{equation}
        \phi_{\FF_4}\left(\Triangle{-1}{1}{1}\right)
        = \frac{1}{\{1, 2\}},
    \end{equation}
    \begin{equation}
        \phi_{\FF_4}\left(\Triangle{1}{-1}{-1}\right)
        = \frac{1}{\{1\} \{2\}}.
    \end{equation}
    \end{subequations}
    \end{multicols}
    \noindent In~\cite{CHN16}, a presentation of $\FF_4$ is described wherein its generators
    satisfy eight nontrivial relations of degree $2$. These relations can be obtained by
    replacing each $\Z$-clique appearing
    in~\eqref{equ:relation_FF4_1}--\eqref{equ:relation_FF4_8} by its image by
    $\phi_{\FF_4}$. For this reason, $\phi_{\FF_4}$ uniquely extends to an operad morphism.
    Moreover, because the image of $\phi_{\FF_4}$ contains all the generators of $\FF_4$,
    this morphism is surjective. Finally, again by~\cite{CHN16}, the Hilbert series of
    $\FF_4$ satisfies~\eqref{equ:Hilbert_series_FF4}, so that $\Oca_{\FF_4}$ and $\FF_4$
    have the same dimensions. Therefore, $\phi_{\FF_4}$ is an operad isomorphism.
\end{proof}

Hence, Proposition~\ref{prop:construction_FF4} shows that the operad $\FF_4$ can be built
through the construction $\Cli$. Observe also that, as a consequence of
Proposition~\ref{prop:construction_FF4}, all suboperads of $\FF_4$ defined in~\cite{CHN16}
that are generated by a subset of~\eqref{equ:generators_FF4} can be constructed by the
clique construction.

\subsection[Bicolored noncrossing configurations]
    {Operad of bicolored noncrossing configurations} \label{subsec:operad_bnc}
The \Def{operad of bicolored noncrossing configurations} $\BNC$ is an operad defined
in~\cite{CG14}. Let us describe this operad.

A \Def{bicolored noncrossing configuration} is a noncrossing configuration $\Cfr$ where each
labeled arc is either \Def{thick} (drawn as a thick line) of \Def{dotted} (drawn as a dotted
line) and such that all dotted arcs are diagonals. For instance,
\begin{equation}
    \begin{tikzpicture}[scale=.8,Centering]
        \node[CliquePoint](0)at(-0.3,-0.95){};
        \node[CliquePoint](1)at(-0.8,-0.58){};
        \node[CliquePoint](2)at(-1.,-0.){};
        \node[CliquePoint](3)at(-0.8,0.59){};
        \node[CliquePoint](4)at(-0.3,0.96){};
        \node[CliquePoint](5)at(0.31,0.96){};
        \node[CliquePoint](6)at(0.81,0.59){};
        \node[CliquePoint](7)at(1.,0.01){};
        \node[CliquePoint](8)at(0.81,-0.58){};
        \node[CliquePoint](9)at(0.31,-0.95){};
        \draw[CliqueEdgeGray](0)--(1);
        \draw[CliqueEdgeGray](1)--(2);
        \draw[CliqueEdgeGray](2)--(3);
        \draw[CliqueEdgeGray](3)--(4);
        \draw[CliqueEdgeGray](4)--(5);
        \draw[CliqueEdgeGray](5)--(6);
        \draw[CliqueEdgeGray](6)--(7);
        \draw[CliqueEdgeGray](7)--(8);
        \draw[CliqueEdgeGray](8)--(9);
        \draw[CliqueEdgeGray](9)--(0);
        \draw[CliqueEdgeBlue](0)--(1);
        \draw[CliqueEdgeBlue](1)--(7);
        \draw[CliqueEdgeBlue](3)--(5);
        \draw[CliqueEdgeBlue](6)--(7);
        \draw[CliqueEdgeBlue](8)--(9);
        \draw[CliqueEdgeRed](1)--(5);
        \draw[CliqueEdgeRed](1)--(9);
    \end{tikzpicture}
\end{equation}
is a bicolored noncrossing configuration of size~$9$. For $n \geq 2$, $\BNC(n)$ is the
linear span of all bicolored noncrossing configurations of size $n$. Moreover, $\BNC(1)$ is
the linear span of the singleton containing the only polygon of size $1$ where its only arc
is unlabeled. The partial composition of $\BNC$ is defined graphically as follows. For
bicolored noncrossing configurations $\Cfr$ and $\Dfr$ of respective arities $n$ and $m$,
and $i \in [n]$, the bicolored noncrossing configuration $\Cfr \circ_i \Dfr$ is obtained by
gluing the base of $\Dfr$ onto the $i$th edge of $\Cfr$, and then,
\begin{enumerate}[label=(\alph*)]
    \item if the base of $\Dfr$ and the $i$th edge of $\Cfr$ are both unlabeled, the arc
    $(i, i + m)$ of $\Cfr \circ_i \Dfr$ becomes dotted;
    \item if the base of $\Dfr$ and the $i$th edge of $\Cfr$ are both thick, the arc $(i, i
    + m)$ of $\Cfr \circ_i \Dfr$ becomes thick;
    \item otherwise, the arc $(i, i + m)$ of $\Cfr \circ_i \Dfr$ is unlabeled.
\end{enumerate}
For example,
\begin{subequations}
\begin{equation}
    \begin{tikzpicture}[scale=.7,Centering]
        \node[CliquePoint](0)at(-0.49,-0.86){};
        \node[CliquePoint](1)at(-1.,-0.){};
        \node[CliquePoint](2)at(-0.5,0.87){};
        \node[CliquePoint](3)at(0.5,0.87){};
        \node[CliquePoint](4)at(1.,0.01){};
        \node[CliquePoint](5)at(0.51,-0.86){};
        \draw[CliqueEdgeGray](0)--(1);
        \draw[CliqueEdgeGray](1)--(2);
        \draw[CliqueEdgeGray](2)--(3);
        \draw[CliqueEdgeGray](3)--(4);
        \draw[CliqueEdgeGray](4)--(5);
        \draw[CliqueEdgeGray](5)--(0);
        \draw[CliqueEdgeBlue](0)--(1);
        \draw[CliqueEdgeBlue](1)--(2);
        \draw[CliqueEdgeBlue](4)--(5);
        \draw[CliqueEdgeBlue](0)--(5);
        \draw[CliqueEdgeRed](3)--(5);
    \end{tikzpicture}
    \enspace \circ_3 \enspace
    \begin{tikzpicture}[scale=.7,Centering]
        \node[CliquePoint](0)at(-0.58,-0.8){};
        \node[CliquePoint](1)at(-0.95,0.31){};
        \node[CliquePoint](2)at(-0.,1.){};
        \node[CliquePoint](3)at(0.96,0.31){};
        \node[CliquePoint](4)at(0.59,-0.8){};
        \draw[CliqueEdgeGray](0)--(1);
        \draw[CliqueEdgeGray](1)--(2);
        \draw[CliqueEdgeGray](2)--(3);
        \draw[CliqueEdgeGray](3)--(4);
        \draw[CliqueEdgeGray](4)--(0);
        \draw[CliqueEdgeBlue](0)--(1);
        \draw[CliqueEdgeBlue](1)--(4);
        \draw[CliqueEdgeBlue](3)--(4);
        \draw[CliqueEdgeRed](2)--(4);
    \end{tikzpicture}
    \enspace = \enspace
    \begin{tikzpicture}[scale=.7,Centering]
        \node[CliquePoint](0)at(-0.34,-0.93){};
        \node[CliquePoint](1)at(-0.86,-0.5){};
        \node[CliquePoint](2)at(-0.98,0.18){};
        \node[CliquePoint](3)at(-0.64,0.77){};
        \node[CliquePoint](4)at(-0.,1.){};
        \node[CliquePoint](5)at(0.65,0.77){};
        \node[CliquePoint](6)at(0.99,0.18){};
        \node[CliquePoint](7)at(0.87,-0.49){};
        \node[CliquePoint](8)at(0.35,-0.93){};
        \draw[CliqueEdgeGray](0)--(1);
        \draw[CliqueEdgeGray](1)--(2);
        \draw[CliqueEdgeGray](2)--(3);
        \draw[CliqueEdgeGray](3)--(4);
        \draw[CliqueEdgeGray](4)--(5);
        \draw[CliqueEdgeGray](5)--(6);
        \draw[CliqueEdgeGray](6)--(7);
        \draw[CliqueEdgeGray](7)--(8);
        \draw[CliqueEdgeGray](8)--(0);
        \draw[CliqueEdgeBlue](0)--(1);
        \draw[CliqueEdgeBlue](1)--(2);
        \draw[CliqueEdgeBlue](2)--(3);
        \draw[CliqueEdgeBlue](3)--(6);
        \draw[CliqueEdgeBlue](5)--(6);
        \draw[CliqueEdgeBlue](7)--(8);
        \draw[CliqueEdgeBlue](0)--(8);
        \draw[CliqueEdgeRed](2)--(6);
        \draw[CliqueEdgeRed](4)--(6);
        \draw[CliqueEdgeRed](6)--(8);
    \end{tikzpicture}\,,
\end{equation}
\begin{equation}
    \begin{tikzpicture}[scale=.7,Centering]
        \node[CliquePoint](0)at(-0.49,-0.86){};
        \node[CliquePoint](1)at(-1.,-0.){};
        \node[CliquePoint](2)at(-0.5,0.87){};
        \node[CliquePoint](3)at(0.5,0.87){};
        \node[CliquePoint](4)at(1.,0.01){};
        \node[CliquePoint](5)at(0.51,-0.86){};
        \draw[CliqueEdgeGray](0)--(1);
        \draw[CliqueEdgeGray](1)--(2);
        \draw[CliqueEdgeGray](2)--(3);
        \draw[CliqueEdgeGray](3)--(4);
        \draw[CliqueEdgeGray](4)--(5);
        \draw[CliqueEdgeGray](5)--(0);
        \draw[CliqueEdgeBlue](1)--(2);
        \draw[CliqueEdgeBlue](1)--(5);
        \draw[CliqueEdgeBlue](2)--(3);
        \draw[CliqueEdgeBlue](3)--(4);
        \draw[CliqueEdgeBlue](4)--(5);
        \draw[CliqueEdgeRed](1)--(4);
    \end{tikzpicture}
    \enspace \circ_5 \enspace
    \begin{tikzpicture}[scale=.7,Centering]
        \node[CliquePoint](0)at(-0.58,-0.8){};
        \node[CliquePoint](1)at(-0.95,0.31){};
        \node[CliquePoint](2)at(-0.,1.){};
        \node[CliquePoint](3)at(0.96,0.31){};
        \node[CliquePoint](4)at(0.59,-0.8){};
        \draw[CliqueEdgeGray](0)--(1);
        \draw[CliqueEdgeGray](1)--(2);
        \draw[CliqueEdgeGray](2)--(3);
        \draw[CliqueEdgeGray](3)--(4);
        \draw[CliqueEdgeGray](4)--(0);
        \draw[CliqueEdgeBlue](0)--(1);
        \draw[CliqueEdgeBlue](1)--(3);
        \draw[CliqueEdgeBlue](1)--(4);
        \draw[CliqueEdgeBlue](2)--(3);
        \draw[CliqueEdgeBlue](3)--(4);
        \draw[CliqueEdgeBlue](0)--(4);
    \end{tikzpicture}
    \enspace = \enspace
    \begin{tikzpicture}[scale=.7,Centering]
        \node[CliquePoint](0)at(-0.34,-0.93){};
        \node[CliquePoint](1)at(-0.86,-0.5){};
        \node[CliquePoint](2)at(-0.98,0.18){};
        \node[CliquePoint](3)at(-0.64,0.77){};
        \node[CliquePoint](4)at(-0.,1.){};
        \node[CliquePoint](5)at(0.65,0.77){};
        \node[CliquePoint](6)at(0.99,0.18){};
        \node[CliquePoint](7)at(0.87,-0.49){};
        \node[CliquePoint](8)at(0.35,-0.93){};
        \draw[CliqueEdgeGray](0)--(1);
        \draw[CliqueEdgeGray](1)--(2);
        \draw[CliqueEdgeGray](2)--(3);
        \draw[CliqueEdgeGray](3)--(4);
        \draw[CliqueEdgeGray](4)--(5);
        \draw[CliqueEdgeGray](5)--(6);
        \draw[CliqueEdgeGray](6)--(7);
        \draw[CliqueEdgeGray](7)--(8);
        \draw[CliqueEdgeGray](8)--(0);
        \draw[CliqueEdgeBlue](1)--(2);
        \draw[CliqueEdgeBlue](1)--(8);
        \draw[CliqueEdgeBlue](2)--(3);
        \draw[CliqueEdgeBlue](3)--(4);
        \draw[CliqueEdgeBlue](4)--(5);
        \draw[CliqueEdgeBlue](4)--(8);
        \draw[CliqueEdgeBlue](5)--(7);
        \draw[CliqueEdgeBlue](5)--(8);
        \draw[CliqueEdgeBlue](6)--(7);
        \draw[CliqueEdgeBlue](7)--(8);
        \draw[CliqueEdgeRed](1)--(4);
    \end{tikzpicture}\,,
\end{equation}
\begin{equation}
    \begin{tikzpicture}[scale=.7,Centering]
        \node[CliquePoint](0)at(-0.50,-0.87){};
        \node[CliquePoint](1)at(-1.00,-0.00){};
        \node[CliquePoint](2)at(-0.50,0.87){};
        \node[CliquePoint](3)at(0.50,0.87){};
        \node[CliquePoint](4)at(1.00,0.00){};
        \node[CliquePoint](5)at(0.50,-0.87){};
        \draw[CliqueEdgeGray](0)--(1);
        \draw[CliqueEdgeGray](1)--(2);
        \draw[CliqueEdgeGray](2)--(3);
        \draw[CliqueEdgeGray](3)--(4);
        \draw[CliqueEdgeGray](4)--(5);
        \draw[CliqueEdgeGray](5)--(0);
        \draw[CliqueEdgeBlue](0)--(1);
        \draw[CliqueEdgeBlue](1)--(2);
        \draw[CliqueEdgeBlue](2)--(5);
        \draw[CliqueEdgeBlue](3)--(5);
        \draw[CliqueEdgeBlue](4)--(5);
        \draw[CliqueEdgeRed](0)--(2);
    \end{tikzpicture}
    \enspace \circ_3 \enspace
    \begin{tikzpicture}[scale=.7,Centering]
        \node[CliquePoint](0)at(-0.59,-0.81){};
        \node[CliquePoint](1)at(-0.95,0.31){};
        \node[CliquePoint](2)at(-0.00,1.00){};
        \node[CliquePoint](3)at(0.95,0.31){};
        \node[CliquePoint](4)at(0.59,-0.81){};
        \draw[CliqueEdgeGray](0)--(1);
        \draw[CliqueEdgeGray](1)--(2);
        \draw[CliqueEdgeGray](2)--(3);
        \draw[CliqueEdgeGray](3)--(4);
        \draw[CliqueEdgeGray](4)--(0);
        \draw[CliqueEdgeBlue](0)--(4);
        \draw[CliqueEdgeBlue](0)--(1);
        \draw[CliqueEdgeBlue](1)--(2);
        \draw[CliqueEdgeBlue](2)--(3);
        \draw[CliqueEdgeBlue](3)--(4);
        \draw[CliqueEdgeRed](1)--(3);
        \draw[CliqueEdgeRed](1)--(4);
    \end{tikzpicture}
    \enspace = \enspace
    \begin{tikzpicture}[scale=.7,Centering]
        \node[CliquePoint](0)at(-0.34,-0.94){};
        \node[CliquePoint](1)at(-0.87,-0.50){};
        \node[CliquePoint](2)at(-0.98,0.17){};
        \node[CliquePoint](3)at(-0.64,0.77){};
        \node[CliquePoint](4)at(-0.00,1.00){};
        \node[CliquePoint](5)at(0.64,0.77){};
        \node[CliquePoint](6)at(0.98,0.17){};
        \node[CliquePoint](7)at(0.87,-0.50){};
        \node[CliquePoint](8)at(0.34,-0.94){};
        \draw[CliqueEdgeGray](0)--(1);
        \draw[CliqueEdgeGray](1)--(2);
        \draw[CliqueEdgeGray](2)--(3);
        \draw[CliqueEdgeGray](3)--(4);
        \draw[CliqueEdgeGray](4)--(5);
        \draw[CliqueEdgeGray](5)--(6);
        \draw[CliqueEdgeGray](6)--(7);
        \draw[CliqueEdgeGray](7)--(8);
        \draw[CliqueEdgeGray](8)--(0);
        \draw[CliqueEdgeBlue](0)--(1);
        \draw[CliqueEdgeBlue](1)--(2);
        \draw[CliqueEdgeBlue](2)--(3);
        \draw[CliqueEdgeBlue](2)--(8);
        \draw[CliqueEdgeBlue](3)--(4);
        \draw[CliqueEdgeBlue](4)--(5);
        \draw[CliqueEdgeBlue](5)--(6);
        \draw[CliqueEdgeBlue](6)--(8);
        \draw[CliqueEdgeBlue](7)--(8);
        \draw[CliqueEdgeRed](0)--(2);
        \draw[CliqueEdgeRed](3)--(5);
        \draw[CliqueEdgeRed](3)--(6);
    \end{tikzpicture}\,.
\end{equation}
\end{subequations}

We now consider the unitary magma $\Mca_{\BNC} := \{\Unit, \Att, \Btt\}$ wherein operation
$\Op$ is defined by the Cayley table
\begin{small}
\begin{equation}
    \begin{tabular}{|c||c|c|c|} \hline
        $\Op$ & \; $\Unit$ \; & \; $\Att$ \; & \; $\Btt$ \; \\ \hline \hline
        $\Unit$ & $\Unit$ & $\Att$ & $\Btt$ \\ \hline
        $\Att$ & $\Att$ & $\Att$ & $\Unit$ \\ \hline
        $\Btt$ & $\Btt$ & $\Unit$ & $\Btt$ \\ \hline
    \end{tabular}\, .
\end{equation}
\end{small}
In other words, $\Mca_{\BNC}$ is the unitary magma wherein $\Att$ and $\Btt$ are idempotent,
and $\Att \Op \Btt = \Unit = \Btt \Op \Att$. Observe that $\Mca_{\BNC}$ is a commutative
unitary magma, but, since
\begin{equation}
    (\Btt \Op \Att) \Op \Att = \Unit \Op \Att = \Att
    \ne
    \Btt = \Btt \Op \Unit = \Btt \Op (\Att \Op \Att),
\end{equation}
the operation $\Op$ is not associative.

Let $\phi_{\BNC} : \BNC \to \NC\Mca_{\BNC}$ be the linear map defined in the following way.
For a bicolored noncrossing configuration $\Cfr$, $\phi_{\BNC}(\Cfr)$ is the noncrossing
$\Mca_{\BNC}$-clique of $\NC\Mca_{\BNC}$ obtained by replacing all thick arcs of $\Cfr$ by
arcs labeled by $\Att$, all dotted diagonals of $\Cfr$ by diagonals labeled by $\Btt$, all
unlabeled edges and bases of $\Cfr$ by edges labeled by $\Btt$, and all unlabeled diagonals
of $\Cfr$ by diagonals labeled by $\Unit$. For instance,
\begin{equation}
    \phi_{\BNC}\left(
    \begin{tikzpicture}[scale=.7,Centering]
        \node[CliquePoint](0)at(-0.49,-0.86){};
        \node[CliquePoint](1)at(-1.,-0.){};
        \node[CliquePoint](2)at(-0.5,0.87){};
        \node[CliquePoint](3)at(0.5,0.87){};
        \node[CliquePoint](4)at(1.,0.01){};
        \node[CliquePoint](5)at(0.51,-0.86){};
        \draw[CliqueEdgeGray](0)--(1);
        \draw[CliqueEdgeGray](1)--(2);
        \draw[CliqueEdgeGray](2)--(3);
        \draw[CliqueEdgeGray](3)--(4);
        \draw[CliqueEdgeGray](4)--(5);
        \draw[CliqueEdgeGray](5)--(0);
        \draw[CliqueEdgeBlue](1)--(2);
        \draw[CliqueEdgeBlue](1)--(3);
        \draw[CliqueEdgeBlue](0)--(5);
        \draw[CliqueEdgeRed](1)--(4);
    \end{tikzpicture}
    \right)
    \enspace = \enspace
    \begin{tikzpicture}[scale=.7,Centering]
        \node[CliquePoint](0)at(-0.49,-0.86){};
        \node[CliquePoint](1)at(-1.,-0.){};
        \node[CliquePoint](2)at(-0.5,0.87){};
        \node[CliquePoint](3)at(0.5,0.87){};
        \node[CliquePoint](4)at(1.,0.01){};
        \node[CliquePoint](5)at(0.51,-0.86){};
        \draw[CliqueEdge](0)edge[]node[CliqueLabel]{$\Btt$}(1);
        \draw[CliqueEdge](1)edge[]node[CliqueLabel]{$\Att$}(2);
        \draw[CliqueEdge](2)edge[]node[CliqueLabel]{$\Btt$}(3);
        \draw[CliqueEdge](3)edge[]node[CliqueLabel]{$\Btt$}(4);
        \draw[CliqueEdge](4)edge[]node[CliqueLabel]{$\Btt$}(5);
        \draw[CliqueEdge](0)edge[]node[CliqueLabel]{$\Att$}(5);
        \draw[CliqueEdge](1)edge[bend right=30]node[CliqueLabel]{$\Btt$}(4);
        \draw[CliqueEdge](1)edge[bend right=30]node[CliqueLabel]{$\Att$}(3);
    \end{tikzpicture}\,.
\end{equation}

\begin{Proposition} \label{prop:construction_BNC}
    The linear span of $\UnitClique$ together with all noncrossing $\Mca_{\BNC}$-cliques
    without edges nor bases labeled by $\Unit$ forms a suboperad of $\NC\Mca_{\BNC}$
    isomorphic to~$\BNC$. Moreover, $\phi_{\BNC}$ is an isomorphism between these two
    operads.
\end{Proposition}
\begin{proof}
    We denote the subspace of $\NC\Mca_{\BNC}$ described in the statement of the proposition
    by $\Oca_\BNC$. First of all, it follows from the definition of the partial composition
    of $\NC\Mca_{\BNC}$ that $\Oca_\BNC$ is closed under the partial composition operation.
    Hence, and since $\Oca_\BNC$ contains the unit of $\NC\Mca_{\BNC}$, $\Oca_\BNC$ is an
    operad. Second, observe that the image of $\phi_{\BNC}$ is the underlying space of
    $\Oca_\BNC$ and, from the definition of the partial composition of $\BNC$, one can check
    that $\phi_{\BNC}$ is an operad morphism. Finally, since $\phi_{\BNC}$ is a bijection
    from $\BNC$ to $\Oca_\BNC$, the statement of the proposition follows.
\end{proof}

Hence, Proposition~\ref{prop:construction_BNC} shows that the operad $\BNC$ can be
built through the noncrossing clique construction. Moreover, observe that in~\cite{CG14}, an
automorphism of $\BNC$ called \Def{complement} is considered. The complement of a bicolored
noncrossing configuration is an involution acting by modifying the labels of some of its
arcs. Under our setting, this automorphism translates simply into the map $\Cli\theta :
\Oca_\BNC \to \Oca_\BNC$ where $\Oca_\BNC$ is the operad isomorphic to $\BNC$ described in
the statement of Proposition~\ref{prop:construction_BNC} and $\theta : \Mca_{\BNC} \to
\Mca_{\BNC}$ is the unitary magma automorphism of $\Mca_{\BNC}$ satisfying $\theta(\Unit) =
\Unit$, $\theta(\Att) = \Btt$, and~$\theta(\Btt) = \Att$.

Moreover, it is shown in~\cite{CG14} that the set of bicolored noncrossing configurations of
arity $2$ is a minimal generating set of $\BNC$. Thus, by
Proposition~\ref{prop:construction_BNC}, the set
\begin{equation}
    \left\{
        \Triangle{\Att}{\Att}{\Att},
        \Triangle{\Att}{\Att}{\Btt},
        \Triangle{\Att}{\Btt}{\Att},
        \Triangle{\Att}{\Btt}{\Btt},
        \Triangle{\Btt}{\Att}{\Att},
        \Triangle{\Btt}{\Att}{\Btt},
        \Triangle{\Btt}{\Btt}{\Att},
        \Triangle{\Btt}{\Btt}{\Btt}
    \right\}
\end{equation}
is a minimal generating set of the suboperad $\Oca_\BNC$ of $\NC\Mca_{\BNC}$ isomorphic to
$\BNC$. As a consequence, all the suboperads of $\BNC$ defined in~\cite{CG14} which are
generated by a subset of the set of generators of $\BNC$ can be constructed by the
noncrossing clique construction. This includes, among others, the magmatic operad, the free
operad on two binary generators, the operad of noncrossing plants~\cite{Cha07}, the
dipterous operad~\cite{LR03,Zin12}, and the $2$-associative operad~\cite{LR06,Zin12}.

\section{Conclusion and perspectives}
In this article we have completed the study of the clique construction introduced
in~\cite{Cliques1} by focusing on the suboperad $\NC\Mca$ of $\Cli\Mca$ of noncrossing
$\Mca$-cliques. As noticed in the previous sections, $\NC\Mca$ has a particular status among
the suboperads of $\Cli\Mca$ because $\NC\Mca$ is the smallest suboperad of $\Cli\Mca$ that
contains all elements of arity $2$ (the $\Mca$-triangles) and is the biggest binary
suboperad of $\Cli\Mca$. This operad is also a Koszul operad when $\Mca$ is a finite unitary
magma.

An open question concerns the Koszul dual $\NC\Mca^!$ of $\NC\Mca$.
Section~\ref{sec:dual_NC_M} contains results about this operad, such as a description of its
presentation and a formula for its dimensions. We have also established the fact that, as
graded vector space, $\NC\Mca^!$ is isomorphic to the linear span of all noncrossing dual
$\Mca$-cliques. To obtain a realization of $\NC\Mca^!$, it is now enough to endow this last
space with an adequate partial composition. Finding such a composition is worth to obtain.

\MakeReferences

\end{document}